\documentclass[11pt]{amsart}
\usepackage{amssymb ,amsfonts}
\usepackage{xypic}
\usepackage{amsthm}
\usepackage[all]{xy}
\usepackage{amscd}
\usepackage{amsthm}

\newtheorem{thm}{Theorem}[section]
\newtheorem{cor}[thm]{Corollary}
\newtheorem{pro}[thm]{Proposition}
\newtheorem{lem}[thm]{Lemma}

\theoremstyle{definition}
\newtheorem{defn}[thm]{Definition}

\newtheorem{exmp}[thm]{Example}

\newtheorem{con}[thm]{Construction}

\newtheorem{rem}[thm]{Remark}

\newcommand{\p}{\text{pro--}}

\newcommand{\mm}{\mathcal{M}}

\newcommand{\ee}{\mathcal{A}}
\newcommand{\pp}{\mathcal{P}}

\newcommand{\bb}{\mathcal{B}}

\newcommand{\ff}{F_{\infty}}

\newcommand{\oo}{\mathcal{O}}
\newcommand{\dd}{\mathcal{D}}
\newcommand{\E}{\mathcal{E}}
\newcommand{\cc}{\mathcal{C}}
\newcommand{\ww}{\mathcal{W}}
\newcommand{\jj}{\mathcal{J}}
\newcommand{\rarr}{\rightarrow}
\newcommand{\ra}{\rightarrow}
\newcommand{\uu}{\mathcal{U}}
\newcommand{\uuu}{\mathcal{X}}
\newcommand{\vv}{\mathcal{V}}
\newcommand{\kk}{\mathcal{K}}
\newcommand{\ttt}{\mathcal{T}}
\newcommand{\sss}{\mathcal{S}}

\newcommand{\op}{\text{op}}
\newcommand{\holim}{\text{holim}}

\newcommand{\colim}{\text{colim}}
\newcommand{\hocolim}{\text{hocolim}}
\newcommand{\lie}{\text{Lie}}
\newcommand{\disc}{\text{dsc}}

\newcommand{\zz}{\mathbb{Z}}
\newcommand{\ess}{\text{essentially levelwise}}

\newcommand{\rr}{\mathbb{R}}

\newcommand{\fin}{\text{fnt}}
\newcommand{\free}{\text{free}\, }
\newcommand{\cofree}{\text{cofree}\, }
\newcommand{\Ho}{\text{Ho}}

\newcommand{\good}{\text{Illman}}
\newcommand{\fib}{\text{hofib}}

\newcommand{\sing}{\text{sing} \, }
\newcommand{\coef}{\mathcal{G}}
\newcommand{\comp}{\mathcal{F}_G}
\newcommand{\stable}{\text{st}}
\renewcommand{\lim}{\text{lim}}

\newcommand{\dfn}{\textbf}
\newcommand{\mdfn}[1]{#1}
\newcommand{\map}{\text{Map}}
\newcommand{\col}{\colon\thinspace}
\newcommand{\cont}{\text{cont}}
\newcommand{\midvee}{\textstyle\bigvee}

\makeatletter
\let\c@equation\c@thm
\makeatother \numberwithin{equation}{section}

\title{Equivariant homotopy theory for pro--spectra}
\author{ Halvard Fausk }  \email{fausk@math.uio.no}
\address{Department of  Mathematics,  University of Oslo, 1053
Blindern, 0316 Oslo, Norway } \subjclass{ Primary 55P91;
Secondary 18G55}

\begin{document}

\begin{abstract}
 We extend the theory of equivariant orthogonal
spectra from finite groups to profinite groups, and more generally
from compact Lie groups to compact Hausdorff groups.
 The $G$--homotopy theory  is   ``pieced  together'' from the
 $G/U$--homotopy theories for suitable quotient groups
  $G/U$ of $G$;  a motivation is the way
   continuous group cohomology of a profinite group is built out
of the cohomology of its finite quotient groups.
 In the model category of equivariant spectra
 Postnikov towers
  are studied  from a
general perspective. We introduce pro--$G$--spectra and construct
various model structures on them. A key property of the model
structures is that pro--spectra are
 weakly equivalent to their  Postnikov towers.
We   discuss two versions of a model structure
with ``underlying weak equivalences''. One of the versions only
makes sense for pro--spectra.
 In the end we  use  the theory to
 study  homotopy fixed points of pro--$G$--spectra.
\end{abstract}

\maketitle

\section{Introduction}

This paper is devoted to the exploration of some aspects of equivariant
homotopy theory of $G$--equivariant orthogonal spectra when $G$ is a
profinite group. We develop the theory sufficiently to be able to
construct homotopy fixed points of $G$--spectra in a natural way.
A satisfactory theory of $G$--spectra, when $G$ is a profinite
group, requires the generality of pro--$G$--spectra. The results
needed about model structures on pro--categories are presented in
two
 papers joint with Daniel Isaksen \cite{ffi}
\cite{tfi}. Most of the theory also works for compact Hausdorff
groups and discrete groups.

  We start out by considering  model
structures on $G$--spaces. This is needed as a starting point for
the model structure on $G$--spectra. A set of closed subgroups of
$G$ is said to be a collection if it is closed under conjugation. To
any collection $\cc$ of subgroups of $G$, we construct a model
structure on the category of $G$--spaces such that a $G$--map $f$ is
a weak equivalence if and only if $f^H$ is a
underlying weak equivalence for $H \in \cc$.

 The collections of subgroups of $G $ that
play the most important role in this paper are the cofamilies,
 i.e.~collections of subgroups that are closed under passing to larger
subgroups.
The example to keep in mind is the cofamily of open subgroups in a
profinite group.

We present the foundation for the theory of orthogonal $G$--spectra,
indexed on finite orthogonal $G$--representations, with minimal
assumptions on the group $G$ and the collection $\cc$. Most of the
results extend easily from the theory developed for compact Lie
groups by Michael Mandell and Peter May \cite{mma}. We include enough details to
make our presentation readable, and provide new proofs when the
generalizations to our context are not immediate. Equivariant
$K$--theory and stable equivariant cobordism theory both extend from
compact Lie groups to general compact Hausdorff groups. A
generalization of the Atiyah--Segal completion theorem is studied in
\cite{fg}.

 Let $R$ be a symmetric monoid in the
category of orthogonal $G$--spectra indexed on a universe of
$G$--representations. In Theorem \ref{spectramodel} the category of
$R$--modules, denoted $\mm_R$, is given a stable model structure;
 the weak equivalences are maps whose $H$--fixed points are
stable equivalences for all $H$ in a suitable collection $\cc$. For
example $\cc$ might be the smallest cofamily containing all normal
subgroups $H$ of $G$ such that $ G/H$ is a compact Lie group. A
stable $G$--equivariant theory of spectra, for a profinite group $G$,
is also given by Gunnar Carlsson in \cite{car}.

We would like to have a notion of ``underlying weak equivalences'' even
when the trivial subgroup is not included in the
collection $\cc$. We consider a more general framework.
 In  Theorem \ref{model2} we show   that for two reasonable
 collections, $\ww$ and $\cc$,
 of subgroups of $G$ such that $W U $ is in
$\cc$, whenever $W \in \ww$ and $U \in \cc$, there is a model
structures on $\mm_R$ such that the cofibrations are retracts of
relative $\cc$--cell complexes and the weak equivalences are maps $f$ such that
$\Pi_{*}^W (f) = \colim_{U \in \cc} \pi_{*}^{WU} (f)$ is an
isomorphism for every $W \in \ww$. For example, $\cc$ can be the
collection of open subgroups of a profinite group $G$ and $\ww$ the
collection, $\{ 1 \}$, consisting of the trivial subgroup in $G$.

In the rest of this introduction we assume that $\Pi^{U}_n (R) = 0$
whenever $n < 0$ and $U \in \cc$. We can then set up a good theory
of Postnikov sections in $\mm_R$. The Postnikov sections are used in
our construction of the model structures on $\p \mm_R$. Although we
are mostly interested in the usual Postnikov sections that cut off
the homotopy groups at the same degree for all subgroups $W \in
\ww$, we give a general construction that allows the cutoff to take
place at different degrees for different subgroups.

In Theorem \ref{promodelcategory} we construct a stable model
structure, called the Postnikov
 $\ww$--$\cc$--model structure, on $\p \mm_R$. It can be thought
 of as the localization of the strict model structure on $\p
 \mm_R$, where we invert all maps from a pro--spectrum to
  its levelwise Postnikov tower,  regarded as a pro--spectrum.
Here is one characterization of the weak equivalences: the class of
weak equivalences in the
 Postnikov $\ww$--$\cc$--model structure is   the
class of pro--maps that are isomorphic to a levelwise map $\{ f_s
\}_{s \in S}$ such that $f_s$ becomes arbitrarily highly connected
(uniformly with respect to the collection $\ww$) as $s$ increases
\cite[3.2]{ffi}.

In Theorem \ref{thm:AH-spectralsequence} we give an
Atiyah--Hirzebruch spectral sequence. It is constructed using the
Postnikov filtration of the target pro--spectrum. The spectral
sequence has good convergence properties because
 any pro--spectrum can be recovered from its Postnikov
tower in our model structure.

The category $\p \mm_R$ inherits a tensor product from $\mm_R$. This
tensor structure is not closed, and it does not give a well-defined
tensor product on the whole homotopy category of $\p \mm_R$ with the
Postnikov $\ww$--$\cc$--model structure.

The Postnikov $\ww$--$\cc$--model structure on $\p \mm_R$ is a stable
model structure. But the associated homotopy category is
not an axiomatic stable homotopy category in the sense of
Hovey--Palmieri--Strickland \cite{hps}.

We discuss two model structures on $\p \mm_R$ with two different
notions of ``underlying weak equivalences''. Let $G$ be a finite
group and let $\cc$ be the collection of all subgroups of $G$. There
are many different, but Quillen equivalent, $\ww$--$\ee$--model
structures on $\mm_R$ with $\ww = \{1 \}$ and $ \{1 \} \subset
\ee \subset \cc$. Two
extreme model structures are the cofree model structure, with $\ee =
\cc$, and the free model structure, with $\ee = \ww = \{ 1 \}$. The
cofibrant objects in the free model structure are retracts of
relative $G$--free cell spectra.

Now let $G$ be a profinite group and let $\cc$ be the collection of
all open subgroups of $G$.
In this case the situation is more complicated. The $\{1 \}$--weak
equivalences are maps $f$ such that $\Pi^{  1 }_* (f) =
\colim_{U \in \cc} \,  \pi_{*}^{U} (f)$ is an isomorphism. We call these
maps the $\cc$--underlying weak equivalences.
The Postnikov $\{ 1\}$--$\cc$--model structure on $\p
\mm_R$ is the closest we  get  a cofree model structure. It is
given in Theorem \ref{thm:Cfreemodelstructure}. Assume  $G$ is a nonfinite
profinite group.  Certainly, it is not
sensible to have a model structure with cofibrant objects relative
free $G$--cell complexes, because $S^n \wedge G_+$ is equivalent to
a point. In $\p \mm_R$, unlike $\mm_R$, we can form an arbitrarily
good approximation to the free model structure by letting the
cofibrations be retracts of levelwise relative $G$--cell complexes
that become ``eventually free''. That is, as we move up the inverse
system of spectra, the stabilizer subgroups of the relative cells
become smaller and smaller subgroups in the collection $\cc$. The
key idea is that the cofibrant replacement of the constant
pro--spectrum $\Sigma^{\infty} S^0$ should be the pro--spectrum
\[ \{  \Sigma^{\infty} EG/N_+ \},\] indexed by  the normal subgroups
 $N$ of $G$ in
$\cc$, ordered by inclusions. We use the  theory of
filtered model categories, developed in \cite{ffi}, to construct the
free model structure on $\p \mm_R$. This $\cc$--free model structure
is given in Theorem \ref{thm:freemodelstructure}.

 The $\cc$--free and $\cc$--cofree  model structures
on $\p \mm_R$ are Quillen adjoint, via the identity functors, but there
are fewer weak equivalences in the free than in the cofree model
structure. Thus, we actually get two different homotopy categories.
We relate this to the failure of having an internal hom functor in the
pro--category. Let
 $\Ho ( \p \mm_R )$ denote  the homotopy category of $\p \mm_R$ with
 the Postnikov $\cc$--model structure. Assume that $X$ is cofibrant and
 that $Y $ is fibrant in the Postnikov $ \cc$--model structure on $\p \mm_R$.
 Then Theorem
\ref{thm:free-cofree} says that the homset of maps from $X$ to $Y$
in the homotopy category of the $\cc$--free model structure on $\p
\mm_R$ is:
\[ \Ho ( \p \mm_R )\,  \left( X \wedge \{ EG/N_+ \} ,  Y \right)  \]
while the homset in the homotopy category of the $\cc$--cofree model
structure on $\p \mm_R$ is:
\[ \Ho ( \p \mm_R ) \, ( X , \hocolim_N \, F ( EG/N_+ , Y )
 , \] where the colimit is taken levelwise.

The Postnikov model structures are well-suited for studying homotopy
fixed points. For definiteness, let $G$ be a profinite group, let
$\cc$ be the collection of open subgroups of $G$, and let $R$ be a
non-equivariant $\sss $--cell spectrum with trivial homotopy groups in
negative degrees. The homotopy fixed points of a pro--$G$--spectrum
$\{ Y_t \}$ is defined to be the $G$--fixed points of a fibrant
replacement in the Postnikov $\cc$--cofree model structure. It is
equivalent, in the Postnikov model structure on $R $--spectra, to
the pro--spectrum
\[ \hocolim_N \, F ( EG/N_+ , P_n Y_t )^G
\] indexed on $n$ and $t$. The spectrum associated to the homotopy
fixed point pro--spectrum (take homotopy limits) is
equivalent to
\[\holim_{t,m } \, \hocolim_N \, F ( (EG/N)^{ (m) }_+ ,  Y_t )^G .\]
These expressions resemble the usual formula for homotopy fixed
points.

The appropriate notion of a ring spectrum in $\p \mm_R$ is a monoid
in $\p \mm_R$. This is more flexible
 than a pro--monoid.
The second formula for homotopy fixed point spectra shows that if
$Y$ is a (commutative) fibrant monoid in $\p \mm_R$ with the strict
$\cc$--model structure, then the associated homotopy fixed point
spectrum is a (commutative) monoid in $\mm_R$.

Under reasonable assumptions there is an iterated homotopy fixed
point formula. This appears to be false if one defines homotopy
fixed points in the $\cc$--strict model structure on $\p \mm_R$. We
obtain a homotopy fixed point spectral sequence as a special case of
the Atiyah--Hirzebruch spectral sequence.

The explicit formulas for the homotopy fixed points, the good
convergence properties of the homotopy fixed point spectral
sequence, and the iterated homotopy fixed point formula are all
reasons for why it is convenient to work in the Postnikov
$\cc$--model structure.

A general theory of homotopy fixed point spectra for actions by
profinite groups was first studied by Daniel Davis in his
Ph.D.~thesis \cite{dav2}. His theory was inspired by a homotopy
fixed point spectral sequence for $E_n$, with an action by the
extended Morava stabilizer group, constructed by Ethan Devinatz and
Michael Hopkins \cite{dho}. We show that our definition of homotopy
fixed point spectra agrees with Davis' when $G$ has finite virtual
cohomological dimension.
 Our theory applies to the
example of $E_n$ above, provided we follow Davis and use the
``pro--spectrum $K (n)$--localization''
of $E_n$ rather than (the $K(n)$--local spectrum)
$E_n$ itself.

\subsection{Acknowledgements} The model theoretical foundation for
this paper is joint work with Daniel Isaksen. I am grateful to
him for many discussions on the foundation of this paper.
I am also grateful  to  Andrew Blumberg and especially
Daniel Davis for helpful comments on the paper.

\section{Unstable equivariant theory}
 We associate to a collection, $\ww$, of closed
subgroups of $G$ a model structure on the category of based
$G$--spaces. The weak equivalences in this model structure are maps $f$
such that the $H$--fixed points map $f^H$ is a non-equivariant weak
equivalence for each $H \in \ww$.
\subsection{$G$--Spaces}
We work in the category of compactly generated weak Hausdorff
spaces.
 Let $G$ be a topological group. A $G$--space $X$ is
a topological space together with a continuous left action by $G$.
The \mdfn{stabilizer} of $x \in X$ is $\{ g \in G \ | \ g x = x \}$. This
is a closed subgroup of $G$ since it is the preimage of the diagonal
in $X \times X$ under the map $g \mapsto x \times g x $. Let $Z$ be
any subset of $X$. The stabilizer of $Z$  is the intersection of the
stabilizers of the points in $Z$, hence a closed subgroup of $G$.
Similarly, for any subgroup $H$ of $G$ the \mdfn{ $H$--fixed
points}, $X^H = \{ x \in X \ | \ h x = x \text{ for each } h \in H
\} $, of a $G$--space $X$ is a closed subset of $X$. The stabilizer
of $X^H$ contains $H$ and is a closed subgroup of $G$. So $X^H =
X^{\overline{H}}$, for any subgroup $H$ of $G$, where $\overline{H}$
denotes the closure of $H$ in $G$. Hence, we consider closed
subgroups of $G$ only.  The $H$--fixed point functor commutes
with  pushout along a closed inclusion.

A based $G$--space is a $G$--space together with a $G$--fixed
basepoint.
 We denote the category of based
$G$--spaces and basepoint preserving continuous $G$--maps by \mdfn{$G \ttt
$}.  The category of based
$G$--spaces $G\ttt$ is complete and cocomplete.

We denote the category of based $G$--spaces and continuous basepoint
preserving maps by \mdfn{$\ttt_G$}. The space of continuous maps is given a
$G$--action by $ (g \cdot f) (x) = g f (g^{-1} x )$ (and topologized
as the Kellyfication of the compact open topology). The action of $G$
on $\ttt_G (X , Y )$ is continuous, since the adjoint of the action map, $ G
\times \ttt_G (X ,Y ) \times X \rarr Y$,  is
continuous.  The corresponding categories of unbased $G$--spaces are denoted $G \uu$
and $\uu_G$.

The category $ G\ttt$ is a closed symmetric tensor category, where
$S^0$ is the unit object, the smash product $ X \wedge Y$ is the
tensor product, and the $G$--space $ \ttt_G (X ,Y )$ is the internal hom functor.

 Define  a functor $G \uu \rarr G \ttt$ by attaching a disjoint basepoint, $X
\mapsto X_+$. This functor is  left adjoint to the forgetful
functor $ G \ttt \rarr G \uu$.
The morphism set $ G \uu ( X ,Y )$ is naturally a retract of $G \ttt
( X_+ , Y_+ )$. More precisely, we have that
\[ G \ttt ( X_+ , Y_+ ) =
\textstyle\coprod_{ Z } G\uu ( Z , Y )  \] where the sum is over all
open and closed $G$--subsets $Z$ of $X$. Let $f \col X_+ \rarr Y_+ $
be a map in $G\ttt$. Then the corresponding unbased map is $f|Z \col
Z \rarr Y$ where $Z= X_+ - f^{-1} (+)$.

\subsection{Collections of subgroups of $G$}
This paper is  mostly concerned with cofamilies of subgroups.

\begin{defn}
\label{defn:collection} A  collection  $\ww$ of subgroups of
$G$ is a nonempty set of closed subgroups of $G$ such that if $H \in
\ww$, then $g H g^{-1} \in \ww $ for any $g \in G$. A collection
$\ww$ is a \mdfn{normal collection} if for all $ H \in \ww$ there
exists a $K \in \ww$ such that $K \leq H$ and $K$ is a normal
subgroup of $G$.
\end{defn}

\begin{defn} \label{defn:normal}
A collection $\ww$ of subgroups of $G$ is a  cofamily  if $K
\in \ww$ implies that $L \in \ww$ for all subgroups $ L \geq K$. A
collection $\cc$ of subgroups of $G$ contained in a cofamily $\ww$
is a  family in $\ww$,  if, for all $K \in \cc$ and $H \in
\ww$ such that $ H \leq K$, we have that $H \in \cc$.

 Let $\ww$ be a collection of subgroups of $G$.
The smallest cofamily of closed subgroups of $G$ containing $\ww$ is
called the {cofamily closure}  of $\ww$ and is denoted
\mdfn{$\overline{\ww}$}.
 A cofamily is called a {normal} cofamily
 if it is the cofamily closure of a collection of normal
 subgroups of $G$. \end{defn}
 We now give some important
cofamilies.

\begin{exmp} \label{exmp:fnt}
 The collection of all  subgroups $U$  of $G$ such that
$G/U $ is finite and discrete is a cofamily. This collection of
 subgroups is closed under finite intersection  since $ G/
U\cap V \leq G/ U \times G / V$. A finite index subgroup of $G$ has
only finitely many
 $G$--conjugate subgroups of $G$. Hence, if $U$ is a finite index
subgroup of $G$, then $\cap_{g \in G} g U g^{-1} $ is a normal
subgroup of $G$ such that $G/ \cap_{g \in G} g U g^{-1} $ is a
finite discrete group. Let  $\fin (G)$  be the collection of
all normal subgroups $U$ of $G$ such that $G/U$ is a finite discrete
group.
\end{exmp}
\begin{exmp} \label{exmp:discrete}
Define  $\disc (G)$  to be the collection of all normal
subgroups $U$ of $G$ such that $G/U$ is a discrete group. This
collection is closed under intersection. We call a collection that
is contained in the cofamily closure of $\disc (G)$ a  discrete
collection  of subgroups of $G$.
\end{exmp}
\begin{exmp}  \label{exmp:lie}
Let  $\lie (G)$  be the collection of all normal subgroups $U$
of $G$ such that $G/U$ is a compact Lie group. This collection is
closed under intersection since a closed subgroup of a compact Lie
group is a compact Lie group. We call a collection that is contained
in the cofamily closure of $\lie (G)$ a  Lie collection  of
subgroups of $G$.
\end{exmp}

\begin{lem} \label{egk} Let $G$ be a compact Hausdorff group, and
let $K$ be a closed subgroup of $G$. Then $ \{ U \cap K \, | \, U
\in \lie (G) \}$ is a subset of $\lie (K)$, and for every $H \in
\lie (K)$ there exists a $U \in \lie (G)$ such that $U \cap K
\subset H$.
\end{lem}
\begin{proof}
Let $U \in \lie (G)$. The subgroup $U \cap K$ is in $\lie (K)$ since
$ K/ K\cap U $ is a closed subgroup of the compact Lie group $G /U$.

Let $H$ be a subgroup in $\lie (K)$. We have that $\cap_{U \in \lie
(G)} U = 1 $ by Corollary \ref{Neil}. Hence $U \cap K / H $ for $U
\in \lie (G)$ is a collection of closed subgroups of the compact Lie
group $K / H$ whose intersection is the unit element. Since $\lie
(G)$ is closed under finite intersections, the descending chain
property for closed subgroups of a compact Lie group \cite[1.25,
Exercise 15]{tom} gives that there exists a $U \in \lie (G)$ such that $U
\cap K$ is contained in $H$.
\end{proof}

We order $\fin (G)$ and $\lie (G)$ by inclusions. We recall the
following facts.
\begin{pro}
A topological group $G$ is a profinite group precisely when
  \[ G \rarr \lim_{U \in
\fin (G)} G/U \] is a homeomorphism. A topological group $G$ is a
compact Hausdorff group precisely when
\[ G \rarr \lim_{U \in \lie (G)} G/U \]
 is a homeomorphism.
\end{pro}
\begin{proof} These  facts   are  well-known.
The second claim is   proved in Proposition \ref{continuous}.
\end{proof}

Even though we are mostly interested in actions by profinite groups,
we find it natural to study actions by compact Hausdorff groups
whenever possible.

\subsection{Model structures on the category of $G$--spaces}

We associate to a collection $\ww$ of closed subgroups of $G$ a
model structure on the category of based $G$--spaces.

\begin{defn} \label{defn:Wweakequivalencesspaces}
Let $f \col X \rarr Y$ be a map in $G \ttt$. The map $f$ is said to
be a \mdfn{$\ww$--}\mdfn{equivalence} if the underlying unbased maps
$ f^U \col X^U \rarr Y^U$ are weak equivalences for all $U \in \ww$.
\end{defn}

\begin{defn}  Let $p \col E \rarr B$ be a map in $G\ttt$.
We say that $p$ is a \mdfn{$\ww $--fibration} if the underlying
unbased maps $p^U \col E^U \rarr B^U$ are Serre fibrations for all
$U \in \ww$.
\end{defn}

We next define the generating cofibrations and generating acyclic
cofibrations. We use the conventions that $S^{-1}$
is the empty set and $D^0 $ is a point.
\begin{defn} \label{generators}
   Let \mdfn{$ \ww I$} be the set  of maps
 \[   \{ (G/U \times S^{n-1})_+ \rarr (G/U \times
D^{n})_+ \} ,\] for $n \geq 0 $ and $U \in \ww$. Let \mdfn{$ \ww J
$} be the set of maps \[ \{(G/U \times D^{n})_+ \rarr (G/U \times
D^{n}\times [0,1])_+ \} ,\] for $n \geq 0 $ and $U \in \ww$.
\end{defn}

 The following model structure is called the \mdfn{$\ww$--model
 structure
on $G \ttt$}. For the definition of relative cell complexes see
\cite[10.5]{hir}.
\begin{pro} \label{modelstructurepaces}
 There is a proper model structure on $G\ttt$ with weak
equivalences $\ww$--weak equivalences, fibrations $\ww$--fibrations,
and cofibrations retracts of relative $\ww I$--cell complexes. The set
$\ww I$ is a set of generating cofibrations and $\ww J $ is a set of
generating acyclic cofibrations.
\end{pro}
\begin{proof} A  map $p \col E \rarr B$ in $G \ttt$ is a
$\ww$--fibration if and only if it has the right lifting property
with respect to all maps in $\ww J$. A map $f$ is a $\ww$--acyclic
fibration if and only if it has the right lifting property with
respect to all maps in $\ww I$. This follows from the corresponding
non-equivariant result and by the fixed point adjunction
\cite[2.4]{hov}.
 The verifications of the model structure axioms follow
 as in  \cite[Section 2.4]{hov}.
The model structure is both left and right proper. This follows from
the corresponding non-equivariant results since pullbacks commute
with fixed points and since pushouts along closed inclusions also
commute with fixed points.
\end{proof}
An alternative way to set up the model structure on $G \ttt $ is
given in \cite[Section III.1]{mma}. Let \mdfn{${\ww} G \ttt $}, or simply
\mdfn{${\ww} \ttt $}, denote $G \ttt$ with the $\ww$--model
structure, and let \mdfn{$Ho ({\ww} G \ttt )$} denote its homotopy
category.

\begin{pro} \label{homotopyclasses}
Let $X$ be a retract of a $\ww I $--cell complex, and let $Y$ be a
$G$--space. Then the set $ Ho ({\ww} \ttt) ( X ,Y )$ is isomorphic
to the set of based $G$--homotopy classes of maps from
$X$ to $Y$.
\end{pro}
\begin{proof} All objects are fibrant and a
retract of a $\ww I $--cell complex is cofibrant in the $\ww$--model
structure. The cylinder object of (a cofibrant object) $X$ in the
$\ww$--model structures is $X \wedge [0,1]_+$.
\end{proof}
The next result has also been proved by Bill Dwyer \cite[4.1]{dwy}.
Note that a $G$--cell complex $X$ is a $\ww I$--cell complex if and
only if all its isotropy groups are in $\ww$.
\begin{cor} \label{dwyer}
Let $X$ and $Y$ be $\ww I $--cell complexes. If a map $f \col X
\rarr Y$ is a $\ww$--weak equivalence, then $f$ is a based $G$--homotopy
equivalence.
\end{cor}

To get a topological model structure on our model category we need
some assumptions on the collection $\ww$.
\begin{defn} \label{good}  Let $\E$ and $\ww$ be two collections of
subgroups of $G$.  Then $\E$ is called a  $\ww$--\good\ collection if
$(G/ U \times G/H)_+$ is a $\ww I $--cell complex for any  $U \in \ww$
and $H \in \E$.  A collection $\ww$ of subgroups of $G$ is called an
\good\ collection if $\ww$ is a $\ww$--\good\ collection.
\end{defn}
In particular, if $\E $ is a $\ww$--\good\ collection, then
$U \cap H \in \ww $, for all $ U \in \ww$ and $ H \in \E$,
because  $U \cap H $ is an isotropy group of $G/ U \times G /H$.
The collection $\{ G \} $ is a $\ww$--\good\ collection of
subgroups of $G$ for any collection $\ww$ of subgroups of $G$.

\begin{lem} \label{lem:goodsmash}
If $\ww$ is a discrete or a Lie collection of subgroups of $G$ and
$\ww$ is closed under intersection, then $\ww$ is an \good\
collection of subgroups of $G$.
\end{lem}
\begin{proof}  The statement is clear when $\ww$ is contained in
$\overline{\disc (G)}$. When $\ww$ is contained in $\overline{\lie
(G)}$, then the claim follows from a result of Illman \cite{ill}.
\end{proof}
The next lemma shows, in particular, that if $\ww$ is an \good\
collection of subgroups of $G$, then the smash product of two
$\ww I$--cell complexes is again a $ \ww I$--cell complex.
\begin{lem} \label{productcomplexes}
Let $\E$ and $\ww$ be collections of subgroups of $G$ such that
$\E$ is a $\ww$--\good\ collection. If  $ X$ is  a $\ww I $--cell
complex and $ Y$  is a $ \E I$--cell complex, then
$X \wedge Y$ is (homeomorphic to) a $\ww I$--cell complex.
\end{lem}
\begin{proof} It suffices  to show that \[ (S^{n-1 } \times S^{m-1}
\times G/ U \times G/ U' \rarr D^{n } \times D^{m} \times G/ U
\times G/ U')_+ \] is a relative $\ww I$--cell complex, for
$ U \in \ww$ and $U' \in \E$. This is so since $(G/ U \times G/ U')_+$
is a relative $\ww I$--cell complex.
\end{proof}
We follow the treatment of a topological model structure given in
\cite[Section III.1]{mma}. Note that the $G$--fixed points of the mapping
spaces in $\ttt_G$ are the mapping spaces in $ G \ttt$. Let $\mm_G$
be a category enriched in $G \ttt$. Let $G \mm $ be the $G$--fixed
category of $\mm_G$. Simplicial structures are defined in
\cite[9.1.1,9.1.5]{hir}. We modify the definition of a simplicial
structure by model theoretically enriching $\mm_G$ in the model
category $\ww \ttt$ instead of the model category of simplicial
sets.

Let  $ i \col  A \rarr X$ and $p \col  E \rarr B$
be two maps in $\mm_G$. Let \[\mm_G ( i^* , p_* ) \col \mm_G ( X ,E
) \rarr \mm_G ( A ,E ) \times_{ \mm_G ( A ,B ) } \mm_G (X , B ) \]
be the $G$--map induced by precomposing with $i$ and composing with
$p$.

\begin{defn}
Let $\mm_G$ be enriched over $G\ttt$.  A model structure on $G \mm$
is said to be \mdfn{$\E$--topological} if it is $G$--topological (see
\cite[9.1.2]{hir}) and the following holds:
\begin{enumerate}
\item  There is  a tensor functor  $X
 \square T$ and a cotensor functor  $F_{\square} (T ,X )$ in $\mm_G$, for
$X \in \mm_G$ and $T \in \ttt_G$, such  that there are natural
isomorphisms of based $G$--spaces
\[ \mm_G ( X \square T  , Y   ) \cong \ttt_G ( T , \mm_G (
 X , Y )) \cong \mm_G ( X   , F_{\square} (T ,Y) ),  \]
for $X ,Y \in \mm_G$ and $ T \in \ttt_G$.
\item  The map
$\mm_G ( i^* , p_* )$ is a $\E$--fibration in $G \ttt$ whenever
$i$ is a cofibration and $p$ is a fibration in $  G \mm$,
and if $i$ or $p$  in addition is a weak equivalence in $  G \mm$,
then $\mm_G ( i^* , p_* )$ is  a $\E$--equivalence in $G \ttt$.
\end{enumerate}
\end{defn}

\begin{rem} \label{rem:topms} The $G$--fixed
points of $\mm_G ( i^* , p_* )$ is $G \mm (i^* ,p_* )$. So if $\{ G
\} \in \E$, then a $\E$--topological model structure on $G \mm$
gives a topological model structure.
\end{rem}

If $\ww$ is an \good\ collection, then the following Lemma implies
the pushout--product axiom \cite[2.1,2.3]{ss}.
\begin{lem} \label{lem:topstructure}
Let $\E$ be  a $\ww$--\good\ collection.  Assume that
$ f \col  A \rarr B$ is in $\E I$ and $g \col X \rarr Y$ is in
$\ww I$, then $f \Box g \col (A \wedge Y) \cup_{A \wedge X}
( B \wedge X) \rarr B \wedge Y$ is a $\ww$--cofibration.
Moreover, if $f$ is in $\E J $ instead of $\E I$ or $g$ is in $\ww J$
instead of  $\ww I$, then $f \Box g$ is a $\ww-$acyclic cofibration.
\end{lem}
\begin{proof} This reduces to our  assumption on $\E$ and  $\ww$.
 See also \cite[II.1.22]{mma}.
\end{proof}

\begin{pro} \label{pro:topmodstructure}
Let $\E$ be a $\ww$--\good\ collection  of subgroups of $G$.
Then the $\ww$--model structure on $G \ttt$, from
Proposition \ref{modelstructurepaces}, is a $\E$--topological
model structure. \end{pro}
\begin{proof}  This  follows from
\cite[III.1.15--1.21]{mma} and Lemma \ref{lem:topstructure}.
\end{proof}

The $\ww$--model structure on $ G \mm $ is a topological model
structure for any collection $\ww$ by Remark \ref{rem:topms}.
\begin{rem} \label{simplicialstructure}
A based topological model category $\mm$ has a canonical based
simplicial model structure. In the topological model structure
denote the mapping space by $ \text{Map } (M , N)$, the tensor by $
M \square X$, and the cotensor by $ F_{\square} ( X , M)$. Here $X $
is a based space, and $M$ and $N$ are objects in $\mm$. The singular
simplicial set functor, $\text{sing} $, is right adjoint to the geometric
realization functor $|-|$. The corresponding based simplicial
mapping space is given by $ \sing (\text{Map } (M , N)) $. The
simplicial tensor and cotensor are $ M \square |K|$ and
$ F_{\square} ( |K| , M)$, respectively, where $K$ is a
based simplicial set and $M$ and $N$ are objects in $\mm$. We use
that $ |K \wedge L| \cong |K| \wedge |L|$.

A based simplicial structure gives rise to an unbased simplicial
structure. We get a unbased simplicial structure by forgetting the
basepoint in the based simplicial mapping space, and by adding a
disjoint basepoint to unbased simplicial sets in the definition of
the tensor and the cotensor. Hence we can apply results about
(unbased) simplicial model structures to a topological model
category.
\end{rem}

\subsection{Some change of groups results for spaces}
\label{sec:changeofgroups-spaces}

Let $\phi \col G_1 \rarr G_2$ be a continuous  group homomorphism
between compact Hausdorff groups.
Let $j \col G_2 \ttt \to G_1 \ttt$ be defined by restricting the
$G_2$--action to a $G_1$--action along $\phi$.  This functor  has a
 left adjoint given by sending $X$ to
$(G_2)_+ \wedge_{G_1} X$ and a right adjoint given by sending $ X$ to
$\ttt_{G_1} ( (G_2)_+ , X )$.

In general
these three functors do not behave well with respect to the model
structures on the categories of $G_1 $--spaces and $G_2 $--spaces.
We give some conditions that guarantee that they are  Quillen adjoint functors.
Let $\ww_1$ be a collection of subgroups of $G_1$ and  let $\ww_2$ be a
collection of subgroups of $G_2$. Let  $\phi \ww_1 $ be the smallest
collection of subgroups of $G_2$ containing (the closures of) $\phi (H)$,
for all $H \in \ww_1$. Let  $ \phi^{-1} \ww_2$ be the smallest collection
of subgroups of $G_1$ containing $ \phi^{-1} (K)$, for all $K \in \ww_2$.

\begin{lem} The functor $j \col  \ww_2 G_2 \ttt \to \ww_1 G_1 \ttt$
is a right Quillen adjoint functor if $ \phi  \ww_1 \subset \ww_2$ and a
left Quillen adjoin functor if,  in addition,  $ \phi^{-1} \ww_2 \subset \ww_1$
and $ \ww_2 \subset \overline{\lie (G_2)}$. \end{lem}

\begin{proof} The left adjoint functor $ (G_2)_+ \wedge_{G_1} X$
is a Quillen left adjoint functor  if it respects the generators
in Definition \ref{generators}.  Hence the first claim follows
since $ (G_2)_+ \wedge_{G_1} G_1 /H_+ \cong (G_2 / \overline{\phi (H)})_+ $.
Let $K $ be in $ \overline{\lie (G_2)}$. Then the restriction of
$G_2 / K$ along $\phi$ is a $ G_1$--cell complex with stabilizers
in $ \phi^{-1} (gKg^{-1})$, for $g \in G_2$ \cite{ill}.
Hence $ j$ respects cofibrations  if  $ \phi^{-1} \ww_2 \subset \ww_1$.
The $H$--fixed points of $j (X)$ is $X^{ \phi (H)}$.
Hence $j$ respects acyclic cofibrations if in addition $ \phi  \ww_1 \subset \ww_2$.
\end{proof}
Any group homomorphism between compact Hausdorff  groups is a
composite of a surjective identification homomorphism followed
by a closed inclusion of a subgroup.  So for compact Hausdorff groups
it suffices to consider  these two types of group homomorphisms.
Let $K$ be a  subgroup of $G$. Then the  forgetful
functor $G \ttt \to K \ttt$ has a left adjoint given by sending $X$ to
$G_+ \wedge_K X$ and a right adjoint given by sending $ X$ to
$\ttt_K ( G_+ , X )$.  Let $N$ be a normal subgroup of $G$.
Then the   functor $ G / N \ttt \to G\ttt$ has a
left adjoint  given by the $N$--orbit functor and a
right adjoint  given by the $N$--fixed point functor.

\begin{exmp} \label{exmp:changeofgroups} Let $K$ be a subgroup of $G$.
The forgetful functor from $ \overline{\lie (G)} \, G\ttt$ to
$\overline{\lie (K)} \, K \ttt$
 is both a left and a right Quillen adjoint functor if $K$ is
in $\overline{\lie (G)}$. It is neither a left nor a right Quillen
adjoint functor if $K$ is not in $\overline{\lie (G)}$.

Let $N$ be a normal subgroup of $G$. Then the functor from
$\overline{\lie (G/N)} \, G/N \ttt$ to $\overline{\lie (G)} \,
G\ttt$ is both a left and a right Quillen adjoint functor.
\end{exmp}

\section{Orthogonal $G$--Spectra}
 Equivariant orthogonal spectra for compact Lie groups were
 introduced by Michael Mandell and Peter May in
\cite{mma}. We generalize their theory to allow more general groups.
We develop the theory with minimal assumptions on the collection of
subgroups used to define cofibrations and weak equivalences. We follow
Chapters 2 and 3 of their work closely.

\subsection{$\jj^{\vv}_G$--spaces}
We define universes of $G$--representations.
\begin{defn} \label{def:universe}  A \mdfn{$G$--universe $\uu$} is a
 countable infinite
direct sum $\oplus^{\infty}_{i =1} \uu ' $ of a real $G$--inner
product space $\uu '$ satisfying the following: (1) the
one-dimensional trivial $G$--representation is contained in $\uu'$;
(2) $\uu$ is topologized as the union of all finite dimensional
$G$--subspaces of $\uu$ (each with the norm topology); and (3) the
$G$--action on all finite dimensional $G$--subspaces $V$ of $\uu$
factors through a compact Lie group quotient of $G$.
\end{defn}
If  $G$ is a compact Hausdorff group, then
the $G$--action on a finite dimensional $G$--representation
factors through a compact Lie group quotient of $G$ by Lemma
\ref{Lierepresentation}. This is not true in general (consider
the representation $\mathbb{Q}/ \mathbb{Z} < S^1$).
We only use the finite dimensional $G$--subspaces of $\uu$, so one
might as well assume that $\uu'$ is a union of such.

\begin{defn} Let \mdfn{$S^V$} denote the one-point
compactification of a finite dimensional $G$--representation $V$.
\end{defn}
The last assumption in Definition \ref{def:universe} is added to
guarantee that spaces like $S^V$ have  the homotopy type of a
finite $\overline{\lie (G)} I$--cell complex.

\begin{defn}
If the $G$--action on $\uu$ is trivial, then $\uu$ is called a
trivial universe. If each finite dimensional orthogonal
$G$--representations is isomorphic to a $G$--subspace of $\uu$, then
$\uu$ is called a  complete $G$--universe.
\end{defn} All compact Hausdorff groups have a complete universe.
However, it might not be possible to find a complete universes  with a
countable dimension. Traditionally, the universes have often been assumed
to have countable dimension \cite[Definition IX.2.1]{ala}.

\begin{rem}
A complete $G$--universe suffices to construct a sensible
equivariant homotopy theory for
compact Hausdorff groups with the weak equivalences determined by
the cofamily closure of $\lie (G)$.  For example, there are
Spanier--Whitehead duals of suspension spectra of finite
cell--complexes with stabilizers in $ \overline{\lie (G)}$ (see Proposition \ref{pro:MRtensormodelcat}).
Transfer maps can then be constructed as in \cite[XVII. Section 1]{ala}.
\end{rem}

We recall some definitions from \cite[Chapter II]{mma}.
\begin{defn}  Let $ \uu$ be a universe.
An  indexing representation  is a finite dimensional $G$--inner
product subspace of $\uu$. If $V$ and $W$ are two indexing
representations and $ V \subset W$, then the orthogonal complement
of $V $ in $W$ is denoted by $W - V$. The collection of all real
$G$--inner product spaces that are isomorphic to an indexing
representation in $\uu$ is denoted \mdfn{$\vv (\uu) $}.
\end{defn} When $\uu$ is understood, we write $\vv$ instead of
$\vv (\uu)$ to make the notation simpler.

\begin{defn} \label{defn:jvg}
Let \mdfn{$\jj^{\vv}_G $} be the unbased topological category with
objects $V \in \vv$ and morphisms linear isometric isomorphisms. Let
\mdfn{$G \jj^{\vv} $} denote the $G$--fixed category $( \jj^{\vv}_G )^G$.
\end{defn}
\begin{defn} A continuous $G$--functor $ X\col   \jj^{\vv}_G \rarr  \ttt_G$
is called a \mdfn{ $\jj^{\vv}_G$--space}. (The induced map on hom spaces is a
continuous unbased $G$--map.)
Denote the category of $\jj^{\vv}_G$--spaces and (enriched) natural
transformations by \mdfn{$ \jj^{\vv}_G \ttt$}. Let \mdfn{$ G
\jj^{\vv} \ttt$} denote the $G$--fixed category $ ( \jj^{\vv}_G \ttt)^G$.
\end{defn}

\begin{defn} \label{defn:smash}
Let \mdfn{$\sss^{\vv}_G$}$\col \jj^{\vv}_G \rarr \ttt_G $ be the
$\jj^{\vv}_G$--space defined by sending $V$ to the one point
compactification $S^V$ of $V$.
For simplicity we sometimes  denote $\sss^{\vv}_G$ by \mdfn{$ \sss$}.
\end{defn}

The external smash product \[ \overline{\wedge} \col \jj^{\vv}_G
\ttt \times \jj^{\vv}_G \ttt \rarr ( \jj^{\vv}_G \times \jj^{\vv}_G
) \ttt \] is defined to be $ X \overline{\wedge } Y ( V , W ) =
 X (V) \wedge  Y ( W)$ for $ X , Y \in \jj^{\vv}_G$ and $V ,W \in
 \vv$. The direct sum of finite dimensional real $G$--inner product
spaces gives  $\jj^{\vv}_G$ the structure of a symmetric tensor category.
A topological left Kan extension gives an internal smash product on
$\jj^{\vv}_G \ttt$ \cite[21.4, 21.6]{mms}. We give an explicit
description of the smash product. Let $W$ be a real $N$--dimensional
$G$--representation in $\vv (\uu)$. Choose $G$--representations
$V_n$ and $V_{n}'$ of dimension $n$ in $\vv (\uu)$ for
$n = 0,1,\ldots , N$. For example let $V_n = V'_{ n}$
be the trivial $n$--dimensional $G$--representation, $\rr^n$.
Then we have a canonical equivalence
\[ X \wedge Y (W) \cong  \vee_{n=0}^{N} \jj^{\vv}_G
( W, V_n \oplus V_{N-n}') \wedge_{\mathcal{O} (V_n)
\times \mathcal{O} (V_{N-n}')} X (V_n) \wedge Y( V_{N-n}') .
\] The internal hom from $X$ to $Y$ is the $\jj^{\vv}_G$--space
\[ V \mapsto \jj^{\vv}_G \ttt ( X(-) , Y ( V \oplus -)) \]
given by the space of continuous natural transformation of
$\jj_{G}^{\vv} \ttt$--functors. The internal smash product and the
internal hom functor give $\jj^{\vv}_G \ttt$ the structure of a closed
symmetric tensor category \cite[II.3.1,3.2]{mma}. The unit object is
the functor that sends the indexing representation $V$ to $S^0$ when
$V=0$, and to a point when $V \not= 0$. By passing to fixed points
we also get a closed symmetric tensor structure on $G \jj^{\vv}
\ttt$.

\subsection{Orthogonal $R$--modules}

 For the definition of
monoids and modules over a monoid in tensor categories see
\cite[VII.3 and 4]{mac}. The functor $\sss^{\vv}_G$ is a strong
symmetric functor.  Hence the
$\jj^{\vv}_G$--space $ \sss^{\vv}_G$ is naturally a symmetric monoid in
$G \jj^{\vv} \ttt$.  The following definition is from
\cite[II.2.6]{mma}.

\begin{defn} \label{defn:orthogonalspectrum}
An \mdfn{orthogonal $G$--spectrum} $X$ is a $\jj^{\vv}_G$--space
$X \col    \jj^{\vv}_G \rarr  \ttt_G$  together with a left
module structure over the symmetric monoid  $\sss^{\vv}_G$
in $G \jj^{\vv}_G \ttt$. Denote the category of $G$--spectra by
\mdfn{$ \jj^{\vv}_G \sss$}. Let \mdfn{$ G \jj^{\vv} \sss$} be the
$G$--fixed category $ ( \jj^{\vv}_G \sss )^G$.
 \end{defn}

The smash product and internal hom functors of orthogonal spectra are
the smash product and internal hom functors of $\sss^{\vv}_G$--modules.
 So the category of orthogonal $G$--spectra, $
\jj^{\vv}_G \sss$, is itself a closed symmetric tensor category with
$\sss^{\vv}_G$ as the unit object \cite[II.3.9]{mma}. The fixed point
category  $ G \jj^{\vv} \sss$ inherits a closed tensor structure from
$ \jj^{\vv}_G \sss$. Explicit formulas for the tensor and internal hom
functors are  obtained from the formulas after Definition \ref{defn:smash}
and \cite[II.3.9]{mma}.  The monoid $\sss^{\vv}_G $ is symmetric so
a left  $\sss^{\vv}_G$--module has a natural right module structure.

\begin{defn}
We call a monoid $R$ in $ G \jj^{\vv} \sss $ an \mdfn{algebra}. We
say that $R$ is a commutative algebra, or simply a \mdfn{ring}, if
it is a symmetric monoid in $ G \jj^{\vv} \sss $. We sometimes add:
orthogonal, $G$, and spectrum, to avoid confusion.
\end{defn}

Let $R$ be an orthogonal algebra  spectrum. We assume that the
basepoint of  each based space $R(V)$ is nondegenerated. This
assumption can be circumvented for the stable model structure when $\ww$
is an \good\ collection \cite[Theorem III.3.5, Section III.7]{mma}.

\begin{defn}  An $R$--module is a left $R$--module in the category
of orthogonal spectra. Let \mdfn{$\mm_R$} denote the category of $R$--modules.
\end{defn} The category of $R$--modules is complete and cocomplete. If
$R$ is a commutative monoid, then the category $\mm_R$ is a closed
symmetric tensor category. A monoid $T$ in the category of $R$--modules
is called an $R$--algebra. Any $R$--algebra is an $\sss$--algebra.
Let $T$ be an $R$--algebra. Then the category of $T$--modules, in the
category of $R$--modules, is equivalent to the category of $T$--modules,
in the category of $\sss $--modules, when $T$ is regarded as
an $\sss $--algebra.

We now give  pairs of adjoint functors between orthogonal $G$--spectra
and $G$--spaces.  The \mdfn{$V$--evaluation functor}
\[ \Omega^{\infty}_V \col   \jj^{\vv}_G \sss \rarr \ttt_G \]
is given by $X \mapsto X (V)$.
We abuse language and let $\Omega^{\infty}_V$ also denote the functor
$\Omega^{\infty}_V$ precomposed with the forgetful functor from $R$--modules
to orthogonal spectra.
There is a   left adjoint, denoted  \mdfn{$\Sigma^{R}_V$},
of the $V$--evaluation functor $\Omega^{\infty}_V \col   \mm_R \rarr
\ttt_G $.  The $R$--module $\Sigma^{R}_V Z$, for a $G$--space $Z$,
sends $W \in \vv (\uu)$ to
\begin{equation} \label{eq:R-suspension}
\Sigma^{R}_V Z (W) =  Z \wedge  \oo (W)_+ \wedge_{ \oo (
W - V )} R (W-V ) \end{equation} when $W \supset V$, and to a point
otherwise \cite[4.4]{mms}. This functor is called the
\mdfn{ $V$--shift desuspension spectrum functor} and is also denoted
$F_V$ and $\Sigma_{V}^{\infty}$ (when $R=\sss$) in \cite{mma}. When $V=0$ we
denote this functor by \mdfn{$ \Sigma^{\infty}_R$}. We have that
$\Sigma^{R}_V Z \cong \Sigma_{V}^{\sss} Z \wedge R$.

\subsection{Fixed point and orbit spectra}
We define fixed point and orbit spectra.  The results on adjoint
functors and change of universes from
\cite[Section V.1]{mma} extend to our setting. The change of groups
results in Subsection \ref{sec:changeofgroups-spaces} extends to
the model structures on the category of orthogonal spectra
(constructed in later sections). We do not make those results explicit.

Let $X$ be an orthogonal spectrum and let $H$ be a
subgroup of $G$. Then the quotient $X/H$ is defined to be $X/H (V) =
X(V)/H$ with structure maps \[ X/H (V) \wedge S^W \rarr X(V)/H
\wedge S^W/H \cong ( X(V) \wedge S^W )/H \rarr X( V \oplus W)/H .\]
The $H$--orbit spectrum is an $N_G H$--spectrum with trivial $H$--action.
The orbit spectrum $X /H$ restricted to the the universe $\uu^H$
is an $N_G H / H$--spectrum.

Let $H$ be a   subgroup of $G$. Let $X$ be a $\uu$--spectrum. Then
the $H$--fixed point spectrum $X^H$ indexed on $\uu^H$ is defined
by $ X^H (V) = (X(V))^H$ for $V \in \vv (\uu^H)$, and the structure
map is \[X^H (V) \wedge S^W \cong (  X(V) \wedge S^W )^H \rarr X^H (
V \oplus W) . \] This is a $N_G H /H$--spectrum.
One can also define geometric fixed point spectra as in \cite[Section V.4]{mma}.

Let $N$  be a closed normal subgroup of $G$, and let $  G/N \mm$ be
the category of $G$--spectra indexed on $\uu^N$.  The $N$--orbit and
$N$--fixed point functors are left and right adjoint functor to the
restriction functor $ G/N \mm \rarr G \mm$, respectively \cite[V.3]{mma}.

\begin{rem} \label{rem:restriction} Let $H$ be a subgroup of $G$.  There is a restriction
functor $G \mm_R \rarr H \mm_R$.
A $G$--spectrum is regarded as an $H$--spectrum indexed on
$ \vv (\uu ) |H$.  This is not a spectrum indexed on $\vv (\uu | H)$
because not all indexing $H$--representations in $\uu$ are obtained
as the restriction of an indexing  $G$--representation in $\uu$.
There is a change of indexing functor associated to the full
inclusion $\vv ( \uu) |H  \rarr \vv ( \uu|H )$.  They define the
same homotopy theory so for simplicity we use $\vv ( \uu )|H$ in
this paper. See \cite[V.2]{mma}.
\end{rem}

\subsection{Examples of orthogonal $G$--Spectra} \label{rem:subspectra}
Let $T \col \ttt_G \rarr \ttt_G $ be a continuous $G$--functor. Then
we define the corresponding $\jj^{\vv}_G $--space by $T \circ
\sss^{\vv}_G \col \vv (\uu) \rarr \ttt_G$. This $\jj^{\vv}_G $--space is
given an orthogonal $G$--spectrum structure by letting $S^W \wedge T (S^V)
\rarr  T( S^{V \oplus W}) \cong T (S^V \wedge S^W ) $ be the adjoint
of the map
\[ S^W \rarr \ttt_G ( S^V , S^W \wedge S^V) \stackrel{T}{\rarr}
\ttt_G (T( S^V), T( S^W \wedge S^V) ) \] where the first map
is a $G$--map adjoint to the identity on $S^W \wedge S^V$.

We can define a $G$--equivariant orthogonal $K$--theory spectrum for a compact
Hausdorff group $G$. If $X$ is a compact $G$--space, then $K_G (X) $
is the Grothendieck construction on the semiring of isomorphism
classes of finitely generated real bundles on $X$. The Atiyah--Segal
completion theorem generalizes to compact Hausdorff groups if we
make use of a suitable completion functor \cite{fg}.

Let $G$ be a compact Hausdorff group. We define a Thom spectrum as
$TO_G (V) = \colim_U TO_{G/U} (V)$ where the colimit is over $U \in
\lie (G)$ such that $V$ has a trivial $U$--action. For more details
see \cite{fg}.

\subsection{The levelwise $\ww$--model structure on orthogonal
$G$--Spectra} \label{Wmodelstructure}
We make some minor modifications to the discussion of model
structures in \cite[III]{mma}.

Let $R$ be  an algebra.
The category of $R$--modules can be described as the category of
continuous $\dd$--spaces for an appropriate diagram category $\dd$.
The objects are the same as those of $ \jj^{\vv}_G$, but the
morphisms are more elaborate. See \cite[Section II.4]{mma} and
\cite[Section 23]{mms}  for details. Interpreted as
a continuous diagram category in $G \ttt$, we give $\mm_R$ the
model structure with weak equivalences and fibrations inherited from
the $\ww$--model structure
on $G \ttt$ \cite[11.3.2]{hir}. Recall Definition \ref{generators}.

\begin{defn} Let \mdfn{$\Sigma^{\infty}_R \ww  I$} denote
the collection of $\Sigma_{V}^R \, i$, for all $i \in \ww I$
and all indexing representations $V$ in $\uu$. Let
\mdfn{$\Sigma^{\infty}_R \ww J$} denote the collection of
$\Sigma_{V}^R \, j$, for all $j \in \ww J$ and all indexing
representation $V$ in $\uu$.
\end{defn}

\begin{lem} \label{lem:hurwicscofibration}
If $f \col X \rarr Y$ is a relative $\Sigma^{\infty}_R \ww I$--cell
complex, then $ X(V ) \rarr Y(V)$ is a $G$--equivariant Hurewicz
cofibration (satisfies the homotopy extension property), for every
$V \in \vv$.  All Hurewicz cofibrations are closed inclusions.\end{lem}
\begin{proof} The adjunction between $\Sigma_{V}^R$ and  $\Omega^{\infty}_V$
gives that the maps in the classes $\Sigma^{\infty}_R \ww I$ and
$\Sigma^{\infty}_R \ww J$ are Hurewitz cofibrations.
Hence relative $\ww$--cell complexes, and their retracts, are Hurewitz
cofibrations. The last claim follows  since our spaces are weak
Hausdorff. \end{proof}

\begin{lem} \label{lem:spectratopstructure}
Let $\E$ be  a $\ww$--\good\ collection.  Assume that
$ f \col  A \rarr B$ is in $\Sigma^{\infty}_R \E I$
and $g \col X \rarr Y$ is in $\Sigma^{\infty}_R \ww I$, then
$f \Box g \col (A \wedge Y) \cup_{A \wedge X} ( B \wedge X)
\rarr B \wedge Y$ is a relative  $\Sigma^{\infty}_R \ww I$--cell
complex. Moreover, if  $f$ is in $\Sigma^{\infty}_R \E J $ instead
of $\Sigma^{\infty}_R \E I$ or $g$ is in $\Sigma^{\infty}_R \ww J$
instead of  $\Sigma^{\infty}_R \ww I$,  then $f \Box g$
is a relative $\Sigma^{\infty}_R \ww J$--cell complex.
\end{lem}

\begin{proof} This reduces to the analogue for spaces, Lemma
\ref{lem:topstructure}, since $\Sigma^{R}_V $ respects colimits
and smash products from spaces to $R$--modules. \end{proof}

The following model structure on the category of orthogonal
$G$--spectra is called the \mdfn{ levelwise $\ww$--model structure}.
\begin{pro}
\label{pro:stablems} Let $\ww$ be a collection of subgroups of $G$.
 Then the  category of $R$--modules has a compactly
generated proper model structure with levelwise $\ww$--weak
equivalences and levelwise $\ww$--fibrations (as
$\jj^{\vv}_G$--diagrams). The cofibrations are generated by
$\Sigma^{\infty}_R \ww I$, and the acyclic cofibrations are
generated by $\Sigma^{\infty}_R \ww J$. If $\E$ is  a $\ww$--\good\
collection, then the model structure is $\E $--topological.
\end{pro}
\begin{proof}
The source of the maps in $\Sigma^{\infty}_R \ww I$ and $ \Sigma^{\infty}_R \ww J$ are small  since $\Omega^{\infty}_V$ respects sum. Let $f$ be a relative  $ \Sigma^{\infty}_R \ww J $--cell complexes. Then  $  \Omega^{\infty}_V f$, for any indexing representation $V$ in $\uu$, is the colimit of a sequence of equivariant homotopy  equivalences that are (closed) Hurewicz cofibrations, hence a  levelwise $\ww$--equivalences of $G$--spaces. It  follows that  $\mm_R$ inherits a model structure from $G \ttt$ via the set of right adjoint functors $\Omega^{\infty}_V$, for  indexing representations $V$ in $\uu$. The model structure is right proper since $\ww \ttt$ is right proper. It is left proper by Lemma \ref{lem:hurwicscofibration} because fixed points respects pushout along a closed inclusion and since pushout of a weak equivalence along  a closed Hurewicz cofibration of spaces is again a weak equivalence.
The last claim follows from Lemma \ref{lem:spectratopstructure}.
See also \cite[6.5]{mms}.
\end{proof}

\begin{defn} The $V$--loop space $\Omega^V Z $ of a $G$--space
$Z$ is $ \ttt_G (  S^V , Z ) $. A  spectrum $X$ is called a
\mdfn{$\ww$--$\Omega$--spectrum} if the adjoint of the structure maps,
$ X (V) \rarr \Omega^{W-V} X( W)$ are  $\ww$--equivalences of
spaces for all pairs $V \subset W$ in $\vv (\uu )$.
\end{defn}

\section{The stable $\ww$--model structure on orthogonal $G$--spectra}
We define stable equivalences between orthogonal $G$--spectra and
localize the levelwise $\ww$--model structure with respect to these equivalences
\cite[III.3.2]{mma}.  Recall that if $Z$ is a based $G$--space,
then $ \pi^{H}_n ( Z) $  denotes the $n$-th homotopy group of $Z^H$.
We say that a map of based spaces, $f \col X \rarr Y$, is a
\mdfn{ based $\ww$--equivalence} if $ \pi_{n}^U (f) $ is an
isomorphism for all $U \in \ww $ and all $n \geq 0 $.
Let $f \col X \rarr Y$ be a based $G$--map.
If $f$ is a unbased $\ww$--equivalence
as defined in Definition \ref{defn:Wweakequivalencesspaces}, then
$f$ is a based $\ww$--equivalence.
If $f$ is a based $\ww$--equivalence, then
$\Omega f \col \Omega X \rarr \Omega Y$ is a unbased $\ww$--equivalence.

\begin{defn}  The \mdfn{$n$-th  homotopy group}
of an orthogonal $G$--spectrum $X$ at a subgroup $H $ of $G$ is
\[ \pi^{H}_n (X) = \colim_{V } \, \pi^{H}_n ( \Omega^V X (V) ) \]
for $n \geq 0$, and \[ \pi^{H}_{-n} (X) = \colim_{V \supset
\mathbb{R}^n } \, \pi^{H}_0 ( \Omega^{V-\mathbb{R}^n} X (V) ) \] for
$n \geq 0$, where the colimit is over indexing representations in
$\uu$. A map $f \col X \rarr Y$ of orthogonal $G$--spectra is a
\mdfn{stable ${\ww}$--equivalence} if $\pi^{H}_n (f) $ is an
isomorphism for all $H \in \ww $ and all $n \in \zz$.
\end{defn}

We follow our program of giving model structures to the category of
orthogonal $G$--spectra with minimal assumptions on the collection of
subgroups used. We need to impose conditions on the  stabilizers of
the indexing representations.  Let $\stable (\uu )$ be the collection
of stabilizers  of elements in the $G$--universe
$\uu$. We need  $\stable (\uu )$ to be  a $\ww$--\good\
collection.  If $\uu$ is a trivial
$G$--universe, then $\stable (\uu )$ is a
$\ww$--\good\ collection for any collection $\ww$ (see
Definition \ref{good}). If
$\uu$ is a complete universe, then $\stable (\uu )$   is a $\ww$--\good\
collection if $\ww$ is a family in the cofamily closure of $\lie
(G)$ (see Definition \ref{defn:normal}). If $\ww$ is an \good\
collection of subgroups of $G$ containing $G$ itself, then $\stable (\uu )$ is
$\ww$--\good\ if and only if $\stable (\uu )$ is contained in $\ww$.

The category $\mm_R$ is a closed symmetric tensor category. We
follow \cite[2]{ss} when considering the interaction of model
structures and tensor structures. A model structure is said to be
tensorial if the following pushout--product axiom is valid.

\begin{defn} \label{def:monoids}
\mdfn{Pushout--product axiom} \cite[2.1]{ss}: Let $f_1 \col X_1
\rarr Y_1$ and $f_2 \col X_2 \rarr Y_2 $ be cofibrations. Then the
map from the pushout, $\text{P}$, to $Y_1 \wedge Y_2$ in the
diagram \[ \xymatrix{ X_1 \wedge X_2 \ar[r]^{f_1 \wedge 1}
\ar[d]_{1 \wedge f_2} & Y_1 \wedge
X_2 \ar[d] \ar[rdd]^{1 \wedge f_2}  &  \\
X_1 \wedge Y_2 \ar[r] \ar[rrd]_{f_1 \wedge 1} & \text{P}
 \ar[rd] & \\
 &  &  Y_1 \wedge Y_2 , } \] is again a cofibration.  If, in
addition, one of the maps $f_1$ or $f_2$ is a weak equivalence, then
$\text{P} \rarr Y_1 \wedge Y_2 $ is also a weak equivalence.
\end{defn}

\begin{defn} \mdfn{Monoid axiom} \cite[2.2]{ss}: Any acyclic cofibration
tensored with an arbitrary object in $\mm$ is a weak equivalence.
Moreover, arbitrary pushouts and transfinite compositions of such
maps are weak equivalences.
\end{defn}

The following model structure on $\mm_R$ is called the (stable) $\ww$--model structure.
We sometimes denote $\mm_R$ together with the $\ww$--model structure
by \mdfn{$ \ww \mm_R$}.

\begin{thm} \label{spectramodel}
Let $R$ be an algebra.
Assume that  $\stable (\uu) $ is a   $\ww$--\good\ collection
of subgroups of $G$.  The category of $R$--modules is a compactly generated
proper  model category such that the weak
equivalences are the stable $\ww$--equivalences, the cofibrations
are retracts of relative $\Sigma_{R}^{\infty} \ww I$--cell
complexes, and a  map $f \col X \rarr Y$ is a fibration if and only if the map $f (V) \col
X(V) \rarr Y(V)$ is a $\ww$--fibration and the obvious map from $X(V)
$ to the pullback of the
diagram  \[ \xymatrix{ & \Omega^W X (V \oplus W) \ar[d] \\
Y(V) \ar[r] & \Omega^W Y(V \oplus W) } \] is a unbased
$\ww$--equivalence of spaces for all $V, W \in \vv$.  A map $f \col X
\rarr Y$ is an acyclic fibration if and only if $f$ is a levelwise
acyclic fibration.

Assume $R$ is symmetric. If
$\E$ is  a $\ww$--\good\ collection, then the model structure is
$\E $--topological. If $\ww$ is an \good\ collection, then the
model structure satisfies the pushout--product axiom and the monoid axiom.
\end{thm}

Remarks \ref{rem:topms} and \ref{simplicialstructure} imply
that the $\ww$--model structure is simplicial.
The fibrant spectra are exactly the $\ww$--$ \Omega$--spectra.
If $X$ is a cofibrant $\sss$--module, then a fibrant replacement is given by
\begin{equation} \label{fibrantreplacement} Q X (V) =  \hocolim_W \, \Omega^W X(V \oplus W) \end{equation} together with the natural transformation $ 1 \rarr Q$.

The stable homotopy group $ \pi^{H}_n ( X)$ is isomorphic to
$\pi^{H}_n ( Q X' )$, where $X'$ is a cofibrant replacement of $X$
in the category of $\sss $--modules with the levelwise $\ww $--model structure.
Hence we get the following.
\begin{lem} \label{corep} Let $H $ be in $ \ww$.
The stable homotopy group $\pi_{n}^H $ is corepresented by $
\Sigma^{R}_{\mathbb{R}^{-n}} G/H_+ $, for $n \leq 0 $, and by $
\Sigma_{R}^{\infty} G/H_+ \wedge S^n $, for $n \geq 0$, in the
homotopy category of $\mm_R$ with the stable $\ww$--model structure.
Hence $\pi_{n}^H $ is a homology theory which satisfies the colimit axiom.
\end{lem}

The \mdfn{colimit axiom} says that $\colim_a \, \pi_{* }^H ( X_a )
\rarr \pi_{* }^H (X)$ is an isomorphism, where the colimit is over
all finite subcomplexes $X_a$ of the cell complex $X$.

\begin{pro}
\label{pro:MRtensormodelcat} Let $\uu$ be a complete $G$--universe.
Let  $\ww$ be a family in $\overline{ \lie (G )}$.
Then the dualizable objects in the homotopy category of
$\ww \mm_R $ (with the $\ww$--model structure) are precisely
retracts of $\vv ( \uu )$--desuspensions of finite $\ww$--cell
complexes.
\end{pro}
\begin{proof}
The proof in \cite[XVI 7.4]{ala} goes through with modifications to
allow for general $R$--modules instead of $\sss$--modules.
\end{proof}

\subsection{Verifying the model structure axioms}

\begin{lem} \label{eq} Let $ X$, $Y$ and $Z$ be based $G$--spaces.
Let $\E$ be a $\ww$--\good\ collection.
If $Z$ is a $\E I$--cell complex and $f \col  X \rarr Y$
is a $\ww$--equivalence, then $ \ttt_G ( Z, X ) \rarr \ttt_G (Z ,Y)$
is a $\ww$--equivalence.
\end{lem}
\begin{proof} Since $  S^k \wedge S^{n}_+$ is a based
CW--complex and $(G/H \times G/L)_+$ is a  $\ww I$--cell complex,
for all $H \in \ww$ and $ L \in \E$, it follows that
$ S^k \wedge G/ H_+ \wedge Z$ is a  $\ww I$--cell based complex.
The smash product and internal hom adjunction  gives that
\[ \pi^{H}_{k} (\ttt_G (Z , X)) \cong
[ S^k \wedge G/ H_+ \wedge Z ,X ]_G \] where the square brackets
are $G$--homotopy classes of maps (see Proposition \ref{homotopyclasses}).
Since $ \pi^{K}_k ( f) \cong [  S^k \wedge G/ K_+ , f ]_G $
is an isomorphism for all $K\in \ww$ and $k \geq 0 $, the result  follows by using
the higher lim spectral sequence.
\end{proof}

\begin{cor}
\label{cor:level-pi-star} Assume  $\stable (\uu )$ is a $\ww$--\good\
collection of subgroups of $G$. A levelwise $\ww$--equivalence of $G$--spectra is
a stable $\ww$--equivalence.
\end{cor}
\begin{proof}  By Definition \ref{def:universe}, any finite
dimensional $G$--representation $V$ in the universe $\uu$ is a
$G/U$--representation for a compact Lie group quotient $G/U$ of $G$. Since
the collection $\stable (\uu)$ is $\ww$--\good, it  follows from   Lemma \ref{eq}
that \[ \Omega^{V} f (V' ) \cong  \ttt_G ( S^V , \, f (V' ))  \] is a
$\ww$--equivalence for all $V , V' \in \vv$.
\end{proof}

Note that if $\stable ( \uu )$ is a  $\ww$--\good\ collection  of subgroups of
$G$, $A$ is a based $\ww I$--cell complex, and $V$ is an indexing
representation in $\uu$, then $A \wedge S^V$ is again a based
$\ww I$--cell complex by Lemma \ref{productcomplexes}. The next result,
together with Corollary \ref{cor:level-pi-star}, show that a map
between $\Omega$--$\ww$--spectra is a levelwise $\ww$--equivalence if
and only if it is a $\ww$--equivalence. This fundamental result is
an extension of \cite[III.9]{mma}.
\begin{pro} \label{prostrict} Assume that $\stable ( \uu )$ is a $\ww$--\good\
collection  of subgroups of $G$. Let $f \col X \rarr Y$ be a map of
$\ww$--$ \Omega$--spectra.  If \[ f_* \col  \pi_{*}^{H} ( X ) \rarr
\pi_{*}^{H} (Y) \] is an isomorphism for each $H \in \ww$, then
$f (V) \col  X (V)  \rarr Y(V)$ is a  $\ww$--equivalence for
all indexing representations $V \subset \uu$.
\end{pro}
\begin{proof} We first prove that
\[ f(V)_* \col  \pi_{*}^{H} ( X (V) ) \rarr \pi_{*}^{H} (Y(V)) \]
is an isomorphism for all $H \in \ww$.
Let $Z$ be the homotopy fiber of $f$. It is again an $\Omega$--$ G$--spectrum.
 We want to show that $\pi_{*}^{H} ( Z ) =
0 $ for all $H \in \ww$, implies that $ \pi_{*}^H ( Z(V)) = 0$ for
any indexing representations $V$ and any $H \in \ww$. Fix an
indexing representation $V$ and a normal subgroup $N \in \lie (G)$
such that $N$ is a finite intersection of elements in $\stable  (\uu )$ and
$N$ acts trivially on $V$. With these choices $ (\Omega^V
Z(V))^H = \Omega^V ( Z(V)^H)$ for all $H \leq N$. Hence $ \pi_{* +
|V|}^{H} ( Z(V))$ is isomorphic to $\pi_{*}^H ( \Omega^V Z(V) ) $
for all $H \leq N$ in $\ww$. Since $Z $ is an $\Omega$--$ G$--spectrum,
an easy argument gives that $\pi^{H}_* ( Z(V)) = 0$ for all $H \in
\ww$ such that $H \leq N$ \cite[III.9.2]{mma}.

We now prove the result for subgroups $H$ in $\ww$ that are not
necessarily contained in $N$. Fix a subgroup $H \in \ww$. Assume by
induction that $\pi_{\ast}^{K} (Z(V)) = 0$ for all subgroups $K \in
\ww$ such that $ K $ is properly contained in $ H $. If $L $ is a stabilizer of
$V$, then $H \cap gLg^{-1}$ is in $\ww$ for all $g \in G$ since $\stable  (\uu )$ is a
$\ww$--\good\ collection. The argument given in \cite[Section III.9]{mma}
implies that $\pi_{\ast}^{H} (Z(V)) = 0$. We now justify that we can
make the inductive argument. The quotient group $H / H \cap N $ is
isomorphic to $H \cdot N / N$, which is a subgroup of the compact
Lie group $G/N$. Hence the partially ordered set of closed subgroups
of $H$ containing $H \cap N$ satisfies the descending chain
property. We have that $ \pi_{*}^{K } ( Z(V)) = 0 $ for all $K \leq
H \cap N$ in $\ww$. We start the induction with the subgroup $H \cap
N$ which by assumption is in $\ww$. For more details see
\cite[Section III.9]{mma}.

Since $ X (V )^H   $ and $ Y (V )^H   $ are weakly equivalent to
$  \Omega^{\rr} X ( V \oplus \rr)^H $ and
$  \Omega^{\rr} Y ( V \oplus \rr)^H $, respectively, and
$\pi_{k +1}  \col \pi_{k+1} (X ( V \oplus \rr)^H  )
\rarr \pi_{k +1}  ( X ( V \oplus \rr)^H ) $   is an isomorphism
for each $H \in \ww$, $ k \geq 0$ and indexing representation $V$,
it follows that $ f (V)$ is a $\ww$--equivalence.
So $f$ is a levelwise   $\ww$--equivalence.
\end{proof}

A set of generating cofibrations is $\Sigma^{\infty}_R \ww I$. We
give a set of generating acyclic cofibrations. Let \[ \lambda_{V,W}
\col \Sigma^{R}_{V \oplus W } S^W \rarr \Sigma_{V}^R S^0 \] be the
adjoint of the map
\[ S^W \rarr  ( \Sigma_{V}^R S^0) ( V \oplus W ) \cong \oo
( V \oplus W)_+ \wedge_{\oo (W)} R(W) \] given by sending an element
$w$ in $S^W$ to $e \wedge \iota (W)(w)$ where $e$ is the identity map in
$\oo (V \oplus W )$, and $\iota \col S^0 \rarr R $ is the unit map. Let
\mdfn{$k_{V ,W}$} be the map from $\Sigma^{R}_{V \oplus W } S^W$ to
the mapping cylinder, $M \lambda_{V ,W}$, of $\lambda_{V ,W}$. Let
\mdfn{$\ww K$} be the union of $\Sigma^{\infty}_R \ww J$ and the set
of maps of the form $ i \Box k_{V ,W}$ for $i \in \Sigma^{\infty}_R \ww I$ and
indexing representations $V,W$ in $\uu$.  The box denotes  the pushout--product map.
The map $ \lambda_{V,W}$ is a $\ww$--equivalence and $ k_{V ,W}$ is a levelwise $\stable (\uu )$--cofibration  for all indexing  representations $V,W$ in $\uu$
\cite[Lemma III.4.5]{mma}. Hence all maps in $\ww K$ are both cofibrations
and $\ww$--equivalences.  The set $\ww K$ of maps in $\mm_R$ is a set of generating acyclic cofibrations.

As in \cite[III.4.8]{mma} and \cite[9.5]{mms}, the
following characterization of the maps that satisfy the right
lifting property with respect to $\ww K$ follows since the
$\ww$--model structure on $G \ttt$ is $\stable (\uu)$--topological.

\begin{pro} \label{fibchar}  A map $p \col  E \rarr B$ satisfies the right lifting
property with respect to $\ww K$ if and only if $p$ is a levelwise
fibration and the obvious map from $E(V) $ to the pullback of the
diagram  \[ \xymatrix{ & \Omega^W E (V \oplus W) \ar[d] \\
B(V) \ar[r] & \Omega^W B (V \oplus W) } \] is an unbased
$\ww$--equivalence of spaces for all $V, W \in \vv (\uu)$.
\end{pro} \begin{proof}
If $k$ is a $\stable (\uu )$--cofibration, and $p$ is a levelwise $\ww$--fibration,
then by Proposition \ref{pro:topmodstructure} $ (\mm_R)_G ( k^* ,
p_* )$ is a $\ww$--fibration of spaces.
Hence the proof of \cite[Proposition III.4.8]{mma} gives the result.
\end{proof}

The proof of Theorem \ref{spectramodel} is similar  to the proofs in
\cite[Section III.4]{mma},\cite[III.7.4]{mma} and \cite[Section 9]{mms}. Note that
proofs of these results use a few lemmas, given in \cite{mma}, that are not made
explicit in this paper.  Properness follows as in \cite[9.10]{mms}.
The model structure is $\E $--topological as in \cite[9.9]{mms}.

\begin{rem}
An alternative construction  of this model structure is provided by the
general framework of Olivier Renaudin \cite{ren}.  The functorial
fibrant replacement functor for the stable $\ww$--model structure
gives a homotopy idempotent functor in the levelwise  model structure on $\mm_R$.
The description of the fibrations follows from the  description of
the fibrant replacement (of a cofibrant object)  \ref{fibrantreplacement}
and a  result of  Bousfield \cite[Theorem 9.3]{bou2}.
\end{rem}

If  $\ww$ is an \good\ collection of subgroups
of $G$,  then the $\ww   $--model
structure on $\mm_R$ is a  tensor
model structure that satisfies the monoid axiom.
 This follows as in \cite[Section III.7]{mma}

\subsection{Positive model structures}

We give some brief remarks about other model categories of spectra.
 Prespectra are defined by replacing the
category ${\vv (\uu)}$, in Definitions \ref{defn:jvg} and
\ref{defn:orthogonalspectrum}, by the  smaller category consisting of
 the indexing representations and their inclusions.
 There is a stable $\ww$--model structure on
the category of
 prespectra. This model category is Quillen equivalent to
 the stable  $\ww$--model structure on
 $G$--orthogonal spectra \cite[III.4.16]{mma}.

We can also consider model structures on the category of algebras.
We need to remove some of the cofibrant and acyclic cofibrant
generators to make sure the free symmetric algebra construction
takes acyclic cofibrant generators to stable $\ww$--equivalences. Let
$\Sigma^{R}_+ \ww I $ and $\Sigma^{R}_+ \ww J$ consist of all
$V$--desuspensions of elements in $\ww I$ and $\ww J$ by indexing
representations $V$ in $\uu$ such that $V^G \not= 0$. The positive
levelwise $\ww$--model structure on the category of
orthogonal spectra is the model structure obtained by replacing
$\Sigma^{\infty}_R \ww  I$ and
$\Sigma^{\infty}_R \ww J$ by $\Sigma^{R}_+ \ww I $ and $\Sigma^{R}_+
\ww J$, respectively. The positive stable $\ww$--model structure on
orthogonal spectra is obtained by replacing $\ww K$ by the set $\ww
K_+$ consisting of the union of $\Sigma^{R}_+ \ww J $ and the maps $i
\square k_{V,W}$ with $ i \in \ww I$ and $V^G \not= 0$.
 The discussion of the positive model structure goes
through as in \cite[Sections III.5, III.8]{mma}.

 \begin{pro} Let $R$ be a commutative monoid in the
category of $G$--orthogonal spectra. Let $\ww$ be an \good\
collection containing $\stable ( \uu)$.  Then there is a compactly
generated $\ww$--topological model structure on the category of
$R$--algebras such that the fibrations and weak equivalences are created
in the underlying positive $\ww$--model category of orthogonal
$G$--spectra. The same applies to the category of commutative $R$--algebras.
\end{pro}

\subsection{Fibrations}
We consider the behavior of   fibrations in the stable model structure
under   restriction  to subgroups. See Remark \ref{rem:restriction}.

\begin{defn} \label{defn:intersection}
Let $\ww$ be a collection of subgroups of $G$, and let
$K$ be a subgroup of $G$. The intersection \mdfn{ $ K \cap \ww$} is
defined to be the collection of all subgroups $H \in \ww$ such that
$ H \leq K$. \end{defn}  If $\ww$ is a Lie collection of subgroups
of $G$ that is closed under intersections and $\stable_G (\uu) $ is
a $\ww$--\good\ collection of subgroups of $G$, then, for every
$K \in \ww$, the collection of $K$--stabilizers $\stable_K ( \uu |K)
$ is a $ K \cap \ww $--\good\ collection of subgroups of $K$.

\begin{lem} \label{lem:easyrestrictfib}
Let $\stable (\uu) $ be a $\ww$--\good\ collection  of subgroups of $G$ and let
$K \in \ww$. Let $f \col X \to Y$ be a fibration  in $ \ww$--$G \mm_R$. Then
$f$ regarded as a map of $K$--spectra  is a fibration  in $(K \cap \ww )$--$ K
\mm_R$.  \end{lem}
\begin{proof}
This follows from the explicit description of fibrations  in
Theorem \ref{spectramodel}. (Alternatively, check that $
G\wedge_K -$ is left Quillen adjoint to the forgetful functor from
$G$--spectra to $K$--spectra.)
\end{proof}

Lemma \ref{lem:easyrestrictfib} need not remain true when the
subgroup $K$ is not in $\ww$. For applications in Section
\ref{sec:htpfixedpt} we give some conditions  that guarantee that
the result remains true even when $ K \not\in \ww$.

\begin{lem} \label{lem:restrictfib}
Let $\stable (\uu) $ be a $\ww$--\good\ collection of subgroups
of a compact Hausdorff group $G$. Let $f \col X \rarr Y$ be a
fibration in $\ww  \mm_R$. Assume that both $X$ and $Y$ are
$\ww$--$\sss$--cell complexes.  Let $K$ be any closed subgroup of
$G$, and let $\ww'$ be a  collection of subgroups of
$K$ such that $\ww' \ww \subset \ww$ and $\stable_K (\uu | K) $
is a $\ww'$--\good\ collection of subgroups of $K$.
Then $f$, regarded as a map of $K$--spectra, is a
fibration in the $\ww' $--model structure on $ K \mm_R$ (indexed on
$\vv ( \uu )  |K$, and $R$ regarded as a $K$--spectrum).
\end{lem}
Note that $X$ and $Y$ are  required to be $\ww$--$ S $--cell complexes
not just  $\ww$--$ R$--cell complexes. This holds if they are
$\ww$--$ R $--cell complexes and $R$ is a $\ww$--$ S $--cell complex.
\begin{proof} Let $ f \col X \rarr Y$ be a $\ww$--fibration between
$\ww$--cofibrant objects in $\ww G \mm_R$. Since $\stable_K (\uu | K) $
is a $\ww'$--\good\ collection, it suffices, by Theorem
\ref{spectramodel} and Proposition \ref{fibchar}, to show that
for any $L\in \ww'$: (1)  the map $f
(V)^L \col X (V)^L \rarr Y(V)^L$ is a fibration, for $V \in \vv (\uu )| K$,
and  (2) the map from $ X (V)^L$ to the pullback of the diagram
\begin{equation}
\label{eq:htpullback}  \xymatrix{  & Y (V)^L \ar[d] \\
( \Omega^{W} X (V \oplus W))^L \ar[r] & (\Omega^{W} Y (V \oplus
W))^L } \end{equation} is a weak  equivalence of spaces, for $V , W \in
\vv (\uu )| K$.

We now prove that maps from  a compact space $C$ to the $L$--fixed points
of $X(V)$ and $Y(V)$ factor through the $UL$--fixed points of $X(V)$
and $Y(V)$ for some $U \in \ww$. Since $X$ and $Y$ are $\ww$--$
\sss $--cell complexes, and a  map from a compact space $C$ into a
$\ww$--$\sss $--cell
complex factors through a finite sub cell complex, it suffices
to verify the claim for individual cells.
Recall, from \ref{eq:R-suspension}, that
$\Sigma^{\infty}_{V'} G/H_+ \wedge D_{+}^n (V) $ is the space $G/H_+ \wedge D_{+}^n
 \wedge \oo ( V )_+ \wedge_{ \oo
(V - V' )} S^{ V - V' }$, for $ V \supset V'$, and a point
otherwise. This is a finite  $\ww$--cell complex by
Lemma \ref{productcomplexes}  and Illman \cite{ill}.
We conclude that a map from a compact space $C$ into the $L$--fixed points of
$X(V) $ and $Y (V)$ factor through $X ( V)^{UL} $ and $Y ( V)^{UL}
$, respectively, for some $U \in \ww$.

We are now ready to prove   claim (1). Let \begin{equation}
\label{eq:fib}
\xymatrix{ D^{n}_+ \ar[r] \ar[d]_j & X(V)^L \ar[d]^{f (V)^L} \\
(D^{n} \times I)_+ \ar[r] & Y(V)^L } \end{equation} be a diagram of
based spaces. There exists a $U \in \ww$ such that the map from $j $
to $f(V)^L $ factors through $f(V)^{UL} $. Since $f(V)^{UL} $ is a
fibration we get a lift in the diagram \ref{eq:fib}. Hence $ f
(V)^L$ is a fibration.

The proof of claim (2) is similar. We note that a map from a
compact space $C$ to $(\Omega^W X ( V \oplus W))^L$, composed with
the inclusion into $\Omega^W X ( V \oplus W)$, is adjoint to a based
map from $ C_+ \wedge S^W$ to $ X ( V \oplus W)$. Hence it factors
through $X (V \oplus W)^{U'}$, for some $U' \in \ww$. By choosing a
smaller $U \leq U'$ such that $U$ acts trivially on $W$, the map
from $C$ factors through $ (\Omega^W X ( V \oplus W))^{UL}$. Hence
to check that the map from $X(V)^L$ to the pullback of
\ref{eq:htpullback} is a weak equivalence, it suffices to check this
on all $UL$--fixed points for $U \in \ww$. This follows by our
assumptions.
\end{proof}

\begin{lem} \label{lem:restrict-n-equivalence}
Let $\stable (\uu) $ be a $\ww$--\good\ collection of subgroups of a compact
Hausdorff group $G$. Let $f \col X \rarr Y$ be a  (co--)$n$--equivalence
in $ \mm_R$ between $\ww$--$ \sss$--cell complexes $X$ and $Y$ which are also
fibrant objects in the $\ww$--model structure on $G \mm_R$.
Let $K$ be any closed subgroup of $G$, and let $\ww'$ be
a collection of subgroups of $K$ such that $\ww'
\ww \subset \ww$ and $\stable_K (\uu | K) $ is $\ww'$--\good.
Then $ f$ regarded as a map of $K$--spectra is a
(co--)$n$--equivalence  in the $\ww' $--model structure on $ K \mm_R$.
\end{lem}
\begin{proof}
The spectra $X$ and $Y$  are also fibrant in the $\ww'$--model
structure on $ K \mm_R$ by Lemma \ref{lem:restrictfib}.
Hence it suffices to prove that  $ X(\mathbb{R}^m  )^L \rarr Y
(\mathbb{R}^m )^L$ is a (co--)$(n-m)$--equivalence for any $L \in \ww'$
and $n \geq m \geq 0$. Since both $X$ and $Y$ are fibrant in the
$\ww$--model structure on $G \mm_R$, it follows that  \[X(\mathbb{R}^m )
\rarr Y (\mathbb{R}^m )\] is a $\ww$--(co--)$(n-m)$--equivalence.
The proof of Lemma \ref{lem:restrictfib} gives the result.
\end{proof}

\section{The $\ww$--$\cc$--model structure on orthogonal $G$--spectra}
Let $R$ be a ring and assume that $\stable (\uu ) $ is a  $\cc$--\good\
collection of subgroups of $G$. We define $K$--equivalences in the
$\cc$--model structure on the category of $R$--modules, $\mm_R$, for
$K$ not necessarily in $\cc$. This is used to construct a model structure
with weak equivalences detected by a collection $\ww$ of subgroups
of $G$ which is not necessarily contained in $\cc$. We start by
briefly describing the resulting  $\ww$--$\cc$--model structure on
$\mm_R$ in the case when $ \ww$ is contained in $\cc$.

Let $H $ be in $\cc$. Then $\pi_{* }^H $ is a homology theory which
satisfies the colimit axiom by Lemma \ref{corep}. The direct sum
\[h =  \textstyle\bigoplus_{K \in \ww , n \in \zz } \, \pi^{K}_n \]
is also a homology theory which satisfies the colimit axiom.
We can now (left) Bousfield localize $\cc\mm_R$ with
respect to the homology theory $h$ \cite{bou} \cite[13.2.1]{hir}.
Hence for any subcollection  $\ww$ in
$\cc$ there is a model structure on $G$--spectra such that the
cofibrations are  retracts of relative
$\cc$--cell complexes and the weak equivalences are maps $f$ such
that $\pi_{n}^H (f) $ is an isomorphism for all $H \in \ww$ and
$n \in \zz$.

\subsection{The construction of $\ww \cc \mm_R$}
Assume that $\cc$ is an \good\ collection of subgroups of $G$.
\begin{defn} \label{kequivalence} Let $K$ be a  subgroup
of $G$ such that the closure $\overline{UK} \in \cc$, for all $U \in
\cc$. The $n$-th stable homotopy group at $K$ is defined to be
\mdfn{$ \Pi^{K}_* ( X) = \colim_{ U \in \cc} \, \pi^{\overline{UK}
}_* ( X) $}.  \end{defn}
The colimit is over the
category with objects  $U$ in $\cc$ and
with morphisms containment of subgroups. The colimit is directed
since $\cc$ is an \good\ collection.
If  $K \in \cc$, then $\Pi_{\ast}^K$ and
$\pi_{\ast}^K$ are  canonically isomorphic functors.
\begin{defn} Let $\cc$ and $\ww$ be two collections of subgroups of
$G$.  Then the product collection \mdfn{$\cc \ww$} has elements the
closure $\overline{UH}$ of the product subgroup $ U H $ in $G$,
for all $U \in \cc$ and all $H \in \ww$.
\end{defn}
The collection $\ww = \{ 1\} $
satisfies $\cc \ww \subset \cc$ for any collection $\cc$. If $\cc$
is a cofamily, then $\cc \ww \subset \cc $ for any collection $\ww$.
\begin{defn}
Let $\ww$ be a collection  of subgroups of $G$ such that
$ \cc \ww \subset \cc$. Then  a map $f$ between
orthogonal spectra is a \mdfn{ $\ww$--equivalence } if $\Pi^{K}_{n}
(f)$ is an isomorphism for all $K \in \ww$ and all integers $n$.
\end{defn}
Directed colimits of abelian groups respect direct sums and exact
sequences. So $\Pi^{K}_* $ is a homology theory which satisfies the
colimit axiom by Lemma \ref{corep}. The direct sum
\[h =  \textstyle\bigoplus_{K \in \ww , n \in \zz } \, \Pi^{K}_n \]
is again a homology theory which satisfies the colimit axiom. Hence
we can Bousfield localize with respect to $h$.
\begin{thm} \label{model2} Let $ \cc$ be an \good\ collection
of subgroups of $G$ such that $ \stable (\uu) \subset \cc$. Let $\ww$ be
any collection of subgroups of $G$ such that $\cc \ww \subset \cc$.
Then there is a cofibrantly generated proper simplicial model
structure on $\mm_R$ such that the weak equivalences are
$\ww$--equivalences and the cofibrations are retracts of relative
$\cc$--cell complexes.
\end{thm}

\begin{proof}
There exists a set
$\mathcal{K}$ of relative $\cc$--$ G$--cell complexes with sources
$\cc$--$G$--cell complexes such that a map $p$ has the right lifting property
with respect to all $h$--acyclic cofibrations with cofibrant source,
if and only if  $p$ has the right lifting
property with respect to $\mathcal{K} $. To find such a set of maps
$\mathcal{K}$ we use the cardinality argument of Bousfield,
 taking into account both the cardinality of $G$ and the cardinality
of $\prod_{V } R (V)$, where the product is over indexing
representations in the universe $\uu$ \cite{bou}. The class of
$h$--equivalences is closed under pushout along $\cc$--cofibrations.
Hence we can apply \cite[13.2.1]{hir} to conclude that if $p$ has the
right lifting property with respect to the maps in the set $\mathcal{K} $,
then it has the right lifting property with respect to all $h $--acyclic
cofibrations. Hence there is a cofibrantly generated left
proper model structure on $\mm_R$ with the specified class of
cofibrations and weak equivalence \cite[4.1.1]{hir}. It remains to
show that the model structure is right proper and simplicial.

The model structure is right
proper by comparing homotopy fibers in   pullback diagrams \cite[9.10]{mms}.
We show that the model structure is simplicial. See Remark
\ref{simplicialstructure}.  The tensor and
cotensor functors are given by $ \Sigma^{\infty}_R | K_+ | \wedge X$
and $F (\Sigma^{\infty}_R | K_+ | , X ) $, respectively, for a
simplicial set $K$ and an $R$--module $X$. The simplicial hom functor
is given by $ \sing G \mm_R (X ,Y )$. The
pushout--product map applied to a simplicial cofibration and a
$\cc$--cofibration in $\mm_R$ is again a $\cc$--cofibration. If the
simplicial cofibration is acyclic, then the pushout--product map is
a $\cc$--acyclic cofibration.
It suffices to show that if $ X_2 \rarr Y_2 $ is a
$\ww$--$\cc$--acyclic cofibration with $\cc$--cofibrant source,
then the map from the pushout of
\[ \xymatrix{ \Sigma^{\infty}_R S^{n-1}_+  \wedge X_2 \ar[r] \ar[d] &
\Sigma^{\infty}_R D^{n}_+ \wedge X_2  \\
\Sigma^{\infty}_R S^{n-1}_+ \wedge Y_2 & } \] to $ \Sigma^{\infty}_R
D^{n}_+ \wedge Y_2 $ is again a $\ww$--$\cc$--acyclic cofibration
\cite[2.3]{ss}. This is the case since our weak equivalences are
given by a homology theory on the homotopy category of the tensor
$\cc$--model structure on $\mm_R$ (see Theorem \ref{spectramodel}).
\end{proof}
This  model structure is called the \mdfn{$\ww$--$\cc $--model structure}
on $\mm_R$. The $\ww$--model structure is the $\ww$--$\ww$--model
structure. We sometimes denote $\mm_R$ together with the $\ww$--$\cc
$--model structure by \mdfn{$ \ww \cc \mm_R$}.
\begin{pro}
\label{pro:c1c2} Let $\cc_1 \subset \cc_2$ be two $\uu$--\good\
collections of subgroups of $G$ containing the trivial subgroup, $1$,
and let $\ww$ be a collection of
subgroups of $G$ such that $ \cc_1 \ww \subset \cc_1$ and $ \cc_2
\ww \subset \cc_2$. Then the identity functors $\ww \cc_1 \mm_R
\rarr \ww \cc_2 \mm_R $ and $\ww \cc_2 \mm_R \rarr \ww \cc_1 \mm_R $
are left and right Quillen adjoint functors, respectively. Hence a
Quillen equivalence.
\end{pro}

\subsection{The $\cc$--cofree model structure on $\mm_R$}
The $\ww$--$\cc$--model structure on $\mm_R$ is of particular interest
when $\ww = \{ 1 \} $.
\begin{defn}
We say that $f$ is a \mdfn{ $\cc$--underlying equivalence} if
\[ \Pi^{1  }_* (f) = \colim_{U \in \cc } \, \pi^{U}_* (f) \] is an
equivalence.
\end{defn}
The name, $\cc$--underlying equivalence, is justified by the next lemma.
\begin{lem} \label{lem:enderlyingequivalence=ccequiv}
Assume that $G$ is a compact Hausdorff group and let $\cc$ be
 the collection
$\overline{\lie (G)}$. Let $R$ be the sphere spectrum $\sss$.
 Then a map  $f \col  X \rarr Y$  between cofibrant
objects (retracts of $\cc$--cell complexes) is a $\cc$--underlying
equivalence if and only if $f$ is a non-equivariant weak
equivalence.
\end{lem}
\begin{proof} This follows as in the proof of Lemma
\ref{lem:restrictfib}.
\end{proof}

\begin{thm} \label{Goerss} Assume that $\uu$ is a
trivial $G$--universe. Let $ \cc$ be an
\good\ collection of subgroups of $G$.
Then there is a cofibrantly generated proper simplicial model
structure on $\mm_R$ such that the weak equivalences are
$\cc$--underlying equivalences and the cofibrations are retracts of
relative $\cc$--cell complexes.
\end{thm}
\begin{proof} This is a special case of Theorem \ref{model2}.
\end{proof}

We refer to this model structure as the \mdfn{$\cc$--cofree model
structure on $\mm_R$}.

\begin{rem} The universe is required to be trivial as part of the
definition of the $\cc$--cofree model structure. When $\{1 \} $ is in
$\cc$, or $\cc$ is a family in $\overline{\lie (G)}$,
then  there is no loss of generality in making this assumption.
In this case the
$\cc$--cofree model categories of $G$--spectra for a universe $\uu$
and its fixed point universe $\uu^G$ are Quillen equivalent \cite[V.1.7]{mma}.
This is so because the homotopy groups $ \pi^{ 1 }_*$ are
isomorphic for the universes $\uu^G$ and $\uu$.
\end{rem}

\section{A digression:
  $G$--spectra for noncompact groups}
In this section we consider an example of a model structure on
orthogonal $G$--spectra where the homotopy theory is ``pieced
together'' from the genuine homotopy theory of the compact Lie
subgroups of $G$. This example is inspired by conversations with
Wolfgang L\"uck. This section plays no role later in the paper.

The model structure we construct below in Proposition
\ref{digression} is in many ways opposite to the model structure (to
be discussed) in Theorem \ref{promodelcategory}: Compact Lie
subgroups versus compact Lie quotient groups, ind-spectra versus
pro--spectra, pro-universes versus ind-universes. The difficulties
here lie in dealing with inverse systems of universes for the
compact Lie  subgroups of $G$.

Let $G$ be a topological group, and let $\uuu$ be a trivial
$G$--universe. Let $R$ be a symmetric monoid in the category of
orthogonal $G$--spectra indexed on $\uuu$. Let $ \mm$ denote the
category of $R$--modules indexed on $\uuu$.

\begin{defn} Let \mdfn{$\comp$} denote the family of
compact Lie subgroups of $G$. \end{defn} If $G$ is a discrete group
or a profinite group then $\comp$ is the family of finite subgroups
of $G$.
The results in this section remain true if each  $\comp$ is
replaced by a subfamily $ \mathcal{F}'_G $ of $ \comp$, such that
$\mathcal{F}'_J \subset \comp'$ whenever $ J<G$.

By Proposition \ref{pro:stablems} there is a cofibrantly generated
model structure on $ \mm$ such that the cofibrations are retracts of
relative $\Sigma_{R}^{\infty} \comp I$--cell complexes and the weak
equivalences are levelwise $\comp$--equivalences. We would like to
stabilize $ \mm$ with respect to $H$--representations for all compact
Lie subgroups $H$ of $G$. An $H$--representation might not be a
retract of a $G$--representation restricted to $H$ (there might not
be any nontrivial $H$--representations of this form).

Our approach is to localize $ \mm$ with respect to stable
$H$--homotopy groups defined using a complete $H$--universe, one
universe for each $H$ in $\comp$.

\begin{defn} An \mdfn{$\comp$--universe}
consists of an $H$--universe $\uu_H$,  for each $H \in \comp$,
such that whenever $H \leq K$, then $\uu_K | H$ is a subuniverse
of $\uu_H $. For any subgroups $H \leq K\leq L$ the three
resulting inclusions of universes are required to be compatible.
\end{defn}

We say that the $\comp$--universe, $\{ \uu_H \}$, is complete if
$\uu_H$ is a complete $H$--universe for each $H \in \comp$.
\begin{lem} There exists  a complete $\comp$--universe. \end{lem}
\begin{proof} Choose a complete $H$--universe $\uu'_H$ for each $H
\in \comp$. Let $\uu_H$ be defined to be \[ \textstyle\bigoplus_{K
\geq H}\, (\uu'_K | H )
\] where the sum is over all $K \in \comp$ that contains $H$.
\end{proof}

Let $H$ be a compact Lie group. Then there is a stable model
structure on orthogonal $H$--spectra, indexed on a trivial
$H$--universe, which is Quillen equivalent to the ``genuine'' model
structure on orthogonal $H$--spectra indexed on a complete
$H$--universe. This is proved in \cite[V.1.7]{mma} (note that the
condition $\vv \subset \vv'$  is not necessary).
Let $H$ be a compact Lie group and let $\vv$ and $\vv'$ be
collections of $H$--representations containing the trivial
$H$--representations. Typically, $\vv$ is the collection $\vv ( \uu)$
of all $H$--representations that are isomorphic to some indexing
representation in an $H$--universe $\uu$. There is a change of
indexing functor
\[ I^{\vv}_{\vv'} \col
 \jj^{\vv'}_G \sss  \rarr \jj^{\vv}_G \sss  \]
defined in \cite[V.1.2]{mma}.
The functor $I^{\vv}_{\vv'}$ induces  an equivalence of categories with $
I^{\vv'}_{\vv}$ as the inverse functor. The functor $I^{\vv}_{\vv'}$
is a strong symmetric tensor functor. These, and other, claims are
proved in \cite[V.1.5]{mma}.

\begin{lem} For each compact Lie subgroup $H$ of $G$ the
functor \[ \pi^{H}_* ( I^{\vv ( \uu_H)}_{\vv ( \uuu)} - )
\] is a homology theory on $\mm$ with the levelwise
$\comp$--$ \uuu$--model
structure, that satisfies the colimit axiom.
\end{lem}

\begin{proof}
This follows since $I^{\vv ( \uu_H)}_{\vv ( \uuu)}$ respects
homotopy colimits and weak equivalences
since $\vv ( \uuu) \subset \vv (\uu_H )$ \cite[V.1.6]{mma}.
The claim follows from Lemma \ref{corep}.
\end{proof}

We localize the stable $\comp$--$ \uuu$--model category with respect to
the homology theory given by \[h = \textstyle\bigoplus_H \,
\pi^{H}_* ( I^{\vv ( \uu_H)}_{\vv ( \uuu)} - ), \] where the sum is
over all $H \in \comp$.

\begin{pro} \label{digression}
Given an $ \comp$--universe $ \{ \uu_H \}$. There is a
cofibrantly generated   stable model structure on $\mm$ such
that the cofibrations are retracts of relative $\Sigma_{R}^{\infty}
\comp I$--cell complexes (for the trivial universe $\uuu$) and the weak equivalences are the
$h$--equivalences. If $\comp$ is an Illman collection of subgroups
of $G$ (see Definition \ref{good}), then the model structure satisfies the
pushout--product axiom.
\end{pro} This model structure is called the
\mdfn{ stable $ \{ \uu_H \}_{H \in \comp }$--model structure} on $\mm$.

\begin{proof} See the proof of Theorem \ref{model2}.
\end{proof}

The cofibrant replacement of $\Sigma_{R}^{\infty} S^0$ in this model
structure (regardless of $\{ \uu_H \}$) is given by $
\Sigma_{R}^{\infty} (E \comp)_+$, where $ E \comp$ is an $ \comp$--cell
complex such that $ (E \comp)^H$ is contractible whenever $H \in
\comp$, and empty otherwise \cite[IV.6]{mma}.

\begin{lem} Assume $G$ is a discrete group.
If $X$ is an arbitrary  $G$--cell complex spectrum, then
$X \wedge \Sigma_{R}^{\infty} (E \comp)_+$ is a cofibrant replacement of $X$.
\end{lem}
\begin{proof} Note that $G/ J_+ \wedge G/H_+ $ is an
$\comp$--cell complex, whenever $H \in \comp$ and $J$ is an arbitrary
subgroup of $G$. The collapse map $(E\comp )_+ \rarr S^0$ induces an
$\comp $--$ \{ \uu_H \} $--equivalence $ X \wedge (E\comp )_+ \rarr X $
\end{proof}

\begin{lem} If $G$ has no compact Lie subgroups besides $ 1 $
(e.g.~torsion-free discrete groups), then the stable $ \{ \uu_1
\}$--model structure $\mm$ is the stable $\{ 1 \} $--model structure.
\end{lem}

\begin{lem} If $G$ is a compact Lie group, then the
stable $ \{ \uu_G |H \}$--model structure on $\mm$ is Quillen
equivalent to the $ \{ \text{all}\}$--$ \uu_G $--model structure on
$\mm$.
\end{lem}

Let $J$ be a   subgroup of $G$. Let $R$ be a
monoid in the category of orthogonal $G$--spectra. Let $\mm_G$ denote
the category of $R$--modules in the category of orthogonal
$G$--spectra indexed on $\uuu$, and let $\mm_J$ denote the category of
$R | J $--modules in the category of orthogonal $J$--spectra indexed
on $\uuu $.  Let $\{ \uu_H \}$
be an $\comp $--universe, and let $ \{ \uu_H \}_{ H \in
\mathcal{F}_J}$   be the $ \mathcal{F}_J$--universe.

Note that the condition in the next Lemma is trivially satisfied if
$G$ is a discrete group.

\begin{lem}  Assume that  $(G/K_+ )| J$ has the structure of a
 $J$--$\mathcal{F}_{J}$--cell complex for
any $K \in \comp$. Give $\mm_G$ the stable $\comp $--$ \{ \uu_H \}$--model
structure, and give $\mm_J$ the stable $\mathcal{F}_J $--$ \{ \uu_H \}_{H \in
\mathcal{F}_J }$--model structure. Then the functor
\[F_{J} ( G_+ , -)  \col \mm_{J} \rarr \mm_{G} \]
is right Quillen adjoint to the restriction functor
\[ \mm_{G} \rarr \mm_{J} . \]
\end{lem}

\begin{proof} The restriction functor from
 $G$--spectra to $J$--spectra respects weak
equivalences by the definition of weak equivalences. Since $G/K_+$ is
a $J$--$\mathcal{F}_{J}$--cell complex for all $K \in \comp$, by our
assumption, the relative $G $--$ \comp$--cell complexes are also
relative $ J $--$ \mathcal{F}_{J}$--cell complexes.
\end{proof}

\begin{lem} \label{lem:digchangeofgroups}
 Assume $G$ is a discrete group.
Give $\mm_G$ the stable $\comp $--$ \{ \uu_H \}$--model structure, and
give $\mm_J$ the stable $\mathcal{F}_J $--$ \{ \uu_H \}_{H \in
\mathcal{F}_J }$--model structure. Then the functor
\[G_+ \wedge_{J}   - \, \col \mm_{J} \rarr \mm_{G} \]
is left Quillen adjoint to the restriction functor
\[ \mm_{G} \rarr \mm_{J} . \]
\end{lem}

\begin{proof} Since $G_+ \wedge_{J} J / K_+ \cong G/K_+$ and the
functor $G_+ \wedge_{J} -$ respects
colimits, it follows that $G_+ \wedge_{J} -$ respects cofibrations.

Let $f \col X \rarr Y$ be a $J$--$\mathcal{F}_{J}$--equivalence. We
observe that $ G_+ \wedge_J X $ is isomorphic to
\[ \textstyle\bigvee_{HgJ \in H \backslash G / J} \, H_+
\wedge_{ g J g^{-1} \cap H }
g X \] as an orthogonal $H$--spectrum. Hence the map $G_+ \wedge_J f
$ is an $\comp $--weak equivalence if each $ H_+ \wedge_{ g J g^{-1}
\cap H } g (f)$ is an $H $--equivalence for $H \in \mathcal{F}_G$.
This follows from \cite[V.1.7,V.2.3]{mma} since $ g (f)$ is a
$K$--equivalence for every $ K \leq g J g^{-1} \cap H $, because $
K\in \comp $ and $ K\leq g J g^{-1}$ implies that $ K\in
\mathcal{F}_{g J g^{-1}} $.
\end{proof}

\begin{lem} \label{lem:changeofgroupsadjoint} Assume $G$ is a discrete group.
Let $X$ be a cofibrant object in $\mm_J$ and
let $Y$ be a fibrant object in $\mm_G$. Then \[
 [G_+ \wedge_{J}  X , Y ]_G \cong  [X , (Y | J)  ]_{J} ,\]
where the first
hom-group is in the homotopy category of the $\{ \uu_H \}$--model
structure on $\mm_R$, and the second hom-group is in the
homotopy category of the $\{ \uu_H \}_{H \in \mathcal{F}_J}$--model
structure on $\mm_J$.
\end{lem}
In particular, if $X$ and $Y$ are $G$--spectra and $H \in \comp$,
then \[ [ G / H_+ \wedge X , Y ]_G \cong [X , Y ]_{H} . \]

\begin{rem} A better understanding of the fibrations
 would be useful. They are completely understood when
$G$ is a compact Lie group \cite[III.4.7,4.12]{mma}. Calculations in
the stable $ \{ \uu_H \} $--homotopy theory reduces to calculations
in the stable homotopy categories for the compact Lie subgroups of
$G$ (via a spectral sequence). This follows from Lemma
\ref{lem:changeofgroupsadjoint} using a cell filtration of a cofibrant
replacement of the source by a $\comp$--cell complex.
\end{rem}

\section{Homotopy classes of maps between suspension spectra}

We first give a concrete description of the set of morphisms between
suspension spectra in the $\ww$--stable homotopy category on $\mm_R$.
Then we prove some results about vanishing of negative stems. These results are needed
Section \ref{sec:t-modelstructures}.

Recall that {$\ww \ttt$} denotes $G \ttt$ with the $\ww$--model
structure.
\begin{lem} \label{lem:suspensionhomotopygroups}
Let $X$ and $Y$ be two based $G$--spaces. Then there is a natural
isomorphism
\[ Ho (\ww \mm_R) ( \Sigma^{\infty}_R X ,\Sigma^{\infty}_R  Y )
\cong Ho (\ww \ttt ) ( X , \hocolim_W \, \Omega^W (R (W) \wedge Y )_{\text{cof}} )
\] where the subscript, $\text{cof}$, indicates a $\ww$--cofibrant replacement.
\end{lem}
\begin{proof}
Recall the description of $ \Sigma^{\infty}_R$ in
\ref{eq:R-suspension}. The functors $\Sigma^{\infty}_R $ and
$\Omega^{\infty}_0$ are a Quillen adjoint pair. The result
follows by replacing $ X$  by a cofibrant object, $X_{\text{cof}}$,
in $\ww \ttt$, and  replacing $\Sigma^{\infty}_R Y$ by a cofibrant and
fibrant object as in \ref{fibrantreplacement}.
\end{proof}

\begin{cor} \label{cor:suspensionhomotopygroups}
Let $X$ and $Y$ be two based $G$--spaces. Then there is a natural
isomorphism
\[ Ho (\ww \mm_{\sss}) ( \Sigma^{\infty} X ,\Sigma^{\infty} Y )
\cong Ho (\ww \ttt ) ( X , \hocolim_W \, \Omega^W \Sigma^W Y_{\text{cof}} ) . \]
\end{cor}
In particular, if $X $ is a finite $\ww$--cell complex, then
\[ Ho (\ww \mm_{\sss}) ( \Sigma^{\infty} X ,\Sigma^{\infty} Y )
 \cong \colim_W \, Ho (\ww \ttt ) ( X
\wedge S^W , Y \wedge S^W) . \]
We next show
that the negative stable stems
are zero.  In what follows homotopy means usual
homotopy (a path in the space of maps).
\begin{lem} \label{mapbetweengenerators}  Let $V$ be a finite
dimensional real $G$--representation with $G$--action factoring
through a Lie group quotient of $G$. Let $X$ be a based $G$--space,
and let $n > 0$ be an integer. Then any based $G$--map
\[   S^V \rarr  S^V \wedge X \wedge S^n \] is $G$--null-homotopic.
\end{lem}
\begin{proof}  Assume  the action on $S^V$ factors through  a
compact Lie group quotient $G/K$. The problem reduces to show that $
S^V \rarr S^V \wedge X^K \wedge S^n $ is $G/K$--null homotopic for
all $n > 0$. Hence we can  assume that $G$ is a compact Lie group.
By Illman's triangulation theorem  $S^V$ is a finite $G$--cell
complex \cite{ill}. We choose a $G$--CW structure on $S^V$.
Let $(G /H_i \times D^{n_i} )_+ $ be a cell of $S^V$. We compare the
real manifold dimensions, denoted $\dim$, of the $H_i$--fixed points of
$S^{V}$ and the cells in $S^V$. This gives that $n_i = \dim( V^{H_i}
) - \dim ( N_G H_i / H_i )$. To prove the Lemma it suffices to show
that any given map $ f \col S^V \rarr S^V \wedge X \wedge S^n $
extends over the cone $S^V \wedge I$ of $S^V$. There is a sequence
\[ S^V = Y_{-1} \rarr Y_0 \rarr Y_1 \rarr \cdots \rarr Y_N =
S^V \wedge I \] where $Y_{n+1} $ is obtained from $Y_n$ by a pushout
\[ \xymatrix{ \midvee {G/H_i}_+ \wedge S^{n_i}
\ar[d] \ar[r] & Y_n \ar[d] \\
\midvee {G/H_i}_+ \wedge D^{n_i +1} \ar[r] & Y_{n+1} } \] where the
wedge sum is over all $i$ such that $n_i =n$, and $N$ satisfies $n_i
\leq N$ for all $i$. Hence it suffices to show that any map $ \vee
{G/H_i}_+ \wedge S^{n_i} \rarr S^V \wedge X \wedge S^n$ is $G$--null
homotopic for all $i$. This is equivalent to showing that $S^{n_i}
\rarr S^{V^{H_i}} \wedge X^{H_i} \wedge S^n$ is null homotopic,
which is true because $ n_i = \dim ( V^{H_i} ) - \dim ( N_G H_i /
H_i ) < \dim (V^{H_i}) +n$.
\end{proof}
\begin{lem} \label{stablemapbetweengenerators}
Let $\uu$ be any $G$--universe, and let $\stable (\uu) $ be a  $\ww$--\good\
collection of subgroups of $G$.
Then we have that \[  \mathrm{Ho} ( \ww \mm_{\sss})  (
\Sigma^{\infty} G/H_+ , \, \Sigma^{\infty} G/K_+ \wedge S^n) =0 , \]
for all $H ,K \in \ww$ and $n >0$.
\end{lem}
\begin{proof}  The space $S^V \wedge G/H_+ $ is homeomorphic to
a finite $\ww$--cell complex by Lemma \ref{productcomplexes} and
compactness of $S^V \wedge G/H_+ $. Corollary
\ref{cor:suspensionhomotopygroups} gives that the group $\text{Ho} (
\ww \mm_{\sss}) ( \Sigma^{\infty} G/H_+ , \, \Sigma^{\infty} G/K_+ \wedge
S^n)$ is isomorphic to
\[ \colim_{ V \in \vv (\uu )} \text{Ho} ( \ww  \ttt )( S^V \wedge
G/H_+ , \, S^V \wedge G/K_+ \wedge S^n ) . \]
It suffices to show that any map $ S^V \wedge G/H_+ \rarr S^V \wedge
G/K_+ \wedge S^n $ is $G$--null homotopic. This is equivalent to
show that $ S^V \rarr S^V \wedge (G/K)_{+} \wedge S^n$ is $H$--null
homotopic. This follows from Lemma \ref{mapbetweengenerators}.
\end{proof}

Lemma \ref{stablemapbetweengenerators} allows the
formation of $\ww$--CW complex approximations.
\begin{lem} \label{cellapproximation}
Let  $\stable (\uu) $ be a  $\ww$--\good\  collection  of
subgroups of $G$. Let $T$ be an $\sss $--module such that $ \pi^{H}_j T =0$,
for $j <n$ and $H \in \ww$. Then there is a cell complex, $T'$,
built out of cells of the form $ \Sigma^{\infty}_{\mathbb{R}^{k'}}
S^{k-1} \wedge G/H_+ \rarr \Sigma^{\infty}_{\mathbb{R}^{k'}} D^{k}
\wedge G/H_+ $, for $k - k' \geq n $ and $ H \in \ww$, and a
$\ww$--weak equivalence $T' \rarr T$.
\end{lem} \begin{proof} The approximation can be constructed
as a $\ww$--cell complex using Lemmas \ref{corep} and \ref{stablemapbetweengenerators}.
\end{proof}
\begin{lem} \label{respecttensor}
 Let $\ww$ be an \good\  collection  of
subgroups of $G$, and assume that  $\stable (\uu) $ is  $\ww$--\good.
 Let $R$ and $T$ be two $\sss $--modules. If $\pi^{H}_i R =
0 $, for $ i <m $ and $H \in \ww$, and $ \pi^{H}_j T =0$, for $j <n$ and
$H \in \ww$, then $\pi^{H}_k (R \wedge T) = 0 $ for $ k < m +n$ and $H
\in \ww$.
\end{lem}
\begin{proof} Replace $R$ and $T$ by $\ww I$--cell complexes
made of cells in dimension greater or equal to $m$ and $n$. This is possible by
Lemma \ref{cellapproximation}. The spectrum analogue
of Lemma \ref{productcomplexes} gives that $R \wedge T$ is again a
$\ww$--cell complex made out of cells in dimension greater or equal
to $m +n$. The result now follows from Lemma
\ref{stablemapbetweengenerators}.
\end{proof}

\begin{pro} \label{AB} Let $\ww$ be an \good\ collection  of
subgroups of $G$ and assume that  $\stable (\uu) $ is   $\ww$--\good.
Let $R$ be a ring spectrum such that $\pi^{H}_n R= 0$ for all
$n <0$ and $H \in \ww$.  Then we have that
\[ \text{Ho} ( \ww \mm_R)
 (\Sigma_{R}^{\infty}  G/H_+ , \, \Sigma_{R}^{\infty} G/K_+
 \wedge S^n )  =0 , \] for all $H , K \in \ww$ and
$n >0$. \end{pro} \begin{proof} The group $\text{Ho} ( \ww \mm_R)
(\Sigma_{R}^{\infty} G/H_+ , \, \Sigma_{R}^{\infty} G/K_+ \wedge S^n
)$ is isomorphic to
\[ \text{Ho} ( \ww \mm_{\sss} ) (\Sigma^{\infty} G/H_+ ,
\Sigma^{\infty} G/K_+ \wedge R \wedge S^n ) . \] The result
follows from Lemma \ref{mapbetweengenerators}
and Proposition \ref{respecttensor}. (Let $T $ be
$ \Sigma^{\infty} G/ K_+ \wedge S^n$.)
\end{proof}

\begin{lem} \label{Rcellapproximation} Let $\ww$ be an \good\ collection  of
subgroups of $G$ and assume that  $\stable (\uu) $ is   $\ww$--\good. Let $R$
be a ring such that $\pi^{H}_n R= 0$, whenever $n < 0$ and $H \in \ww$.
  Let $T$ be an $R$--module such that $ \pi^{H}_j T=0$, for $j <n$ and
$H \in \ww$. Then there is a cell complex, $T'$, built out of cells
of the form $ \Sigma^{R}_{\mathbb{R}^{k'}}
S^{k-1} \wedge G/H_+ \rarr \Sigma^{R}_{\mathbb{R}^{k'}} D^{k}
\wedge G/H_+ $, for $k - k' \geq n $ and $ H \in \ww$, and a $\ww$--weak
equivalence $T' \rarr T$. \end{lem}
\begin{proof} This follows from Lemma \ref{AB} and the proof of
Lemma \ref{cellapproximation}.
\end{proof}
If the universe $\uu$ is trivial and $H$ is a not subconjugated to
$K$ in $G$, then there are no nontrivial maps from $ \Sigma^{\infty}
G/ H_+ \wedge S^n $ to $ \Sigma^{\infty} G/ K_+ \wedge S^m $. We
take advantage of this to strengthen Lemma \ref{AB}.
\begin{pro} \label{Abtrivial} Let  $\uu$ be a  trivial universe.
Let $\ww$ be an \good\ collection of subgroups of $G$. Let $R$ be
ring spectrum such that $\pi^{H}_n R= 0$ for all $n <0$ and $H \in
\ww$.  Then, for each $H, K$ in $\ww$, we have that
\[ \text{Ho} ( \ww \mm_R)
 (\Sigma_{R}^{\infty}  G/H_+ , \, \Sigma_{R}^{\infty}
 G/K_+ \wedge S^n )  =0 , \]
whenever $n >0$ or $H $ is not subconjugated to $ K $.
\end{pro} \begin{proof}
If $n >0$, then the result follows from Proposition
\ref{AB}. If $H $ is not subconjugated to $ K $, then
\[  \text{Ho} (\ww \mm_R)
 (\Sigma_{R}^{\infty}  G/H_+ , \,
 \Sigma_{R}^{\infty} G/K_+ \wedge S^n ) \cong
 \colim_m \, \pi_0 \, \Omega^m ( G/K_+ \wedge
 (R \wedge S^n)( \rr^m  ))^H .\]
This is 0 since $ G/K^H_+$ is the basepoint.
\end{proof}

\subsection{The Segal--tom Dieck splitting theorem}
We consider homotopy groups of suspension spectra.
Let $G$ be a compact Hausdorff group, let $\uu$ be a complete $G$--universe,
and let $\mm_{\sss}$ have the $\overline{\lie (G)}$--model structure.
\begin{pro} \label{pro:tom diecksplitting}
If $Y$ is a $G$--space, then there is an isomorphism of abelian
groups \[ \textstyle\bigoplus_H \, \pi_* (\Sigma^{\infty} ( EW_G H_+ \wedge_{W_G H}
\Sigma^{\text{Ad} (W_G H)} Y^H )) \rarr \pi^{G}_* ( \Sigma^{\infty}_{\sss}
Y )\] where the sum is over all $G$--conjugacy classes of subgroups $H$
in $\overline{\lie (G)}$ and $\text{Ad} (W_G H)$ is the Adjoint representation of $ W_G H$.
\end{pro}
\begin{proof}  We have an isomorphism
\[ \colim_{N \in \lie (G)} \, \colim_{ V \subset  \uu^N} \,\pi_* (
(\Omega^V S^V Y^N )^G ) \rarr \pi^{G}_* (\Sigma^{\infty}_{\sss} Y) \]
where the rightmost colimit is over  indexing representations $V$ in $\uu^N$.
The universe $\uu^N$ is $G/N$--complete. The result follows from the splitting
theorem for compact Lie groups \cite[V.9.1]{lms}.
\end{proof}

If $\uu$ is a complete $G$--universe, then $\uu$ restricted to $K$ is
again a complete $K$--universe \cite[Section 3]{art}. So for any $K \in
\overline{\lie (G)}$, the stable homotopy groups at $K$ of a $G$--space
$Y$ calculated in the $G$--homotopy category are isomorphic to those
calculated in the $K$--homotopy category. Hence the calculation of
the $n$-th stable  homotopy group at $K \in \overline{\lie (G)}$ of a
$G$--space $Y$ reduces to Proposition \ref{pro:tom diecksplitting}
(with $G$ replaced by $K$).

\subsection{Self-maps of the unit object}

The additive tensor category, $\text{Ho} ( \ww \mm_R)$, is naturally
enriched in the category of modules over the ring \[Ho (\ww \mm_R )
( \Sigma_{R}^{\infty} S^0 , \Sigma_{R}^{\infty} S^0 ) \] of self
map of the unit object  in the homotopy
category. Let us denote this ring
by \mdfn{$B^{R}_{\ww}$}, and denote $B^{S^0}_{\ww}$ simply by $B_{\ww}
 $. The ring $B^{R}_{\ww}$ depends on $G$, $\ww$, $R$,
and the $G$--universe $\uu$. If $G \in \ww$, then we can
identify $ B^{R}_{\ww}$ with $\pi^{G}_0 (R) $.

If $A$ is an $R$--algebra, then $ B^{A}_{\ww}$ is an $
B^{R}_{\ww}$--algebra. Since all algebras of orthogonal spectra are
$\sss $--algebras, it is important to understand the ring $ B_{\ww}$.

If $G$ is a compact Lie group, the universe is complete,
and $\ww$ is the collection  of all closed
subgroups of $G$, then $B_{\ww} $ is naturally isomorphic
to the Burnside ring, $A(G)$, of $G$ \cite[XVII.2.1]{ala}.
\begin{lem} Let $\uu$ be a complete $G$--universe and let $\ww $ be
the collection $ \overline{\lie (G)}$. Then the self-maps of $\Sigma^{\infty} S^0$,
in the homotopy category of $\ww \mm$, is naturally isomorphic to
\[ \colim_{ U \in \lie (G) } \, A (G/U), \] where
$A(G/U) \cong \Ho ( \ww \mm ) ( \Sigma^{\infty} G/U_+ , \Sigma^{\infty} S^0 )$
is the Burnside ring of the Lie group $G/U$ and
the maps in the colimit are induced by the quotient maps $ G/U_+
\rarr G/V_+$, for $U < V$ in $\lie (G)$. \end{lem}

In general, it is difficult to determine $B_{\ww}$. For example,
when $G$ is a finite group and $\ww$ is a family, then the proof of
the Segal conjecture gives that the ring $\bb_{\ww}$ is isomorphic
to the Burnside ring $ A(G)$ of $G$ completed at the augmentation
ideal \[\cap_{H \in \ww} \ker ( A(G) \rarr A (H )) , \] where the maps
$ A(G) \rarr A (H)$ are the restriction maps \cite[XX.2.5]{ala}.

We give an elementary observation which shows that different
collections $\ww$ might give rise to isomorphic rings $\bb_{\ww}$.
\begin{lem} \label{lem:different}
Let $N$ be a normal subgroup of a finite group
$G$, and let $\ww_N$ be the family of all subgroups contained in
$N$. If $X \in G \ttt$ has a trivial $G$--action and $Y \in G \ttt$,
then
\[ \Ho (\{ N \} \mm_R ) ( \Sigma^{\infty} X , \Sigma^{\infty}  Y )
\rarr \Ho ( \ww_N \mm_R ) (\Sigma^{\infty}  X ,
 \Sigma^{\infty}  Y ) \] is an isomorphism.
\end{lem}
In particular, $B_{ \{ N \} } $ is isomorphic to $B_{\ww_N}$.
\begin{proof}  Let $X'$ be a cell complex
replacement of $X$ built out of cells with trivial $G$--actions.
The space $E\ww_N$ is $\overline{\lie (G)}$--equivalent to $E (G/N)$.
This is an $\{ N \}$--cell complex. Hence $ X' \wedge EG/N_+ \rarr
X'$ is a cofibrant replacement of $X$ both in the $\{ N \}$ and and
in the $\ww_N$--model categories. A $\ww_N$--fibrant replacement
$Y'$ of $Y $ is also an $\{ N \}$--fibrant replacement.
\end{proof}
\begin{rem}
If $X$ does not have trivial $G$--action, then the $ \Ho (\ww \mm ) (
 \Sigma^{\infty} X  , \Sigma^{\infty} Y )$   are typically different for
the collections $ \{ N \} $ and $\ww_N$.
\end{rem}

\section{Postnikov t-model structures}
\label{sec:t-modelstructures} We modify the construction of the
$\ww$--$\cc$--model structure on $\mm_R$ by considering the
$n$--$ \ww$--equivalences, for all $n$, instead of $\ww$--equivalences.
This is used when we give model structures to the category of
pro--spectra, $ \p \mm_R$, in Section \ref{sec:postmodelpro}.

The homotopy category of a stable model category is a triangulated
category \cite[7.1]{hov}.
We consider t-structures on this triangulated category
together with a lift of the t-structures to the
model category itself. The relationship between $n$--equivalences and
t-structures is given below in Definition \ref{neq} and Lemma
\ref{lem:t-structure}.

\subsection{Preliminaries on t-model categories}

We  recall the terminology of a t-structure \cite[1.3.1]{bbd}
and of a t-model structure \cite[4.1]{tfi}.
\begin{defn}
A homologically graded \dfn{t-structure}
on a triangulated category $\dd$, with shift functor $ \Sigma$,
consists of two full subcategories $\dd_{\geq  0}$ and $\dd_{\leq
 0}$ of $\dd$, subjected to the following three axioms:
\begin{enumerate}
\item   $\dd_{\geq  0} $ is closed under $\Sigma$, and
$\dd_{\leq 0} $ is closed under $\Sigma^{-1}$;
\item for every object $X$ in $\dd$, there is a distinguished
triangle \[ X' \rarr X \rarr X'' \rarr \Sigma X' \] such that $X'
\in \dd_{\geq 0} $ and $X'' \in \Sigma^{-1} \dd_{\leq 0}$; and
\item $\dd ( X ,Y ) = 0$,
whenever $ X \in \dd_{\geq 0} $ and $Y \in \Sigma^{-1} \dd_{\leq
0}$.
\end{enumerate}
For convenience we  also assume   that $\dd_{\geq 0}$ and
$\dd_{\leq 0}$ are  closed under isomorphisms in $\dd$.
\end{defn}
\begin{defn} Let $ \dd_{\geq  n} = \Sigma^{n} \dd_{\geq  0}$, and
let $ \dd_{\leq n} = \Sigma^{n} \dd_{\leq 0}$. \end{defn}

\begin{rem} A  homologically graded t-structure $( \dd_{\geq 0} ,
\dd_{\leq 0} )$ corresponds to a cohomologically graded t-structure
$( \dd^{\leq 0} , \dd^{\geq 0} )$ as follows: $ \dd^{\geq n} =
\dd_{\leq -n} $ and $ \dd^{\leq n} = \dd_{\geq -n} $.
\end{rem}

\begin{defn} \label{neq} The class of
\mdfn{$n$--equivalences} in $\dd$, denoted $W_n$, consists of all
maps $f \col X \rarr Y$ such that there is a triangle
\[ F \rarr X \stackrel{f}{\rarr} Y \rarr \Sigma F \] with
$F \in \dd_{ \geq n} $. The class of \mdfn{co$-n$--equivalences} in
$\dd$, denoted co$W_n$,  consists of  all maps $f$ such that there is
triangle \[ X \stackrel{f}{\rarr} Y \rarr C \rarr \Sigma X \] with
$ C \in \dd_{\leq   n}$.
\end{defn}

If $\dd$ is the homotopy category of a stable model category $\kk$,
then a map $f$ in $\kk$ is called a (co--)$n$--equivalence if the
corresponding map $f$ in the homotopy category, $\dd$, is a
(co--)$n$--equivalence. We use the same symbols $W_n$ and co$W_n$ for the
classes of $n$--equivalences and co--$ n$--equivalence in $\kk$ and
$\dd$, respectively.

\begin{defn} A  t-model
category  is a proper simplicial stable model category $\kk$
equipped with a t-structure on its homotopy category together with a
functorial factorization of maps in $\kk$ as an $n$--equivalence
followed by a co--$n$--equivalence in $\kk$, for every integer $n$.
\end{defn}

T-model categories are discussed in detail in \cite{tfi}. They give
rise to interesting model structures on pro--categories.

\subsection{The $d$--Postnikov t-model structure on $\mm_R$}
We  construct a t-model structure on $\mm_R$ such that the
t-structure on the homotopy category of $\mm_R$ is given by
Postnikov sections. In Section
\ref{pro--structure} we use this t-model structure to produce model
structures on the category of pro--spectra. We allow Postnikov
sections where the cut-off degree of $\pi_{*}^H $ depends on $H$.
See also \cite{lew}. \begin{con} \label{con:t-structure} Assume that $\dd$ is the
homotopy category of a proper simplicial stable model category
$\mm$.  Let $\dd_{\geq 0
}$ be a strictly full subcategory of $\dd$ that is closed under
$\Sigma$. Define $\dd_{\geq n} $ to be $\Sigma^{n} \dd_{\geq 0}$.
Let  $W_n$ be as in Definition \ref{neq}, and lift $W_n$ to $\mm$.
Let $C$ denote the class of cofibrations in $ \mm$, and define $C_n
= W_n \cap C$ and $F_n = \text{inj } C_n$, the class of maps with the
right lifting property with respect to $C_n$. Let $\dd_{\leq n- 1 }$ be
the full subcategory of $\dd$ with objects isomorphic to $\fib (g)$
for all $g \in F_n$.
If there is a functorial factorization of any map
in $ \mm$ as a map in $C_n$ followed by a map in $F_n$, then
$\dd_{\geq 0}$, $\dd_{\leq 0 }$ is a t-structure on $\dd$. Hence the
model category $\mm$, the factorization, and the t-structure
on $\dd$ is a t-model structure on $\mm$ \cite[4.12]{tfi}.
\end{con}

Let $\cc$ and $\ww$ be  collections of subgroups of $G$
such that  $\cc \ww \subset \cc$.  Let $R$ be a ring, and
let $\dd$ be the homotopy category of $\ww \cc \mm_R$.
\begin{defn} A  \mdfn{class
function on $\ww$} is a function $d \col \ww \rarr \zz\cup \{-
\infty, \infty \} $ such that $ d (H) = d (g H g^{-1})$, for all $ H
\in \ww$ and $g \in G$.\end{defn}

\begin{defn} Let $d$ be a class function. Define a full subcategory
of $\dd$ by \[  \dd_{\geq  0}^d = \{ X \ | \ \Pi^{U}_i ( X) =0 \text{ for }
i < d (U) , \ U \in \ww\} . \] Let $ \dd_{\leq -1}^d $ be the full
subcategory of $\dd$ given by Construction \ref{con:t-structure}.
\end{defn}

The next result is needed to get a t-model structure on $\mm_R$.

\begin{lem}
\label{ncon} Any map in $\mm_R$ factors functorially as a map in
$C_n$ followed by a map in $F_n$. Moreover, there is a canonical map
from the $n$-th factorization to the $(n-1)$-th factorization.
\end{lem}
\begin{proof}
The proof is similar to the proof of Theorem \ref{model2}. See also
\cite[Appendix]{fch}.
\end{proof}
\begin{lem}
\label{lem:t-structure} Let $d$ be a class function on $\ww$.
The $\ww$--$ \cc $--model structure on $\mm_R$,
the two classes $\dd_{\geq  0 }^d $ and
$\dd_{\leq 0}^d $, together with the factorization in
Lemma \ref{ncon} is a t-model structure.
\end{lem}
\begin{proof}
This follows from Theorem \ref{model2} and Construction
\ref{con:t-structure} \cite[4.12]{tfi}.
\end{proof}
This t-model structure is called the \mdfn{ $d$--Postnikov t-model
structure } on $\ww \cc \mm_R$. We call the $0$--Postnikov t-model
structure simply the Postnikov t-model structure.

A map $f$ of spectra is an $n$--equivalence with respect to
the $d$--Postnikov  t-structure, as defined in Definition \ref{neq},
if and only if  $ \Pi^{U}_m (f) $ is an isomorphism for $m < d(U) +n $
and $\Pi^{U}_{d(U)+n} (f)$ is surjective for all $U \in \ww$.

\subsection{An example: Greenlees connective K-theory}

To show that there is some merit to the generality of $d$--Postnikov
t-structures, we recover Greenlees equivariant connective $K$--theory as
the $d$--connective cover of equivariant $K$--theory for
a suitable class function $d$. Let $G$ be a compact Lie group, and
let $\ww = \cc$ be the collection of all closed subgroups of $G$.
Let $P_n $ denote the $n$-th Postnikov section
functor, and let $C_n $ denote the $n$-th connective cover functor.

\begin{lem} \label{Greenlees}
Let $G$ be a compact Lie group. Let $d$ be the class function such
that $d (1   ) = 0$ and $ d ( H) = - \infty $ for all $ H \not=
1 $. Then \[ X_{\leq  n  }  = F ( EG_+ , P_{ n } X )  \]
is a functorial truncation functor for the $d$--Postnikov t-model
structure on $\ww \mm_R$.
\end{lem}

The $n$-th $d$--connective cover, $X_{\geq n}$, of $X$ is such that
the left most square in the following diagram is a homotopy pullback square
\begin{equation} \label{eq:Gktheory} \xymatrix{ X_{\geq n}
\ar[r] \ar[d] & X \ar[r] \ar[d] & F ( EG_+ , P_{ n -1} X ) \ar@{=}[d] \\
F ( EG_+ , C_{ n }X ) \ar[r] & F ( EG_+ ,X ) \ar[r] & F ( EG_+ , P_{
n -1} X ) .}\end{equation} In particular, $ (K_G )_{\geq 0}$ is
Greenlees' equivariant connective $K$--theory \cite[3.1]{gre}.

\begin{proof} Axiom 1 of a t-structure is satisfied since \[
\Sigma^{-1} F ( EG_+ , P_{ n } X ) \cong F ( EG_+ , \Sigma^{-1} P_{ n }
X ) \cong F ( EG_+ , P_{ n } (\Sigma^{-1} P_{ n }
X )) . \]
We combine the verification of axioms 2 and 3 of a t-structure. Let
$X_{\geq n}$ denote the homotopy fiber of the natural transformation
$X \rarr F ( EG_+ , P_{ n -1} X )$. Since $ X \rarr F (EG_+ , X)$ is
a non-equivariant equivalence we conclude, using Diagram
\ref{eq:Gktheory}, that $X_{\geq n} $ and $ C_{n} X$ are
non-equivariant equivalent. Hence $ X_{\geq n} \in \dd_{\geq n}$ for
all $X \in \dd$.
If $ Y\in \dd_{\geq n}$ and $X\in \dd$, then \[\dd ( Y, F ( EG_+ ,
P_{ n -1} X)) = 0\] since $ Y\wedge EG_+$ is in $ C_{n} \dd$.
\end{proof}

This example can also be extended to arbitrary compact Hausdorff
groups \cite{fg}.
\subsection{Postnikov sections}
Suppose $d$ is a constant function and $R$ has trivial $\cc$--homotopy
groups in negative degrees. Then there is a useful description of
the full subcategory $ \dd_{\leq 0}$ of the homotopy category $\dd$
of $\ww \cc \mm_R$.

\begin{defn} \label{defn:connective}
We say that a spectrum $R$ is \mdfn{$\cc$--connective} if $\Pi^{U}_n
(R) =0$ for all $n<0$ and all $U \in \cc$.
\end{defn}
\begin{pro}  \label{pro:Dgeq1}
Let $R$ be a $\cc$--connective ring.
Then there is a t-structure $( \dd_{\geq  0} , \dd_{\leq  0} )$ on
the homotopy category $\dd$ of $\ww
\cc \mm_R$ such that:
 \[   \dd_{\geq  0} = \{ X \ | \ \Pi^{U}_i ( X) =0
 \ \mathrm{ whenever } \ i
< 0, \, U \in \ww \} \] and \[ \dd_{\leq 0} \subset  \{ X \ | \ \Pi^{U}_i (
X ) = 0 \ \mathrm{ whenever } \ i >0, \, U \in \ww \} . \]
The inclusion is an equivalence whenever $\ww \subset \cc$.
 \end{pro}
\begin{proof}
Recall that an object $Y$ is in  $ \dd_{\leq -1}$ if and only if
$\dd (X ,Y) =0$ for all $X \in \dd_{\geq 0}$
\cite[1.3.4]{bbd}. Proposition \ref{AB} gives that
$\Sigma_{R}^{\infty}  G/H_+ \wedge
S^n \in \dd_{\geq 0}$, for all $ n \geq 0$ and all $H \in \cc$. This
gives that \[ \dd_{\leq -1 } \subset \{ Y \ | \ \Pi^{U}_i ( Y) = 0
\text{ whenever } U\in \ww, i\geq 0\} .\]

We prove the converse inclusion when $\ww \subset \cc$.
Assume that $X \in \dd_{\geq 0}$ and that
$\Pi^{U}_i ( Y) = 0 $, whenever $ U\in \ww$ and $i \geq 0$. By Lemma
\ref{cellapproximation} there is a $\ww$--cell complex approximation
$X' \rarr X$ such that $X'$ is a cell complex built from cells in
non-negative dimensions, and $ X' \rarr X$ is a $\ww$--isomorphism in
non-negative degrees, hence a $\ww$--equivalence.
This implies  that $\dd ( X , Y) =0$. Since $
\dd ( X , Y) =0$   for all $X \in \dd_{\geq 0}$ we
conclude that $Y \in \dd_{\leq -1}$.
\end{proof}

If a map $g$ is a co$-n$--equivalence in the  Postnikov
t-structure, then   $ \Pi^{U}_m (g) $ is an isomorphism for
$m > n  $ and  $ \Pi^{U}_{n } (g )$ is injective for $U\in \ww$.

When the universe is trivial we can give a similar description of
the t-model structure for more general functions $d$. We say that a
class function $d \col \ww \rarr \zz\cup \{- \infty, \infty \} $ is
\mdfn{increasing} if $ d (H) \leq d (K)$ whenever $ H \leq K$.
\begin{pro} \label{gentpostnikov} Assume the $G$--universe $\uu$ is
trivial, and let $R$ be a $\cc$--connective ring. Let $d$ be an
increasing class function.
Then there is a t-structure $( \dd_{\geq  0} , \dd_{\leq  0} )$
on the homotopy category $\dd$ of $\ww
\cc \mm_R$ such that:
\[   \dd_{\geq  0}^d = \{ X \ | \ \Pi^{U}_i ( X) =0
 \ \mathrm{ whenever } \  i
< d (U) , \, U \in \ww \} \] and \[ \dd_{\leq 0}^d \subset  \{ X \ | \
\Pi^{U}_i ( X ) = 0 \ \mathrm{ whenever } \ i> d (U), U \in \ww \} .
\]
\end{pro}
\begin{proof} This follows from Proposition \ref{Abtrivial}, Construction
\ref{con:t-structure} and the proof of Proposition \ref{pro:Dgeq1}.
 \end{proof}

\begin{defn}  Let $X$ be a spectrum in $  \mm_R$.
The  $n$-th Postnikov section  of $X$
is a spectrum $P_n X$ together with a map $p_n X \col X \rarr P_n X$
such that $ \Pi^{U}_m (P_n X) =0 $, for $m > n$ and   $U
\in \ww$, and $\Pi^{U}_m ( p_n X ) \col \Pi^{U}_m (X) \rarr
\Pi^{U}_m ( P_n X )$ is an isomorphism, for  $m \leq n$ and
$U \in \ww$. A \mdfn{Postnikov system} of $X$ consists of a Postnikov
factorization $p_{n} \col X \rarr P_{n} X$, for every  $n \in \zz $,
together with maps $r_{n} X \col  P_{n } X \rarr P_{n-1} X$,
for all $n \in \zz$, such that $ r_{n} X \circ p_{n} X= p_{n-1}X$.
\end{defn}
Dually, one defines the  $n$-th connected cover  $C_n X \rarr X$ of $X$.
The $n$-th connected cover satisfies $\Pi^{U}_k ( C_n X ) = 0 $, for $ k
< n$, and $\Pi^{U}_k ( C_n X ) \rarr \Pi^{U}_k ( X ) $ is an
isomorphism, for $k \geq  n$.

\begin{defn} A  functorial Postnikov system  on $\mm_R$
consists of functors $ P_n $, for each $n \in \zz$, and natural
transformation $ p_n \col 1 \rarr P_n $ and, $ r_n \col P_n \rarr
P_{n-1}$ such that
$p_n (X)$ and $r_n (X)$, for $n \in \zz$, is a Postnikov system for
any spectrum $X$.
\end{defn}
\begin{pro}  Let $R$ be a  $\cc$--connective ring. Then the
category $\ww \cc \mm_R$ has a functorial Postnikov system. \end{pro}
\begin{proof} This follows from Lemma \ref{ncon}. \end{proof}

\begin{rem} It is often  required that the maps $r_n X$ are
fibrations for every $n$ and $X$. We can construct a functorial
Postnikov tower with this property if we restricted ourself to the
full subcategory $ \dd_{\geq n}$ for some $n$ \cite[Section 7]{tfi}.
\end{rem}

\subsection{Coefficient systems} \label{sec:coeff} In this
subsection we describe the Eilenberg--Mac\,Lane objects in the
Postnikov t-structure. Let $\cc $ be an \good\ collection of
subgroups of $G$ and assume that $\stable (\uu) $ is a $\ww $--\good\
collection. Let $\ww$ be a collection such that $ \cc \ww \subset \cc$.
Let $R$ be a $\cc$--connective ring spectrum.

\begin{defn}
 The \mdfn{heart} of a t-structure $( \dd_{\geq  0} , \dd_{\leq   0} )$
on a triangulated category $\dd$ is the full subcategory $ \dd_{\geq
0} \cap \dd_{\leq 0}$ of $\dd$ consisting of objects that are
isomorphic to objects both in $ \dd_{\geq 0}$ and in $\dd_{\leq 0}$.
\end{defn}  The heart of a t-structure is an abelian
category \cite[1.3.6]{bbd}.

\begin{defn}  An $R$--module $X$ is said to
be an  Eilenberg--Mac\,Lane spectrum  if $\Pi^{U}_n ( X) = 0
$, for all $n \not= 0$ and all $U \in \ww$.
\end{defn}
\begin{lem} Let  $d \col \ww \rarr \zz$ be  the $ 0$--function. If  $\cc$ and $\ww$ are  collections of subgroup of $G$ such that $
\cc \ww \subset \cc $, then  the heart of the homotopy category of $\ww \cc \mm_R$
is contained in the full subcategory consisting of the Eilenberg--Mac\,Lane
spectra. If $\ww \subset \cc$, then the heart is exactly the full
subcategory consisting of the Eilenberg--Mac\,Lane spectra.
\end{lem}
\begin{proof} This follows from Proposition  \ref{pro:Dgeq1}.
\end{proof}

We give a more algebraic description of the heart in
terms of coefficient systems when $ \ww \subset \cc$. Let $\dd $
denote the homotopy category of $\ww \cc \mm_R$.

\begin{defn}  The \mdfn{orbit category,  $\oo$,}
is the full subcategory of $\dd$ with objects $ \Sigma^{\infty}_R
G/H_+$, for $H \in \ww$.
\end{defn}
The orbit category depends on $G$, $\ww ,\cc$, and the $G$--universe $\uu$.
\begin{defn} \label{def:coefficients}
A \mdfn{$\ww$--$ R$--coefficient system} is a contravariant additive
functor from $ \oo^{\op} $ to the category of abelian groups.
\end{defn}

Denote the category of $\ww$--$ R$--coefficient systems by
\mdfn{$\coef$}. This is an abelian category. An object $Y$ in $\dd$
naturally represents a coefficient system given by \[
\Sigma^{\infty}_R G/H_+ \mapsto \dd ( \Sigma^{\infty}_R G/H_+ , Y ) . \]

\begin{defn}
Let $X$ be an $R$--module spectrum. The $n$-th homotopy coefficient system of
$X$, \mdfn{$\pi_{n}^{\ww} (X)$}, is the coefficient system naturally
represented by $ X \wedge \Sigma_{\mathbb{R}^{n}}^R S^0$, for  $n \geq 0$,
 and by $ X \wedge \Sigma_{R}^{\infty} S^{-n} $, for $n \leq 0 $.
\end{defn}
\begin{lem} \label{projective}
There is a natural isomorphism \[ \coef ( \pi_{0}^{\ww} (
\Sigma_{R}^{\infty} G/H_+ ) , M ) \cong M (G /H_+ ) \] for any $H \in
\ww$ and any coefficient system $M \in \coef$.
\end{lem}
\begin{proof}
This is a consequence of the Yoneda Lemma.
\end{proof}

\begin{pro} \label{postnikovhearth} Let $R$ be a $\cc$--connective
ring spectrum. The functor $\pi_{0}^{\ww} $ induces a natural
equivalence from the full subcategory of Eilenberg--Mac\,Lane
spectra in the homotopy category of $\ww \cc \mm_R$, to the category
of $\ww$--$ R$--coefficient systems.
\end{pro}
\begin{proof}
We show that for every coefficient system $M$,
there is a spectrum $HM$ such that $ \pi^{\ww}_0 ( HM )$ is
isomorphic to $M$ as a coefficient system, and furthermore, that
$\pi^{\ww}_0 $ induces an isomorphism $ \dd ( HM , HN ) \rarr
\coef^{\ww} (M , N)$ of abelian groups.

We construct a functor, $H$, from $\coef$ to the homotopy category
of spectra. The natural isomorphism in Lemma \ref{projective} gives
a surjective map of coefficient systems
\[ f(M) \col
\textstyle\bigoplus_{H \in \ww} \textstyle\bigoplus_{M(G/H_+)} \,
\pi^{\ww}_0 ( G/H_+ ) \rarr M .\] This construction is natural in $M$. Let
$C_M$ be the kernel of $f (M)$ and repeat the construction with
$C_M$ in place of $M$. We get an exact sequence \begin{equation}
\label{shortex} \textstyle\bigoplus_{K \in \ww}
\textstyle\bigoplus_{C_M (G/K_+)} \pi^{\ww}_0 ( G/K_+ ) \rarr
\textstyle\bigoplus_{H \in \ww} \textstyle\bigoplus_{M(G/H_+)}
\pi^{\ww}_0 ( G/H_+ ) \rarr M \rarr 0.\end{equation} This sequence
is natural in $M$. We have that $ \textstyle\bigoplus_{H \in \ww}
\textstyle\bigoplus_{M(G/H_+)} \pi^{\ww}_0 ( G/H_+ )$ is naturally isomorphic to
\[ \pi^{\ww}_0 ( \midvee_{H \in \ww} \midvee_{M(G/H_+)} G/H_+) \] and
\[ \coef^{\ww} \left( \pi^{\ww}_0 (G/K_+) , \pi^{\ww}_0 ( \midvee_{H \in \ww}
\midvee_{M(G/H_+)} G/H_+)\right) \] is naturally isomorphic to \[ \dd (
G/K_+ , \midvee_{H \in \ww} \midvee_{M(G/H_+)} G/H_+ ) . \] Hence
there is a map \[ h (M) \col
\midvee_{K \in \ww} \midvee_{C_M (G/K_+)} G/K_+ \rarr \midvee_{H \in
\ww} \midvee_{M(G/H_+)} G/H_+ ,\] unique up to homotopy, so
$\pi_{0}^{\ww} ( h (M )) $ is isomorphic to the leftmost map in the
exact sequence \ref{shortex}.  Let $Z$ be the homotopy
cofiber of $h(M)$.  Proposition
\ref{AB} says that $ \pi^{\ww}_n ( G/K_+) =0$, for all $n <0$ and $K
\in \ww$.  So $\pi_n (Z) = 0$, for $n < 0$, and
there is a natural isomorphism $ \pi_{0}^{\ww} (Z) \cong M$.
Let $HM$ be $P_0 (Z)$, the 0-th Postnikov section of $Z$.
Then  $H M$ is an Eilenberg--Mac\,Lane spectrum and there is
a natural isomorphism \[ \pi^{\ww}_0 ( HM ) \cong M .\]
This proves the first claim.

The map $ Z  \rarr P_0 (Z) = HM$ induces an isomorphism \[ [ HM
, HN ] \rarr [ Z , HN ] . \] This gives an exact sequence
\[ 0 \rarr [ H M , HN ] \rarr [ \midvee_{H \in \ww}
\midvee_{M(G/H_+)} G/H_+ , HN ] \rarr [ \midvee_{K \in \ww}
\midvee_{C(G/K_+)} G/K_+ , HN ] ,
\] where the rightmost map is induced by $h(M)$
and the leftmost map is injective. Applying $\pi^{\ww}_0 $ gives an
isomorphism between the rightmost  map and the map \[ \coef \left(
\textstyle\bigoplus_{H \in \ww , M(G/H_+)} \pi^{\ww}_0 (G/H_+ ) , N
\right) \rarr \coef \left(\textstyle\bigoplus_{K \in \ww , C(G/K_+)}
\pi^{\ww}_0 (G/K_+ ) , N \right) .\] The kernel of this map is $ \coef ( M
, N )$, so $\pi^{\ww}_0 $ induces an isomorphism
\[ [ H M , HN ] \rarr \coef ( M ,N ) \] of abelian
groups. This proves the second claim.
\end{proof}
\begin{rem} When $d$ is not a constant class
function, then the homotopy groups of the objects in the heart need
not be concentrated in one degree. For example the heart of the
t-structure in Lemma \ref{Greenlees} consists of spectra of the form
$ F ( EG_+ , HM)$, where $M$ is an Eilenberg--Mac\,Lane spectrum.
The heart of the Postnikov t-structure on $\dd$ is not well
understood for general functions $d$.
\end{rem}

\subsection{Continuous $G$--modules} \label{sec:contcoeff}
When $\ww \not\subset \cc$ it is harder to describe the full
subcategory of the homotopy category of $\ww \cc \mm$ consisting
of the Eilenberg--Mac\,Lane spectra as a category of
coefficient systems. We give a description of the heart of the
Postnikov t-structure on the homotopy category of the $\overline{\lie
(G)}$--cofree model structure on $\mm_R$ when $G$ is a compact
Hausdorff group.
Let $R^0 $ denote the (continuous) $G$--ring $\colim_U \, \pi^{U}_0
(R)$. The ring $R^0$ has the discrete topology. By a continuous
$R^0 $--$ G$--module we mean an $R^0 $--$ G$--module with the discrete
topology such that the action by $G$ is continuous.
The module $\Pi^{1}_0 (X ) \cong \colim_U \, \pi^{U}_0 (X)
$ is a continuous $R^0 $--$ G$--module, for any $R$--module $X$.
The stabilizer of any element $m$ in a  continuous
$R^0 $--$ G$--module  $ M$ is in $\overline{\fin (G)}$.

\begin{pro} \label{pro:heartcontmodules}
The heart of $\dd$ is equivalent to the category of continuous
$R^0 $--$G$--modules and continuous $G$--homomorphisms between them.
\end{pro}
\begin{proof}
Let $M$ be a continuous $R^0 $--$ G$--module.  We get a canonical
surjective map \[f \col \textstyle\bigoplus_m R^0 [ G/ \stable(m) ] \rarr M
\] where the sum is over all elements $m$ in $M$. The map
$R^0 [ G/ \stable (m)] \rarr M $, corresponding to the summand $m$, is
given by sending the element $(r , g)$ to $r
\cdot gm$. This is a $G$--map since $ g' r \cdot g' gm = g'  ( r gm )$, for $g' \in G$.
Repeating this construction with the kernel of $f$ gives a canonical
right exact sequence of continuous $R^0$--$G$--modules
\begin{equation} \label{eq:contcoeff} \textstyle\bigoplus R^0 [ G/U']
\stackrel{h}{\rarr} \textstyle\bigoplus R^0 [ G/U] \rarr M  \rarr 0 . \end{equation}

We want to realize this sequence at the  level of spectra. We have that \[
\Pi^{1}_0 ( R\wedge G/U_+ ) \cong R^0 [ G/U] \] as
$R^0$--$G$--modules, for all $U \in \overline{\fin (G)}$. The map $h$ is
realized as $\Pi_{0}^{  1  }$ applied to a map \[ h (M) \col \midvee
R \wedge G/U'_+ \rarr \midvee R \wedge G/U_+ .\] Let $Z$ be the homotopy
cofiber of $h(M)$ and let $ HM$  be   $P_0 (Z)$.
The right exact sequence in \ref{eq:contcoeff} is
naturally isomorphic to $\Pi^{1}_0 $ applied to the sequence
\[ \midvee R \wedge G/U'_+ \rarr \midvee R \wedge G/U_+ \rarr HM . \]
Proposition \ref{AB} gives
that $\Pi^{ 1 }_n (HM) =0 $ when $n \not= 0$, and there is a
natural isomorphism $\Pi^{  1  }_0 (HM) \cong M$ of continuous
$R^0 $--$G$--modules, for any continuous $R^0$--$G$--module $M$.
It remains to show that $\Pi^{1}_0$ is a full and faithful functor. The same argument as in the
proof of Proposition \ref{postnikovhearth} applies. This works since
$\dd ( \Sigma^{\infty}_R G/U+ , HM ) $ is naturally isomorphic to $ M^U$.
\end{proof}

\section{Pro--$G$--spectra} \label{pro--structure}

In this section we use the $\ww$--$ \cc$--Postnikov t-model structure
on $\mm_R$, discussed in Section \ref{sec:t-modelstructures}, to
give a model structure on the pro--category, $\p \mm_R$. For
terminology and general properties of pro--categories see for
example \cite{tfi,strict}. We recall the following.
\begin{defn} \label{defn:prodef}
 Let $M$ be a collection of maps in
$\cc$. A levelwise map $g= \{ g_s \}_{s \in S}$ in $\p \cc$ is a
levelwise $M$--map  if each $g_s$ belongs to $M$. A pro--map $f$
is an  essentially levelwise $M$--map  if $f$ is isomorphic,
in the arrow category of $\p \cc$, to a levelwise $M$--map. A map in
$\p \cc$ is a  special $M$--map  if it is isomorphic to a
cofinite cofiltered levelwise map $f = \{ f_s \}_{s \in S}$ with the
property that for each $s \in S$, the map
\[
M_s f\col X_s \rarr \lim_{t <s } X_t \times_{ \lim_{t <s } Y_t} Y_s
\]
belongs to $M$.
\end{defn}

\begin{defn} \label{defn:extensionpro}
Let $F \col \ee \rarr \bb $ be a functor between two categories
$\ee$ and $\bb$. We abuse notation and let $F \col \p \ee\rarr \p
\bb$ also denote the extension of $F$ to the pro--categories given by
composing a cofiltered diagram in $\ee$ by $F$. We say that we apply
$F$ levelwise to pro--$\ee$.
\end{defn}

\subsection{Examples of pro--$G$--Spectra}

We list a few examples of pro--spectra.
\begin{enumerate}
\item
The finite $p$--local spectra $ M_I $ constructed by Devinatz
assemble to give an interesting pro--spectrum $\{ M_I \}$
\cite{dev}. The pro--spectrum is more well behaved than the
individual spectra. This pro--spectrum is important in understanding
the homotopy fixed points of the spectrum $E_n$ \cite{dav}
\cite{dho}.
\item
There is an approach to Floer homology that is based on pro--spectra
\cite{cjs}.
\item
The spectrum $\mathbb{RP}^{\infty}_{-\infty} $, and more generally,
the pro--Thom spectrum associated to a  (virtual)  vector bundle over a space $X$, are
non-constant pro--spectra.
\end{enumerate}

\begin{con} \label{con:EGtower} Let $N$ be a normal subgroup
of $G$ in $\cc$. Let \mdfn{$EG/N$} denote the free contractible
$G/N$--space constructed as the (one sided) bar construction of
$G/N$. Then $EG/N_+$ is a cell complex built out of cells $ (G/N
\times D^m )_+ $ for integers $m \geq 0$. The bar construction gives
a functor from the category with objects quotient groups $ G/N $, of
$G$, and morphisms the quotient maps, to the category of unbased
$G$--spaces. In particular, we get a pro--$G$--spectrum $\{
\Sigma^{\infty} EG/N_+ \}$ indexed on the directed set of normal
subgroups $N \in \cc$ ordered by inclusion. This pro--spectrum plays
an important role in our theory. The notation is slightly ambiguous;
the $N$--orbits of $EG$ are denoted by $(EG)/N$.
\end{con}

\subsection{The Postnikov model structure on pro--$\mm_R$}
\label{sec:postmodelpro}
 The most immediate candidate for a model
structure on $\p \mm_R$ is the strict model structure obtained from
the $\ww$--$\cc$--model structure on $\mm_R$ \cite{strict}. In
this model structure the cofibrations are the essentially levelwise
$\cc$--cofibrations and the weak equivalences are the \ess\ $\ww$--equivalences.

Let $\{ P_n \} $ denote a natural Postnikov
tower (for $ n \geq 0$) \cite[Section 7]{tfi}. The
Postnikov towers in $\mm_R$ extend to $\p \mm_R$ (Definition
\ref{defn:extensionpro}). A serious drawback of the strict model
structure on pro--$ \mm_R$ is that $ Y \rarr
\holim_n \, P_n Y $ is not always a weak equivalence.
We explain why. The homotopy limit in $\p \mm_R$,
$\holim_n \, P_n \{ Y_s \} $, of the Postnikov tower $ P_n \{ Y_s \}$,
is strict weakly equivalent to the pro--object $\{ P_n
Y_s \}$, for any $\{ Y_s \} \in \p \mm_R$.  Let  $Y$ be the
constant pro--object  $ \vee_{n\geq 0} \, H\zz [n] $, where
$H \zz [n]$ is the Eilenberg--Mac\,Lane spectrum of $\zz$
concentrated in degree $n$. Then the  map  $Y \rarr \{ P_n Y \}$
is not a strict pro--equivalence because $ \oplus_{n \geq 0} \, \pi_n $
applied to this map is not  an isomorphism of pro--groups.

In this section we construct an alternative to the
strict model structure to rectify the flaw that
$ Y \rarr \holim_n \, P_n Y$ is not always a weak equivalence.
This model structure has the same class of cofibrations but
more weak equivalences than the strict model structure. The
class of weak equivalences is the smallest class of maps closed
under composition and retract containing both
the strict weak equivalences and  maps
$ Y \rarr \{ P_n Y \}_{n \in \zz}$, for $Y \in \p \mm_R$. This
model structure on $\p \mm_R$ is called the
\mdfn{ Postnikov $\ww$--$\cc$--model structure}.
 We construct the Postnikov $\ww$--$\cc$--model structure
on $\p \mm_R$ from the Postnikov t-structure on $\ww \cc \mm_R$
using a general technique developed in \cite{ffi,tfi}.

The benefits of replacing pro--spaces (and pro--spectra) by their
Postnikov towers  was already made clear by   Artin--Mazur \cite[\S 4]{am}.
Dwyer--Friedlander also made use of this replacement \cite{dwf}.

In the next theorem we give a general model structure, called the
\mdfn{$ d$--Postnikov $\ww$--$\cc$--model structure} on $\p \mm_R$.
The Postnikov $\ww$--$\cc$--model structure referred to above is the
model structure obtained by letting $d$ be the constant class
function that takes the  value $0$.
\begin{thm} \label{promodelcategory}   Let $\cc$ be a $\uu$--\good\
collection of subgroups of $G$ and let $\ww$ be a collection of
subgroups of $G$ such that $ \cc \ww \subset \cc$. Let $R$ be a
$\cc$--connective ring spectrum. Let $d \col \ww \rarr \zz \cup \{
-\infty ,\infty \} $ be a class function.
Then there is a proper simplicial stable model structure on
$\p \mm_R$ such that:
\begin{enumerate} \item the cofibrations are essentially levelwise
retracts of relative $\cc$--cell complexes;
\item The weak equivalences are \ess\
$n$--$ \ww$--equivalences in the $d$--Postnikov  t-model
structure on $\ww  \cc \mm_R$ for all integers $n$; and
\item the fibrations are retracts of special $\ff$--maps.
\end{enumerate}  Here $\ff$ is the
class of all maps that are both $\ww$--$ \cc$--fibrations and
co--$n$--$\ww$--equivalences, for some $n$, in the $d$--Postnikov t-model
structure on $ \ww \cc \mm_R$.
\end{thm}

\begin{proof} This is a consequence of
 \cite[6.3, 6.13]{tfi},
\cite[16.2]{pro}, and Lemma  \ref{lem:t-structure}.
\end{proof}

We consider a particular example. Since $ \{ 1 \} \cc = \cc$ there
is a $d$--Postnikov $\{ 1 \} $--$ \cc$--model structure on $\p \mm_R$.
If $ d$ is $ + \infty$, then the model structure  is the strict model
structure on $\p \mm_R$ obtained from the $\cc$--cofree model structure
on $\mm_R$ and if $d$ is $- \infty$, then all maps are weak equivalences.
If $d$ is an integer, then the model structure is
independent of the integer $d$, so we
omit it from the notation. We call this model structure the \mdfn{
(Postnikov) $\cc$--cofree model structure} on $\p \mm_R$ and denote
it \mdfn{$\cc$--cofree $\p \mm_R$}.
\begin{thm} \label{thm:Cfreemodelstructure}
Let $\uu$ be a trivial $G$--universe. Then there is a model structure
on the category  $\p \mm_R$ such that:
\begin{enumerate} \item the cofibrations are essentially levelwise
retracts of relative $\cc$--cell complexes;
\item the weak equivalences are \ess\
$\cc$--underlying $n$--equivalences, for all integers $n$; and
\item the fibrations are retracts of special $F_{\infty}$--maps.
\end{enumerate}  Here $F_{\infty}$ is the class of all maps that
are both $\{ 1 \}$--$ \cc$--fibrations and $\cc$--underlying
co--$n$--equivalences, for some $n$, in the Postnikov t-model structure on
$\{ 1 \} \cc\mm_R$.
\end{thm}

We return to the general situation. Let \mdfn{$ \dd$} denote the
homotopy category of the $\ww$--$\cc$--model structure on $\mm_R$, and
let \mdfn{$\pp$} denote the homotopy category of the $d$--Postnikov
$\ww$--$\cc$--model structure on $\p \mm_R$.

An alternative description of the weak
equivalences in the Postnikov model structure is given in \cite[9.13]{tfi}.
There is a t-structure on $\pp$ as described in \cite[9.4]{tfi}.

Let \mdfn{$\map$} denote the simplicial mapping space in $\mm_R$
with the $\ww$--$ \cc$--model structure. We give a concrete description
of the homsets in the homotopy category of the $d$--Postnikov $ \ww$--$
\cc$--model structure on $\p \mm_R$.
\begin{pro} \label{homspaces}  Let $X$ and $Y$ be objects in
$\p \mm_R$ such that each $X_a $ is cofibrant and each $d$--Postnikov
section $P_n Y_b$ is fibrant in $\ww\cc\mm_R$. Then the
group of maps from $ X $ to $ Y $, in the homotopy category of the
$d$--Postnikov $\ww$--$\cc$--model structure on $\p \mm_R$, is equivalent to
\[ \pi_0 (\holim_{n , b} \, \hocolim_a \, \map  ( X_a , P_n Y_b )). \]
\end{pro}
\begin{proof} This follows from \cite[8.4]{tfi}.
\end{proof}

Recall that the constant pro--object functor $c \col \mm_R \rarr
\p\mm_R$ is a left adjoint to the inverse limit functor $ \lim \col
\p\mm_R \rarr \mm_R$. The composite functor $\lim \circ c $ is
canonically isomorphic to the identity functor on $ \mm_R$.
\begin{pro} Let $\mm_R$ have the $\ww$--$\cc$--model structure,
and let $ \p \mm_R$ have the $d$--Postnikov $\ww$--$\cc$--model
structure. Then $ c$ is Quillen left adjoint to $ lim$. If $d$ is a
uniformly bounded below class function ($d \geq n$ for some integer
$n$), then the constant pro--object functor
 $c$ induces a full embedding $ c \col  \text{Ho} ( \mm_R ) \rarr
 \text{Ho} ( \p \mm_R) $.
\end{pro}
\begin{proof}  It is clear that $c$ respects cofibrations and
acyclic cofibrations. Let $X$ and $Y$ be in $\mm_R$.
The assumption on $d$ gives that  $ Y \to \holim_{n } \,
P_n Y$ is a $\ww$--weak equivalence.
The result follows from Proposition \ref{homspaces}  since
\[ \holim_{n }\,  \map ( X , P_n Y )) \to
\map ( X , \holim_{n } \, P_n Y)) \] is a weak equivalence of simplicial sets.
\end{proof}

\begin{rem}
If $d$ is a uniformly bounded (above and below) class function on
$\ww$, then a map is an essentially levelwise $(n +d)$--equivalence
for every integer $n$, if and only if it is an essentially levelwise
$n$--equivalence for every integer $n$ (with the constant function with value $0$).
Hence under this assumption the $d$--Postnikov $\ww$--$\cc$--model
structure on $\p \mm_R$ is the same as  the
Postnikov $\ww$--$\cc$--model structure on $\p \mm_R$.
\end{rem}

\begin{rem} \label{rem:adams}
It is not clear if there is an Adams
isomorphism when $G$ is not a compact Lie group or  $\{ 1 \}$ is
not contained in $\ww$. There are no free $G$--cell complexes (that
are cofibrant) so the usual statements does not make sense. One
might try to replace $G$ by $\{ G/N \}$, indexed on normal
subgroups, $N$, of $G$ in $\ww$. The most naive implementation of
this approach does not work.

Assume that $G$ is a compact Hausdorff
group which is not a Lie group, and let $\cc =\ww = \overline{ \lie (G)
}$. Assume in addition that $\uu$ is a
complete $G$--universe. Proposition \ref{pro:tom diecksplitting} and
Proposition \ref{homspaces} applied to the pro--suspension spectrum $ \{
 \Sigma^{\infty} EG/N_+\}$ give that
 $ \pi^{G}_* ( \{ \Sigma^{\infty} EG/N_+\} )$ is $0$.
\end{rem}

\begin{exmp} \label{exmp:equivalenceofcoefficients}
It is harder to be
an essentially levelwise $\pi^{\ww}_n$--isomorphism  than it is to
be an essentially levelwise $\pi^{H}_n$--isomorphism for each $H \in
\ww$ individually. The difference is fundamental as the following
example shows (in the category of spaces, or the category of spectra
indexed on a trivial universe). Let $\ww$ be a normal collection
that is closed under intersection. Let $N$ be a normal subgroup of
$G$. The fixed point space $(EG/N)^H $ is empty, for $ H \nleq N$,
and it is $EG/N$, for $H \leq N$. The pro--map $\{ \Sigma^{\infty}
EG/N_+\} \rarr \{ * \}$, where the first object is indexed on the
directed set of normal subgroups of $G$ contained in $\ww$ ordered
by inclusion, is a $\pi_{n}^H$--isomorphism, for all $H \in \ww$ and
any integer $n$. But this map is typically not an essentially
levelwise $\pi^{\ww}_n$--isomorphism for any $n$. The same conclusion
applies when the universe is not trivial by (an (in)complete
universe version of) Proposition \ref{pro:tom diecksplitting} (see
Remark \ref{rem:adams}).
\end{exmp}

We include a result about fibrations for use in Section \ref{sec:htpfixedpt}.
\begin{lem} \label{lem:pro-restrictfib}
Let $\stable ( \uu) $ be a $\ww$--\good\ collection of subgroups of a compact
Hausdorff group $G$. Let $f \col X \rarr Y$ be a
fibration in $\p  \mm_R$ between  fibrant objects $X$ and $Y$
in the Postnikov $\ww$--model structure on $\p G \mm_R$. Assume in addition that $X$ and
$Y$ are levelwise $\ww$--$\sss$--cell complexes. Let $K$ be any closed subgroup of
$G$, and let $\ww'$ be a  collection of subgroups of $K$
such that $\ww' \ww \subset \ww$ and $\stable_K ( \uu )$ is $\ww' $--\good.
 Then $ f$ regarded as a map of pro--$K$--spectra is a
fibration in the Postnikov $\ww' $--model structure on $ \p K \mm_R$.
\end{lem}
\begin{proof}
This reduces to Lemmas \ref{lem:restrictfib} and
\ref{lem:restrict-n-equivalence}, since fixed points respect the
limits used in the definition of special $F_{\infty}$--maps (see
Definition \ref{defn:prodef}).
\end{proof}

\subsection{Tensor  structures on $\p \mm_R $} \label{sec:tensor}
Let $R$ be a ring. Then the category $\mm_R$ is a closed symmetric
tensor category. The category pro--$\mm_R$ inherits a symmetric
tensor structure. Let $\{X_s \}_{s \in S} $ and $\{Y_t \}_{t \in T}$
be two objects in $\p \mm_R$. \begin{defn} The smash product $\{
X_{s} \}_{s \in S} \wedge \{ Y_{t} \}_{t \in T} $ is defined to be the
pro--spectrum $\{ X_{s} \wedge Y_{t} \}_{s \times t\in S \times T}
$. \end{defn}
The tensor product in $\p \mm_R$ is not closed.
Worse, the smash product does not commute
with direct sums in general \cite[11.2]{tfi}.
The tensor product does not behave well homotopy theoretically for
general $R$--modules.
\begin{defn}
\label{defn:boundedbelow}
A pro--object $ Y$ is  bounded below
if it is isomorphic to a pro--object $ X = \{ X_s \}$ and there
exists an integer $n$ such that each $
*\rarr X_s$ is an $n$--equivalence. A pro--object
$Y$ is  levelwise bounded below  if it is isomorphic to a
pro--object $ X = \{ X_s \}$ and for every $s$ there exists an
integer $n_s$ such that $* \rarr X_s$ is an $n_s$--equivalence.
\end{defn}

The simplicial structure on $\mm_R$ is compatible with the tensor
structure \cite[12.2]{tfi}. If $\cc$ is an \good\ collection of
subgroups of $G$, then $\p \mm_R$, with the strict model
structure, is a tensor model category \cite[12.7]{tfi}.

If $R$ is $\cc$--connective, then Lemma \ref{respecttensor}  gives that
the Postnikov t-structure on $\dd$ is
compatible with the tensor product \cite[12.5]{tfi}.
This implies that the  full subcategory  of $\p  \mm_R$, with the
Postnikov model structure, consisting of essentially bounded below
objects is a tensor model category \cite[12.10]{tfi}.

We can define a pro--spectrum valued hom functor. Let \mdfn{$F$}
denote the internal hom functor in $\mm_R$.
\begin{defn}
\label{defn:partialinnerhom} We extend the definition of $F$ to $\p
\mm_R$ by letting $F(X,Y) $ be the pro--object
\[ \{ \, \colim_{s \in S} \, F (X_s , Y_t )\, \}_{t \in T} .\]
\end{defn}

The pro--spectrum valued hom functor is not an internal hom functor in
general.  The next result shows that under some finiteness assumption
the (derived) pro--spectrum valued hom functor behaves as an internal hom functor
in the homotopy category.
\begin{lem} \label{innerhom}
Let $\{X_s \} $ be a pro--spectrum such that each $X_s$ is a retract
of a finite $\cc $--cell spectrum. Let $Y$ be an essentially
bounded below pro--spectrum and let $Z$ be a pro--spectrum.
Then there is an isomorphism
\[ \pp ( X \wedge Y , Z ) \cong \pp  ( X , F (Y , Z )) \]
\end{lem}
\begin{proof} This follows from Proposition \ref{homspaces}.
\end{proof}

\subsection{Bredon cohomology and group cohomology}
Assume that $ \ww \subset \cc$ or that
$ \ww = \{1 \}$ and $\cc = \overline{\lie (G)}$.
Under these assumptions the heart of the Postnikov t-structures are
equivalent to categories of coefficient systems as described in
Propositions \ref{postnikovhearth} and \ref{pro:heartcontmodules}.
In the latter case the
coefficient systems are continuous (discrete) $G$--modules.

\begin{defn} \label{dfn:Bredoncohomlogy}
The $n$-th \mdfn{Bredon cohomology} of a pro--spectrum $X$ with
coefficients in a pro--coefficient system $\{ M_a \}$, is defined to
be \[ \pp ( X , \Sigma^{\infty}_R S^n \wedge \{ H M_a \} ). \]
\end{defn} The Bredon  cohomology of $X$ with coefficients in $\{ M_a
\}$ is denoted by $ H^n ( X ; \{ M_a \} )$. This is the cohomology
functor with coefficients in the heart, in the terminology of
\cite[2.13]{tfi} (when $ \{ H M_a \}$ is in the heart of $\pp$
\cite[9.11, 9.12]{tfi}).

\begin{lem}
Let $M$ be a constant pro--coefficient system. Then the $n$-th Bredon
cohomology of a pro--spectrum $X = \{ X_s \}$ is
\[ \colim_s \,  H^n ( X_s ; M ) .\]
\end{lem}
\begin{proof} This follows from Proposition \ref{homspaces}.
\end{proof}

\begin{defn} \label{dfn:Bredonhomlogy}
The $n$-th \mdfn{Bredon homology} of an essentially bounded below
pro--spectrum $X$ with coefficients in a pro--coefficient system $\{
M_a \}$, is defined to be \[ \pp ( \Sigma^{\infty}_R S^n , X \wedge \{ H M_a \} ).
\]  \end{defn}
The Bredon homology of $X$ with coefficients in $\{ M_a \}$ is
denoted $ H_n ( X ; \{ M_a \} )$.

\begin{rem}
We get isomorphic groups if we  use the strict model structure instead of
the Postnikov model structure on $\p \mm_R$ to define the Bredon
cohomology and Bredon homology. In the case of cohomology
since $\{ H M_{\alpha} \}$ is bounded above, and in the case of
homology since we are mapping from a constant pro--object.
\end{rem}

Let $R$ be the sphere spectrum $\sss$,
$ \ww = \{1 \}$ and $\cc = \overline{\lie (G)}$. Let $ \{ M_{a} \} $
be a pro--object of discrete $G$--$R^0$--modules, and let $ \{ HM_{a} \}$
be the associated Eilenberg--Mac\,Lane pro--spectrum.

\begin{defn} \label{defn:gpcohomology} The
\mdfn{group cohomology of $G$ with coefficients in $\{ M_{a} \}$} is
the   Bredon cohomology of $ \{ EG/N_+ \}$ with
coefficients in $ \{ HM_{a} \}$. \end{defn}
We denote the $n$-th group cohomology
by \mdfn{$H_{\cont}^n ( G ; \{ M_{a} \} ) $}.
If $ M$ is a constant coefficient system, then we
recover the usual definition of group cohomology as
\[ H_{\cont}^n ( G ;  M  ) \cong \colim_{N}\, H^n ( G/N ;
M^N ), \] where the colimit is over all subgroups $N \in \lie (G)$.
In general, there is a higher lim spectral sequence relating the
group cohomology of a pro--coefficient system $\{ M_{a}
\}$ to the continuous group cohomology of the individual modules $
M_{a }$.

\begin{lem} A short exact sequence of pro--$G$--modules gives a long exact
sequence in group cohomology.
\end{lem}
\begin{proof} This follows from the fact that a short exact sequence
of pro--$G$--modules gives rise to a cofiber sequence, of the corresponding
Eilenberg--Mac\,Lane pro--spectra, in the Postnikov
$\overline{\lie (G)} $--model structure on pro--$\mm_R$.
\end{proof}

\begin{lem} The group cohomology functor in Definition
\ref{defn:gpcohomology}, composed with the functor from towers of
discrete $G$--modules and levelwise maps between them to pro--$G$--coefficient
systems, agrees with Jannsen's group cohomology \cite{jan}.
\end{lem}
\begin{proof}
A comparison to Jannsen's
cohomology follows from the proof of \cite[Lemma 3.30]{jan}.
\end{proof}

\subsection{The Atiyah--Hirzebruch spectral sequence}
\label{sec:AH} A t-structure on a triangulated category gives rise
to an Atiyah--Hirzebruch spectral sequence \cite[10.1]{tfi}. Let
$\cc $ be an \good\ collection of subgroups of $G$ and assume
that $\stable (\uu ) $ is $\ww$--\good.  Assume that $ \ww \subset \cc$ or
$ \ww = \{1 \}$ and $\cc = \overline{\lie (G)}$. We
can relax this assumption to $\cc \ww \subset \cc$ if we work with
objects in the heart instead of coefficient systems. Let $R$ be a
$\cc$--connective ring spectrum. Let $\pp$ denote the homotopy
category of $\p \mm_R$ with the Postnikov $\ww$--$\cc$--model
structure.  Let square brackets denote
homotopy classes in $\pp$. Recall Definition  \ref{defn:boundedbelow}.
\begin{thm} \label{thm:AH-spectralsequence}
Let $X$ and $Y$ be any pro--$G$--spectra. Then there is a spectral
sequence with
\[ E_{2}^{p,q} = H^p ( X , \Pi_{-q}^{\ww}(Y) ) \]
converging to $[X, Y]^{p+q}_G$. The differentials have degree $(r,
-r + 1)$. The spectral sequence is conditionally convergent if:
\begin{enumerate}
\item
$ X$ is bounded below; or
\item  $X$ is levelwise bounded below and $Y$ is a constant
pro--$G$--spectrum.
\end{enumerate}
\end{thm}

The $E_2$--term is the Bredon cohomology of $X$ with coefficients in
the pro--coefficient system $ \Pi_{-q}^{\ww} (Y)$ (Definition
\ref{dfn:Bredoncohomlogy}).
\begin{proof}  This follows from
\cite[10.3]{tfi} and our identification of the heart in Propositions
\ref{postnikovhearth} and \ref{pro:heartcontmodules}. \end{proof}

\begin{lem} If $Y$ is a monoid in the homotopy category of
 $\p \mm_R$ with the
strict model structure obtained from the $ \ww$--$ \cc$--model
structure on $\mm_R$, then the spectral sequence is multiplicative.
\end{lem}
\begin{proof}
By Lemma \ref{respecttensor} the Postnikov t-structure respects the
smash product in the sense of \cite[12.5]{tfi}. The result follows from
\cite[12.11]{tfi}.
\end{proof}

\begin{rem} \label{ref:alternativeAH}
 If $X$ is a bounded below  $CW$--$R$--module, then it is
possible to filter $ F ( X , Y)$ by the skeletal filtration of $X$
instead of the Postnikov filtration of $Y$. The two filtrations give
rise to two spectral sequences. When $Y$ is a constant bounded above
pro--spectrum the two spectral sequences are isomorphic
\cite[App.B]{gm}. However, in  general we
get two different spectral sequences. For a discussion of this see
\cite{dav3}.
\end{rem}

\section{The $\cc$--free model structure on $\p \mm_R$}
Suppose  $\cc$ is an \good\ collection that does not contain the
trivial subgroup  $1 $. Then it is not possible to have a
$\ww$--$\cc$--model structure on $ \mm_R$ such that the cofibrant
objects are free $G$--cell complexes; this is so because $
* \rarr G_+$ is a $\ww$--equivalence for any collection $\ww$
such that $\cc \ww \subset \cc$. When $ \cap_{ H \in \cc} H = 1 $,
it turns out that it is possible to construct a model structure on
$\p \mm_R$ with cofibrations that are arbitrarily close approximations
to $G$--free cell complexes. The class of weak equivalences in this model
structure is contained in the class of weak equivalences in the Postnikov
$ \{ 1\}$--$ \cc$--model structure on $\p \mm_R$.

In this section we assume that $\uu$ is a trivial $G$--universe, and
that $\cc$ is a normal \good\ collection  of subgroups of $G$ (see Definitions
\ref{defn:collection} and \ref{good}). We also assume that $R$ is a
$\cc$--connective ring spectrum.

\subsection{Construction of the $\cc$--free
model structure on $\p \mm_R$} \label{sec:free}
We use the framework of filtered model structures
defined in \cite[4.1]{ffi}.
Define an indexing set  $A$ to be the product of $\cc$, ordered by containment, and the
integers, with the usual totally ordering.
\begin{lem}
There is a proper simplicial filtered model structure on $\mm_R$,
indexed on the directed set $A$, such that:
\begin{enumerate}
\item $C_{U,n} = C_U$ is the class of  retracts of
relative $G$--cell complexes with cells of the form
$\Sigma^{R}_{V} G/H_+ \wedge D^{m}_+ $, for some integer $m$,
indexing representation $V$  and $H \in \cc$ such that $H \leq U$;
\item $F_{U,n}$ is the class of  maps $f$ such that $f^H$ is a
fibration and a co--$n$--equivalence, for each $H \in
\cc$ such that $H \leq U$; and
\item $W_{U,n} = W_n $ is the class of maps $f$ for which  there
exists an $H \in \cc$ such that $f^K$ is an $n$--equivalence for
every $K \in \cc$ such that $K \leq H$.
\end{enumerate} \end{lem}
\begin{proof}
The directed set of classes $C_U$ and $W_{n}$ are decreasing and the
directed set of classes $F_{U,n}$ is increasing. The verification of
the proper filtered model structure axioms is similar to the
verification of the axioms in the case of $G$--spaces. We omit the
details and refer the reader to the detailed discussion given in
\cite[Section 8]{ffi}.
The simplicial structure follows as in the proof of Theorem \ref{model2}.
\end{proof}
Let $F_{\infty}$ denote the union $ \cup_{U,n} F_{U,n}$. The following
model structure on $\p \mm_R$ is a consequence of \cite[Theorems
5.15, 5.16]{ffi}.
\begin{thm} \label{thm:freemodelstructure}
There is a proper simplicial model structure on $\p \mm_R$ such that:
\begin{enumerate} \item the cofibrations are maps that are
retracts of essentially levelwise $C_U$--maps for every $U \in \cc$;
\item the weak equivalences are maps that are \ess\
$W_{n}$--maps for every $n \in \zz$; and
\item the fibrations are special $F_{\infty}$--maps.
\end{enumerate}
\end{thm}
We call this model structure the \mdfn{$\cc$--free model structure}
on $\p \mm_R$. We denote the model category by
 \mdfn{$\cc $--free $\p \mm_R$}.

\begin{lem} \label{lem:freecofibrant}
Let $\{\Sigma_{R}^{\infty} EG/N_+ \}$ be indexed on all $N \in \cc$ which
are normal in $G$. If $X = \{ X_s \}$ is cofibrant in the $\cc$--model structure
on $\p \mm_R$, then $ X \wedge \{\Sigma_{R}^{\infty} EG/N_+ \}$, is
a cofibrant replacement of $X$ in the $\cc$--free model structure on
$\p \mm_R$.
\end{lem}
\begin{proof}
By our assumption each $X_s$ is a retract of a $\cc$--cell complex.
Hence, since $\cc$ is an \good\ collection, $ X_s \wedge
\Sigma_{R}^{\infty} EG/N_+$ is also a retract of a $\cc$--cell
complex built out of $G/H$--cells for $H \leq N$. Since $\cc$ is a
normal \good\ collection, $ \ast \rarr X \wedge
\{\Sigma_{R}^{\infty} EG/N_+ \}_{N \leq U}$ is an \ess\ $C_U$--map
for every $U \in \ww$. We conclude that $ X \wedge
\{\Sigma_{R}^{\infty} EG/N_+ \} $ is a cofibrant replacement of $X$.
 \end{proof}
 In particular, $\{ \Sigma_{R}^{\infty}
 EG/N_+ \} \rarr R $ is a cofibrant replacement of $R $ in the
 $\cc$--free model structure on $\mm_R$.

\subsection{Comparison of the free and the cofree model
structures }

We compare
the $\cc$--cofree model structure on $ \p \mm_R$, given in Theorem
\ref{thm:Cfreemodelstructure}, with the $\cc$--free model structure
on $ \p \mm_R$, given in Theorem \ref{thm:freemodelstructure}.

Clearly a $\cc$--free cofibration is a $\cc$--cofree cofibration,
and a $W_{U,n}$--equivalence is a $\cc$--underlying $n$--equivalence.
Hence the identity functors give a Quillen adjunction
\[  \cc \text{--}
\free \p \mm_R \leftrightarrows \cc \text{--} \cofree \p \mm_R . \]

If $1$ is in $\cc$, then the $\cc$--free model structure on $\p \mm_R$ is
the model structure obtained from the Postnikov  t-model structure on
$\{ 1 \} \{ 1 \} \mm_R$. Hence the $\cc$--free model structure on $\p
\mm_R$ is Quillen equivalent to the $\cc$--cofree
model structure  on $\p \mm_R$ by Proposition \ref{pro:c1c2}.
If  $ 1 $ is not an element in $\cc$, then the $\cc$--free and the
$\cc$--cofree model structures on $\mm_R$ are typically not Quillen
equivalent, as shown by the next example.

\begin{exmp} Let $ f \col \vee_N \Sigma_{\sss}^{\infty} EG/N_+ \rarr \vee_N\,  \Sigma^{\infty}_{\sss} S^0  $ be the sum
of the collapse maps for all normal subgroups $N \in \cc$.
The map $f$ is a $\cc$--underlying
equivalence, but if $\cc$ does not contain a smallest element
(ordered by subconjugation), then $f$ is not a $\cc$--free weak
equivalence by Remark \ref{rem:adams}.
\end{exmp}

Let $X$ be a cofibrant object and let $Y$ be a fibrant object in
$\p \mm_R $ with the  Postnikov $\cc$--model structure.
Then, by Lemma \ref{lem:freecofibrant}, $X \wedge \{
\Sigma_{R}^{\infty} EG/N_+ \} $ is a cofibrant object in the
$\cc$--free model structure on $\p \mm_R$, and $Y$ is a fibrant object in
the $\cc$--free model structure on $\p \mm_R$. The cofibrations in
the $\cc$--cofree model structure and the Postnikov $\cc$--model
structure on $\p \mm_R$ are the same. For the purpose of calculating
mapping spaces in the $\cc$--cofree model structure on $\p \mm_R$ it
suffices to replace $Y$ by a levelwise fibrant replacement in the $
\cc $--cofree model structure on $\mm_R$. This follows from
\cite[5.3,7.3]{tfi} since $Y$ is \ess\ $\ww$--bounded above. We
choose a natural fibrant replacement for bounded above objects in
the $\cc $--cofree model structure on $\mm_R$ and denote it by
adding a subscript $f$. We describe the fibrant replacement functor
in the $\cc$--cofree model structure on $\mm_R$ after the next lemma.
\begin{lem}
\label{lem:groupcohomology}
Let $G$ be a compact Lie group and let $ \cc = \overline{\lie (G)}$.
Let $L < K $ be two   normal subgroups of  $G$.
Let $\tilde{M}$ be a $\cc$--coefficient system. The group of
equivariant (weak) homotopy classes of maps
\[ [ \Sigma^{\infty}_R  EG/K_+ , H\tilde{M} ]_G
\] is isomorphic to the group
cohomology of $G/K$ with coefficients in the $G/K$--module $M (G/K)$.
 \end{lem}
\begin{proof} This follows from the equivariant
 chain homotopy description  of Bredon cohomology.  See for example
\cite[8.1]{art}.
  \end{proof}

Recall that $F$ denotes the internal hom functor in $\mm_R$.
\begin{pro}
\label{pro:cofreefibrant} Assume that $G$ is a compact Hausdorff
group, $\cc = \overline{\lie (G)}$, and  the
universe $ \uu$ is trivial. Let $R$ be a $\cc$--connective ring, and let $Y$
be a $\cc$--bounded above fibrant object in the $\cc$--model
structure on $\mm_R$. Then the two maps
\[  \hocolim_{N}\, F (  \Sigma_{R}^{\infty}  EG/N_+ , Y  )
\longrightarrow \hocolim_{N}\, F ( \Sigma_{R}^{\infty} EG/N_+ , Y_f )
\longleftarrow Y_f \] are weak equivalences in $\cc \mm_R$.
\end{pro}
\begin{proof}  We need to show that both maps induce isomorphisms on
$ \pi_{*}^{\cc}$. For $ K\in \cc$ we get that $\pi_{n}^K$ applied to
the sequence above is isomorphic to
\begin{equation} \label{eq:freecofree} \colim_N \, [
\Sigma_{R}^{\infty} EG / N_+
\wedge \Sigma_{R}^{\infty} G/K_+ , Y ] \longrightarrow \end{equation} \[ \colim_N \, [
\Sigma_{R}^{\infty} EG / N_+ \wedge \Sigma_{R}^{\infty} G/K_+ , Y_f  ] \longleftarrow [
S^0 \wedge \Sigma_{R}^{\infty} G/K_+ , Y_f  ] \]
 where the square brackets denote the homomorphism groups
in the homotopy category of $\cc G \mm_R$. The second map is an
equivalence since \[ \Sigma_{R}^{\infty} G/K_+ \wedge  \Sigma_{R}^{\infty} E G/ N_+
\rarr \Sigma_{R}^{\infty} G/K_+\] is an underlying equivalence of cofibrant objects in
the $ \cc $--cofree model structure on $G \mm_R$, for any normal subgroup $N$ of $G$ in $\cc$.

We now prove that the first map is an isomorphism.
We first consider the case when  $K=G$.
There is a map between conditionally convergent spectral sequences converging to
the first map in \ref{eq:freecofree}. The spectral sequences are
the Atiyah--Hirzebruch spectral sequences in
the Postnikov $\cc $--model structure on $\p G \mm_R$
\cite[10.3]{tfi}; see also  Subsections \ref{sec:contcoeff} and \ref{sec:AH}.
The map between the $E_2$--terms is
\[ \colim_N \,  H^{p}_{\cont} ( G/ N , \Pi^{ N  }_{-q} (  Y) )
\rarr \colim_N \,
H^{p}_{\cont} ( G/N , \Pi^{ N }_{-q} (  Y_f ))  \] induced by $Y \rarr Y_f$.
This map is isomorphic to
\[ \colim_N \,  H^{p}_{\cont} ( G , \Pi^{ N  }_{-q} (  Y) )
\rarr \colim_N \,
H^{p}_{\cont} ( G , \Pi^{ N }_{-q} (  Y_f )).  \]
This is an isomorphism since group cohomology commutes
with directed colimits of continuous $G$--modules
and   $\colim_N \, \Pi^{ N  }_{-q} (  Y)
\rarr \colim_N \, \Pi^{ N  }_{-q} (  Y_f ) $ is
an isomorphism of continuous $G$--modules.
The spectral sequences converges conditionally
since $\{ \Sigma_{R}^{\infty} EG / N_+ \}$ is bounded below.

For general $K \in \cc$, we use  Lemmas \ref{egk} and  \ref{lem:easyrestrictfib} together with  the adjunction  between $G \wedge_K -$ and the restriction functor $ G \mm_R \rarr K \mm_R$ to reduce the problem to the case $G = K$.
\end{proof}

\begin{cor} \label{cofreefibrant}
Let $Y$ be in $\mm_R$. Then  the fibrant replacement,
$ Y_f$, of $Y$ in the $\overline{\lie (G)}$--cofree
model structure on $\mm_R$ is equivalent to \[
\holim_m \, \hocolim_{ N} \, F ( \Sigma^{\infty}_R EG/N_+ , P_m
Y ) \] in the  $\overline{\lie (G)}$--model structure
on $\mm_R$, where the homotopy
colimit is over $N \in \lie (G)$, and the homotopy limit and colimit
are formed in the $\overline{\lie (G)}$--model structure.
\end{cor}

\begin{rem} Note that the proof of Proposition
\ref{pro:cofreefibrant} only uses the Postnikov $\ww$--$\cc$--model
structures on $\p \mm_R$. It does not depend on the existence of the
$\cc$--free model structure, which is technically more sophisticated.
\end{rem}

The next result clarifies the relationship between the $\cc$--free
and the $\cc$--cofree model structures on $\p \mm_R$.
Let $\map $ denote the simplicial mapping space in the $\cc$--model
structure on $\mm_R$. Recall that  $\{ P_n \} $ denotes a natural Postnikov
tower.
\begin{thm} \label{thm:free-cofree}
We assume that $G$ is a compact Hausdorff group, $\cc =
\overline{\lie (G)}$, and the universe is trivial. Let $R$ be a
$\cc$--connective ring.
 Let $\{X_s \}$ and $\{ Y_t \}$ be objects
in $\p \mm_R$ such that each $X_s$ is cofibrant, and each $P_n Y_t$
is fibrant in $\cc \mm_R$. Then the homset in the $\cc$--free model
structure on $\p \mm_R$ is isomorphic to
\[  \pi_0 \, \holim_{t,n} \, \hocolim_{s ,N} \,  \map
( X_s \wedge EG/N_+ , P_n Y_t ) ) , \] and the homset in the $\cc$--cofree model structure on $\mm_R$ is isomorphic to
\[  \pi_0 \, \holim_{t,n} \, \hocolim_{s} \, \map
( X_s , \hocolim_N \, F ( EG/N_+ , P_n Y_t ) ). \]
\end{thm}
\begin{proof}
This follows from the description of mapping spaces in
\cite[5.3,7.3]{tfi} and from Lemma \ref{lem:freecofibrant} and Proposition \ref{pro:cofreefibrant}.
\end{proof}

\begin{rem}
One can also let $\p \mm_{\sss}$ inherit a model structure from the
$\ww$--$\cc$--model structure on $G\mm_{\sss}$ along the (right adjoint) inverse
limit functor \cite[12.3.2]{hir}.
The weak equivalences are
pro--maps $ f \col X \rarr Y$ such that $f \col\lim_s \, X_s \rarr
\lim_t \, Y_t$ are weak equivalences in $\mm_R$. The fibrations are the
pro--maps such that the inverse limits are fibrations in $\mm_R$.
(This follows from the right lifting property.) We have that $c$ and
$\lim$ are a Quillen adjoint pair between $ \mm_R$ and $ \p \mm_R$.
This model structure does not play any role in this paper.
\end{rem}

\section{Homotopy fixed points} \label{sec:htpfixedpt}
We define homotopy fixed points of pro--$G$--spectra for closed subgroups
of $G$, and show that they behave well with respect to iteration.

\subsection{The homotopy fixed points of a pro--spectrum}
Let $G$ be a compact Hausdorff group, $\uu$ a trivial $G$--universe,
and $\cc $ the cofamily closure, $\overline{\lie (G)}$, of $\lie
(G)$. Let $R$ be a   $\cc$--connective $\sss $--cell complex
ring  with a trivial $G$--action. The last  assumption guarantees that we
can apply  Lemma \ref{lem:restrictfib}.
\begin{defn} \label{defn:htpfixedpoints}
Let $Y$ be a pro--$G$--$ R$--module. The \mdfn{ homotopy fixed point
pro--spectrum} \mdfn{$Y^{h G}$}  is defined to be the
$G$--fixed points of a fibrant
replacement of $Y$ in the Postnikov $\overline{\lie (G)}$--cofree
model structure on $\p \mm_R$.
\end{defn} By choosing a   fibrant  replacement functor, $Y \mapsto Y_f$,
we get a homotopy fixed point functor.

\begin{lem} \label{lem:description of htpfixed points}
Let $Y = \{ Y_b \}$ be a   pro--$G$--$R$--module which is levelwise fibrant and bounded above
in the  $\overline{\lie (G)}$--model structure on $ \mm_R$. Then the homotopy fixed point
pro--spectrum $ Y^{hG}$ is weakly equivalent to
\[ \{ (\hocolim_{N \in \lie (G)} \,
F (\Sigma_{R}^{\infty} EG/N_+ ,  Y_b   ))^G \} \] in the (non-equivariant) Postnikov
model structure on $\p \mm_R$.
\end{lem} \begin{proof} This follows from Proposition
\ref{pro:cofreefibrant}. \end{proof}

In particular, if $G$ is a compact Lie group and $Y$ is a fibrant $G$--spectrum, then the $G$--homotopy
fixed point pro--spectrum   $Y^{hG}$ is  equivalent to \[ F (\Sigma_{R}^{\infty} EG_+
, \{ P_n Y \} )^G \] in the Postnikov  model
structure on $\p \mm_R$. The associated spectrum is equivalent to $F (\Sigma_{R}^{\infty} EG_+
,  Y  )^G $ in   $\mm_R$. See Lemma \ref{lem:ordinaryhtpfps}.

If $Y$ is a pro--$G$--spectrum and $K $ is a subgroup of $G$, then we
expect the $K$--homotopy fixed point pro--spectrum to be equipped with
an action by $N_G K/K$. We make the following definition.

\begin{defn}
\label{defn:relativehtpfixedpoints} Let $Y$ be a pro--$G$--$ R$--module,
 and let $K$ be a closed subgroup of $G$.
The \mdfn{ $K$--$G$--homotopy fixed point pro--spectrum $Y^{h_G K}$} of
$Y$ is defined to be
\[ \hocolim_N \, ((Y_f)^{KN}) \] where the colimit is over
$N \in \lie (G)$, and $Y_f$ is a fibrant replacement of $Y$ in the
Postnikov $\overline{\lie (G)}$--cofree model structure on $\p
\mm_R$.
\end{defn}
The pro-spectrum $Y^{h_G K}$ is an
$N_G K/K$--pro--spectrum. If $K \in \overline{\lie (G)}$, then the canonical
map $Y^{h_G K} \rarr (Y_f )^K$ is an equivalence in the Postnikov
$\overline{\lie (N_G K / K)}$--model structure on $\p \mm_R$.
The restriction of $Y_f$ to a subgroup $K \in \overline{\lie (G)}$
 is  a fibrant object in the $\overline{\lie (K)}$--cofree model
 structure on $\p K \mm_R$ by Lemma \ref{lem:easyrestrictfib}.
Hence  $Y^{h_G K} $ is equivalent to $ (Y|K )^{h
K}$ in the Postnikov model structure on pro--$R$--modules, for all
 subgroups $K \in \overline{\lie (G)}$. This need not
be true when $K \not\in \overline{\lie (G)}$. For example, consider
the suspension $G$--spectrum $ \Sigma_{R}^{\infty} G_+$ and $K = \{ 1
\} \not\in \overline{\lie (G)}$. The next lemma  shows that for
certain pro--spectra $Y$, we still have that $Y^{h_G K} $ is equivalent to
$ (Y|K )^{h K}$ even when $K \not\in \overline{\lie (G)}$.

\begin{lem}  \label{lem:comparisonhtpfix}
Let $K \lhd L$ be two  closed subgroups of $G$.
Let $Y$ be a pro--$G$--spectrum that is both fibrant and cofibrant in
the Postnikov $\overline{\lie (G)}$--model structure on $\p
\mm_R$. Let $Y'$ be $\hocolim_{N \in \lie (G)} \,
F (\Sigma_{R}^{\infty} EG/N_+ ,   Y  ) $ (levelwise hocolimit). Then  $
Y^{h_G K}$  is equivalent to $(Y' | L)^K$  in the Postnikov
$\overline{\lie (L /K)}$--model structure on $\p L /K \mm_R$.
\end{lem}

Proposition
\ref{pro:cofreefibrant} says that  $Y'$ is equivalent to $Y_f$  in the
Postnikov $\overline{\lie (G)}$--model structure on $\p
\mm_R$. Note that we don't claim that $Y'$ is fibrant in the
Postnikov $\overline{\lie (G)}$--cofree model structure on $\p
\mm_R$.
\begin{proof}
In the following we use  that  homotopy colimits commute with fixed points.
The pro--spectrum $\hocolim_{N'}(Y')^{JN'} $ is equivalent
to $Y^{h_G J}$ in the Postnikov model structure,
for every  subgroup $J$ of $G$ containing $K$, where the
colim is over  $N' \in \overline{\lie (G)}$. So it suffices
to show that  $ (Y')^{h_G J}$  is equivalent to $(Y' | L)^J$ in the Postnikov
model structure  for $ K \leq J \leq L$.  Consider the map
\[ \hocolim_{N'\in \lie (G)} \,(\hocolim_{N \in \lie (G)}\,
F (\Sigma_{R}^{\infty} EG/N_+ ,  Y  ))^{ J N'}  \rarr \]
\[  (\hocolim_{N \in \lie (G)} \,
F (\Sigma_{R}^{\infty} EG/N_+ , \ Y ))^J . \]
 Levelwise, $Y$ is a bounded above fibrant object in the
 $\overline{\lie (G)}$--model structure on $\p
\mm_R$. Cofinality (take the homotopy
colimit over $ N= N'$), using that the $n$-th skeleta of $ EG/N_+$
are all finite cell complexes (hence small objects), and finally,
by Lemma \ref{lem:restrict-n-equivalence} which says  that
$ \hocolim_{N \in \lie (G)} \, Y^{NH} \rarr Y^{H } $ is an
equivalence of pro--spectra in the Postnikov model structure,
for any closed subgroup $H$ of $G$.  Note  that the internal hom
functor from a small object  respects the
functor that sends a spectrum to an
$\Omega$--spectrum (see \ref{fibrantreplacement}).
\end{proof}

\begin{pro}
\label{pro:compositionoffixedpoints} Let $G$ be a compact Hausdorff
group and let $K $ be a closed normal subgroup of $G$.
Then there is an equivalence of pro--spectra
\[ ( Y^{h_G K} )^{hG/K} \simeq Y^{h G}  \]
in the (non-equivariant) Postnikov model structure on $\p \mm_{R}$.
\end{pro}
\begin{proof} We assume that $  Y $ is   fibrant in the
Postnikov--$\overline{\lie (G)}$--model structure on $\p \mm_R$.
The pro--spectrum  $ ( Y^{ h_G K } )^{h G/K}$ is equivalent to
\begin{equation} \label{eq:sick}
\hocolim_W \, F ( \Sigma_{R}^{\infty} EG/WK_+ , \hocolim_{ N,L} \, F (
\Sigma_{R}^{\infty} EG / N_+ ,  Y  )^{LK} )^G .
\end{equation} We use that the internal hom functor (from a small object) respects the
functor that sends a spectrum to an
$\Omega$--spectrum (\ref{fibrantreplacement}).
Since \[ \hocolim_{ N,L} \, F ( \Sigma_{R}^{\infty}
EG / N_+ ,  Y  )^{LK} \] is levelwise bounded above (since $Y$ is)
in the $ \overline{\lie (G/K)} $--model
structure on $\mm_R$, and $ \Sigma_{R}^{\infty} EG/WK_+$ has a
dualizable $n$-th skeleton for all $n$, \ref{eq:sick} is equivalent
to \[ \hocolim_{ N,L, W}\,
 F ( \Sigma_{R}^{\infty} EG/WK_+ , F (\Sigma_{R}^{\infty} EG / N_+ , Y  )^{LK}
)^G . \] (See Proposition \ref{pro:comparisonhtpfixedpoints}.)
 By cofinality, \ref{eq:sick} is equivalent to the colimit
over $N =L=W$. The fixed point adjunction and internal hom adjunction now give
\[ (\hocolim_{ N} \,
F (  \Sigma_{R}^{\infty} ( EG/KN_+  \wedge EG / N_+ ) ,
Y ))^G . \] This is equivalent to $ Y^{h G}$ by Lemma
\ref{lem:description of htpfixed points}.
\end{proof}

\subsection{Homotopy orbit and homotopy fixed point spectral
sequences} In this subsection we work in the Postnikov
$\overline{\lie (G)}$--model structure on $\p \mm_R$ for a compact
Hausdorff group $G$ and a trivial universe. We denote the homsets in the associated
homotopy category by $ [ X ,Y ]_G$.

\begin{defn} \label{htporbits}
Let $X = \{ X_a \}$ be a $G$--pro--spectrum. Then the \mdfn{$G$--homotopy
orbit pro--spectrum  $X_{hG}$} of $X$ is
\[  X_{\text{cof}}   \wedge_G \{ \Sigma^{\infty}_R E G /N_+ \} =
\{ (  (X_a)_{\text{cof}} \wedge \Sigma^{\infty}_R E G /N_+  )/G \}_{N ,a}, \] where $\{ (X_a)_{\text{cof}} \}_a $ is a cofibrant replacement of $X$.   The \mdfn{Borel
homology} of $X$ is $ \pi_* (  X_{hG}) $.
\end{defn}

\begin{lem} \label{htporbitinvariance}
Let  $G$ be a profinite group and set $\cc = \overline{\lie (G)}$.  Then the $G$--homotopy orbits
of a $G$--pro--spectrum $X$
are the $G$--orbits of a cofibrant replacement of $X$ in the $\cc$--free
model structure on $\p \mm_R$.
Let $ f \col X \rarr Y$ be a weak equivalence in the $\cc$--free model structure on
$\p \mm_R$.  Then $ f_{hG} \col X_{hG} \to Y_{hG}$ is a weak equivalence.
\end{lem}
\begin{proof}  The first claim follows from Lemma \ref{lem:freecofibrant}
since $\overline{\lie (G)}$ is a normal \good\ collection by Example \ref{exmp:lie}.
The acyclic cofibrations are maps that are essentially levelwise
$C_U \cap W_n $--maps, for all $U \in \cc$ and $n \in \zz$  \cite[5.11]{ffi}.
The $G$--orbit  functor is a Quillen left adjoint functor from
$\p G$--$\mm_R$ with the $\cc$--free model structure to $ \p \mm_R$
with the Postnikov model structure. (This follows since
$ (G/ H)/G $ is a point for all $H$.)
Hence the $G$--orbit functor respects weak equivalences between cofibrant
objects \cite[8.5.7]{hir}.
\end{proof}

The next result is
an instance of Theorem \ref{thm:AH-spectralsequence}.

\begin{pro} \label{pro:specseqEG}
Let $X$ and $Y$ be objects in $\p \mm_R$ with the Postnikov
$\overline{\lie(G)}$--model structure. Assume that $X$ is cofibrant
and essentially bounded below.
Then there is a spectral sequence
with \[ E_{2}^{p,q} = H^p (X \wedge \{\Sigma_{R}^{\infty} EG/N_+ \}
;\Pi^{ \overline{\lie (G) } }_{-q} (Y ) ) \] converging
conditionally to $ [ X \wedge \{ \Sigma_{R}^{\infty} EG/N_+ \} , Y
]^{p+q}_G$. If $Y$ is a monoid in the homotopy category of $\p
\mm_R$ with the strict $\overline{\lie (G)}$--model structure, then
the spectral sequence is multiplicative.
\end{pro}

The homotopy orbit and fixed point spectral sequences are special
cases of the spectral sequence in Proposition \ref{pro:specseqEG}.
We first consider the homotopy orbit spectral sequence.

\begin{cor} Let $X$ and $Y$ be two objects in $\p \mm_R$, and assume that
$X$ is cofibrant and bounded below, and that
$Y$ is  a non-equivariant pro--spectrum.
Then there is a spectral sequence
with \[ E_{2}^{p,q} = H^p ( (X \wedge \{\Sigma_{R}^{\infty} EG/N_+
\})/ G ;\Pi^{ 1 }_{-q} (Y) ) \] converging conditionally to $
Y^{p+q} ( X_{hG} )$. If, in addition, $Y$ is a monoid in the
homotopy category of $\p \mm_R$ with the strict $\overline{\lie
(G)}$--model structure,
then the spectral sequence is multiplicative.\end{cor}

We now consider the homotopy fixed point spectral sequence. If each
$X_s$ in $X = \{ X_s \}$ is a retract of a finite $\overline{\lie
(G)}$--cell complex, then the abutment of the spectral sequence in
Proposition \ref{pro:specseqEG} is naturally isomorphic to \[ [ X ,
\hocolim_N \, F (\Sigma_{R}^{\infty} EG/N_+ , Y ) ]^{p+q}_G \] by Lemma
\ref{innerhom}. If, in addition, $X$ is a $\p G/K$--spectrum for some
$K \in \cc$ (made into a $\p G$--spectrum), then the abutment is
isomorphic to  \[ [ X , \hocolim_N \, F
(\Sigma_{R}^{\infty} EG/N_+ , Y )^K ]^{p+q}_{G/K}.\]

\begin{pro}
\label{pro:htpfixedptss} Let $Y$ be in $\p \mm_R$. Then there is a
 spectral sequence with \[ E_{2}^{p,q} =
H^{p}_{\cont} ( G ; \Pi^{1}_{-q}(Y ) )
\] converging conditionally to $ \pi_{-p-q} ( Y^{hG}) $.  If $Y$ is a
monoid in the homotopy category of $\p \mm_R$ with the strict
$\overline{\lie (G)}$--model structure, then the spectral sequence
is multiplicative.
\end{pro}
\begin{proof} This follows from Lemma \ref{innerhom}
 and Proposition \ref{pro:specseqEG}. \end{proof}
A spectral sequence of this type  was first studied by Devinatz and Hopkins
for the spectrum $E_n$, with an action by the extended Morava
stabilizer group \cite{dho}. It has also been studied by Daniel
Davis \cite{dav2, dav, dav3, dav4}.

We combine the homotopy fixed point spectral sequence and
Proposition \ref{pro:compositionoffixedpoints}
to obtain a generalization of the Lyndon--Hochschild--Serre spectral
sequence. A spectral sequence like this was obtained by Ethan
Devinatz for $E_n$ \cite{dev2}. See also \cite[7.4]{dav4}.
\begin{pro} \label{LHS} Let $ K \unlhd L $ be closed
subgroups of $G$.
Let $Y$ be any $\p G$--$ R$--module, where  $R$ is a
$\overline{ \lie (G) }$--connective $\sss$--cell complex ring
with trivial $G$--action.  Then there is a
spectral sequence with \[ E_{2}^{p,q} = H^{p}_{\cont} (
L/K ; \Pi^{1}_{-q} (Y^{h_G K} ) ) \] converging conditionally
to $ \pi_{-p-q} ( Y^{h_G L}) $. If $Y$ is a monoid in the homotopy
category of $\p \mm_R$ with the strict  $\overline{\lie (G)}$--model
structure, then the spectral sequence is multiplicative.\end{pro}
\begin{proof}  Without loss of generality  assume that $Y$ is  both
fibrant and cofibrant in
the Postnikov $\overline{\lie (G)}$--model structure on $\p \mm_R$.
Apply Propositions \ref{pro:compositionoffixedpoints} and
\ref{pro:htpfixedptss}  to  the pro--spectrum $(Y')^{h_G K }$,
where \[Y' = \hocolim_{N \in \lie (G)} \,
F (\Sigma_{R}^{\infty} EG/N_+ ,   Y  ) \] (levelwise hocolimit).
The $L$-th homotopy fixed points of $Y'| L$  is
equivalent to $Y^{h_G L }$ by Lemma \ref{lem:comparisonhtpfix}.
The replacement $Y'$ respects monoids
 as in the proof of  Corollary \ref{cor:respectmult}.
\end{proof}

We now give a more concrete description of the $E_2 $--term of the
homotopy fixed point spectral sequence for certain pro--spectra.
\begin{pro}
Let $K$ be a closed subgroup of $G$. Let $\{Y_m \} $ be a countable
tower in $\p \mm_R$. Then there is a short exact sequence  \[ 0
\rarr \lim^{1}_{m,n} \, H^{p-1}_{\cont} ( G/K ; \Pi^{ 1 }_{-q}( (
P_n Y_m)^{h_G K } ) ) \rarr E_{2}^{p,q} \]
\[ \rarr \lim_{m,n} \, H^{p}_{\cont} ( G/K ; \Pi^{ 1 }_{-q} ((
P_n Y_m)^{h_G K } ) ) \rarr 0  \] where $E_{2}^{p,q}$ is the $E_2 $--term of
the spectral sequence in Proposition \ref{LHS}.
\end{pro}

\subsection{Comparison to Davis' homotopy fixed points}
In this subsection we show that the  homotopy fixed
points defined in  Definition \ref{defn:htpfixedpoints}   agree with the classical
homotopy fixed points when $G$ is a compact Lie group. We also
compare our definition of homotopy fixed points
to a construction by Daniel Davis \cite{dav}. We work in homotopy
categories.
The next lemma says that if $G$ is a compact Lie group, then the
orthogonal $G$--spectrum  associated to the homotopy fixed points in $\p \mm_{\sss}$,
with the strict or with the Postnikov cofree model structures are
equivalent.

\begin{lem} \label{lem:ordinaryhtpfps} Let $G$ be a compact Lie group (or a
discrete group), and let $\cc$ be the collection of all (closed)
subgroups of $G$. Let $X$ be any orthogonal $G$--spectrum. Then the map
\[  F ( \Sigma^{\infty} EG_+ ,  X)^G \rarr
\holim_n \, F ( \Sigma^{\infty} EG_+ , P_n X)^G \] is an
equivalence. \end{lem}
\begin{proof}
Since $\{1 \} \in \cc$, the pro--spectrum $ \{\Sigma^{\infty}
EG/N_+ \}$, indexed on normal subgroups $N$ in $G$, is equivalent to
$ \Sigma^{\infty} EG_+$.
 The spectrum $\Sigma^{\infty} EG_+$ is the homotopy
colimit of the ($G$--cell complex) skeleta $ \Sigma^{\infty} EG_{+}^{(m)}$, for $m
\geq 0$. Hence, $F (\Sigma^{\infty} EG_+ , Z)$ is equivalent to
\[F ( \hocolim_m \, \Sigma^{\infty}
EG_{+}^{(m)} , Z) \simeq \holim_m \, F ( \Sigma^{\infty}
 EG_{+}^{(m)} , Z)
\] for any $G$--spectrum $Z$.  The canonical map
\[ F ( \Sigma^{\infty} EG_{+}^{(m)} , X )^G \rarr F (
\Sigma^{\infty} EG_{+}^{(m)} , P_n X )^G \] is $(n - m -
\text{dim}\, G )$--connected. So
\[ \holim_m \, F ( \Sigma^{\infty} EG_{+}^{(m)} ,
X )^G \rarr \holim_{m , n} \, F ( \Sigma^{\infty} EG_{+}^{(m)} ,
P_n X )^G \] is a weak equivalence. This proves the claim.
\end{proof}

Daniel Davis defines homotopy fixed point spectra for towers of
discrete simplicial Bousfield--Friedlander $G$--spectra for any
profinite group $G$ \cite{dav}. The main difference from our
construction, translated into our terminology, is that he uses the
strict model structure on $\p \mm_{\sss}$ obtained from the
$\overline{\fin (G)}$--cofree model structure on $G \mm_{\sss}$, rather
than the Postnikov $\overline{\fin (G)}$--cofree model structure on
$G \mm_{\sss}$. He shows that if $G$ has finite virtual cohomological
dimension, then his definition of the homotopy fixed points of a
(discrete) pro--$G$--spectrum $\{ Y_b \} $ is equivalent to
\begin{equation} \label{eq:davis2} ( \holim_b \, \text{Tot} \,  \hocolim_N \,
\Gamma^{\bullet}_{G/N} ( Y_b ) )^G ,
\end{equation} where $N$ runs over all open normal subgroups of $G$.
Here $\Gamma^{\bullet}_{G/N} ( Y_b )$ is defined to be the
cosimplicial object given by $ F_{\square } (G/N^{\bullet +1 } , Y_b
)$, where $F_{\square}$ is the cotensor functor, and
$G/N^{\bullet +1 }$ is a simplicial object obtained from a group
\cite{dav2, dav3}. When $G$ has finite virtual cohomological dimension,
one can use \ref{eq:davis2} as the definition of homotopy fixed
points for categories of $G$--spectra other than the category of
discrete simplicial $G$--spectra.

Since Tot is the homotopy inverse limit of the $\text{Tot}_n$, and
$\text{Tot}_n$ is a finite (homotopy) limit, we get that
\ref{eq:davis2} is equivalent to
\[(\holim_{b,n} \,   \hocolim_N \, \text{Tot}_n \Gamma^{\bullet}_{G/N}
 ( Y_b ) )^G
. \] By definition, $\text{Tot}_n \Gamma^{\bullet}_{G/N} ( Y_b )$ is
equivalent to the internal hom $
 F(\Sigma^{\infty}  (EG/N_+)^{(n)} , Y_b )$ in $\mm_{\sss}$,
  where $(EG/N_+)^{(n)}$ denotes the $n$-th skeleton
of $EG/N_+$. Hence Davis' homotopy fixed points are equivalent to
\begin{equation} \label{eq:Davis} \holim_{b,n} \, ( \hocolim_N \, F (
\Sigma^{\infty} (EG/N_+)^{(n)} , Y_b ))^G .
\end{equation}

We compare \ref{eq:Davis} to the spectrum associated to our
definition of homotopy fixed points.

\begin{pro}
\label{pro:comparisonhtpfixedpoints} Let $G$ be a profinite group.
The canonical maps from \[\{ \hocolim_N \, F (\Sigma^{\infty}
(EG/N_+)^{(n)} , Y_b )^G \}_{b,n}\] and
\[ \{   \hocolim_N \, F (
\Sigma^{\infty} EG/N_+ , P_m Y_b )^G \}_{b, m } \] to
\begin{equation} \label{eq:F} \{ \hocolim_N \, F (
\Sigma^{\infty} (EG/N_+)^{(n)} , P_m Y_b )^G \} _{b, m ,n}
\end{equation} are both equivalences in pro--$\mm_{\sss}$ with
the strict model structure.
\end{pro}
\begin{proof} It
suffices to prove the result when $Y$ is a constant pro--spectrum.
Since $P_m Y $ is co--$m$--connected,
the skeletal inclusion gives an equivalence \[
 F ( \Sigma^{\infty} (EG/N_+) , P_m Y ) \rarr
F (\Sigma^{\infty} (EG/N_+)^{(n)} , P_m Y) \] when $ n >m$. Hence
the map from the second expression to \ref{eq:F} is an equivalence.

Since $(EG/N_+)^{(n)}$ only has cells in dimension less than or equal to
$n$ we get that \[
 F ( \Sigma^{\infty} (EG/N_+)^{(n)} ,  Y ) \rarr
F ( \Sigma^{\infty} (EG/N_+)^{(n)} , P_m Y) \] is an equivalence
when $ n < m$.
Hence the map from the first expression to \ref{eq:F} is an equivalence.
\end{proof}
Hence, the spectrum associated to our definition of the homotopy
fixed point pro-spectrum agrees with Davis' definition when $G$ is a
profinite group with finite
virtual cohomological dimension.

\begin{cor} \label{cor:respectmult}
If $Y$ is a   (commutative)  algebra in $\p \mm_{ \sss }$,
then $ Y^{h_G K}$  is a  (commutative) algebra
in $\p \mm_{ \sss }$, for all $K$.
\end{cor}
\begin{proof} By Proposition \ref{pro:comparisonhtpfixedpoints} we
get that $Y^{h_G K}$ is equivalent to \[ \{ \hocolim_U \, (\hocolim_N \, F
( \Sigma^{\infty}  (EG/N_+)^{(n)} , Y_b ))^{UK} \}_{b, n } . \] The
result follows since the pro--category is cocomplete by
\cite[11.4]{tfi},
directed colimits of algebras are created in the underlying category
of modules, and fixed points preserves algebras.
\end{proof}

\appendix
\section{Compact Hausdorff Groups}
In this appendix we recall some well known properties of compact Lie
groups. The relationship between compact Lie groups and
compact Hausdorff groups is analogues to the relationship between
finite groups and profinite groups.

We first note that if $G$ is a compact Hausdorff group, then the
finite dimensional $G$--representations are all obtained from suitable
$G/N$--representations via the quotient map $G \rarr G/N$ where
$G/N$ is a compact Lie group quotient of $G$.

\begin{lem} \label{Lierepresentation}
Let $V$ be a finite dimensional $G$--representation. Then the $G$--action on
$V$ factors through some compact Lie group quotient $G/N$ of $G$.
\end{lem}
\begin{proof}  A $G$--representation $V$ is a group homomorphism
\[\rho \col  G \rarr \text{GL} (V) . \]
The action factor through the image $\rho (G)$. Since $G$ is a
compact group $\rho (G)$, with the subspace topology from $\text{GL}
(V)$, is a closed subgroup of the Lie group $\text{GL} (V)$. Hence
$\rho (G)$ is itself a Lie group. Again, since $G$ is compact
Hausdorff, the subspace topology on $\rho (G)$ agrees with the
quotient topology from $\rho$. Hence $\rho$ gives a  homeomorphism $G /
\ker \rho \cong \rho (G)$, and  $G / \ker \rho$ is a compact Lie group.
\end{proof}

Recall from Example \ref{exmp:lie} that $ \lie (G)$ denotes the
collection of closed normal subgroups $N$ of $G$ such that $G/N$
is a compact Lie group. We consider the inverse system of quotients $G/N$
such that $G/N$ is a compact Lie group. If $G/N$ and $G/K$ are
compact Lie groups, then $G/ N\cap K$ is again a compact Lie group,
since it is a closed subgroup of $ G/N \times G/K$. Hence the
inverse system is a filtered inverse system.

In the next theorem it is essential that we work in the category of
compactly generated weak Hausdorff topological spaces.

\begin{pro} \label{continuous} Let $X$ be a topological space with
a (not necessarily continuous) $G$--action. Then the $G$--action on
$X$ is continuous if and only if the $G$--action on $X/N$ is
continuous for all subgroups $N \in \overline{\lie (G)} $ and the
canonical map \[ \rho \col X \ra \lim_N \,  X/N ,
\] where the limit is over all $N \in \lie (G)$, is a homeomorphism.
\end{pro}
\begin{proof}
Assume that $\rho$ is a homeomorphism. Then the $G$--action on $X$
is continuous since the $G$--action on $\lim_N \, X/N$ is continuous.

We now assume that the $G$--action on $X$ is continuous. We first
show that  \[ \rho \col  X \rarr \lim_N \, X/N \]
is a bijection. The Peter--Weyl theorem for compact Hausdorff groups
implies that there are enough finite dimensional real $G$--representations
to distinguish any two given elements in $G$
\cite[3.39]{ada}. Hence $ \cap_N N$ is $1 $, and $\rho$ is
injective. Now let $ \{ N x_N \} $ be an element in $\lim_N \, X/N$.
Since $G$ is a compact group and since the $G$--action on $X$ is
continuous we get that $N x_N$ is a compact subset of $X$ for every $N \in \lie (G)$.
In particular, $N x_N$ is a closed subset of the compact space $ G x_G $
for every  $N $ in $\lie (G)$. Since $\cap_N N =1$ and
$N x_N \cap L x_L \supset N\cap L x_{N\cap L}$, we conclude that the
intersection of the closed sets $N x_N $, for $N \in \lie (G)$, is a
point. Call this point $x$. We then have that $\rho
(x) = \{ N x_N \}$. So $\rho$ is surjective.

We need to show that $\rho$ is a closed map. This amounts to showing
that for any closed set $A$ of $X$, and for any $N \in \lie (G)$ we
have that $N \cdot A$ is a closed subset of $X$.
 When $A$ is a compact (hence closed) subset
of $X$ this follows since $N \cdot A$ is the image of $ N \times A$
under the continuous group action on $X$. Since we use the compactly
generated topology the subset $N \cdot A$ of $X$ is closed if for
all compact subsets $K$ of $X$ the subset $ (N \cdot A ) \cap K$ is
closed in $X$. This is true since
\[ ( N \cdot A ) \cap K = ( N \cdot ( A \cap ( N \cdot K ))) \cap
K \] and $N \cdot K$ is a compact subset of $X$. Hence $\rho $ is a
homeomorphism.
 \end{proof}

\begin{cor} \label{Neil} Any compact Hausdorff group $G$ is an
inverse limit of compact Lie groups.
\end{cor}
\begin{proof}  This follows from Theorem \ref{continuous} by
letting $X$ be $G$.
\end{proof}

The pro--category of compact Lie groups
is equivalent to the category of compact Hausdorff groups. This
follows since a closed subgroup of a compact Lie group is again a
compact Lie group. The categories are equivalent as
topological categories since both homspaces are compact Hausdorff
spaces.

Groups which are inverse limits of Lie groups have
been studied recently. See for example \cite{hom}.

\begin{cor}
The category  $G \ttt$ is a full subcategory and a retract of the pro--category
of the full subcategory of $G \ttt$ consisting of $G$--spaces
 with a $G$--action factoring through $G/N$ for some $N \in \lie (G)$.
\end{cor}
\begin{proof}
A $G$--space $X$ is sent to the pro--$G$--space $\{ X /N \}$.
The retract map is given by taking the inverse limit. By Theorem
\ref{continuous} the composite is isomorphic to the identity map on
$ G \ttt$. Let $X$ and $Y$ be two $G$--spaces.  Then  the canonical map
\[ G \ttt ( X , Y ) \rarr \lim_L \, \colim_N \, G  \ttt ( X / N , Y /
L ) \] is a bijection.
\end{proof}

\begin{rem}  \label{set} Let $G$ be a profinite group. We
observe that in the category of sets, $X$ is a continuous $G$--set
if and only if $ \colim_N \, X^N \rarr X $ is a bijection. On the
other hand, in the category of  compactly generated weak Hausdorff spaces, $X$ is a
continuous $G$--space if and only if $ X \rarr \lim_N \, X/N$ is a
continuous $G$--space. \end{rem}

\end{document}